\documentclass[preprint,3p,times]{elsarticle}

\usepackage[section]{placeins}
\usepackage{multirow}
\usepackage{makecell}
\usepackage{longtable}
\usepackage{booktabs}
\usepackage{supertabular}
\usepackage{color}
\usepackage{ulem}
\usepackage{colortbl}

\usepackage{soul}
\soulregister\cite7
\soulregister\citep7
\soulregister\citet7
\soulregister\ref7
\soulregister\pageref7

\usepackage{cancel}

\usepackage[bf]{caption2}

\newcommand{\minewB}[1]{{\color{black}{#1}}}

\newcommand{\miold}[1]{\iffalse{#1}\fi}

\usepackage{amssymb}
\usepackage{amsmath}
\newtheorem{theorem}{Theorem}
\newtheorem{lemma}{Lemma}

\newtheorem{remark}{Remark}
\newtheorem{example}{Example}

\journal{Elsevier}

\begin{document}

\begin{frontmatter}

\title{An efficient mapped WENO scheme using approximate constant 
mapping}

\author[a]{Ruo Li}
\ead{rli@math.pku.edu.cn}

\author[b,c]{Wei Zhong\corref{cor1}}
\ead{zhongwei2016@pku.edu.cn}

\cortext[cor1]{Corresponding author}

\address[a]{CAPT, LMAM and School of Mathematical Sciences, Peking
University, Beijing 100871, China}

\address[b]{School of Mathematical Sciences, Peking University,
Beijing 100871, China}

\address[c]{Northwest Institute of Nuclear Technology, Xi'an 
710024, China}

\begin{abstract}
  We present a novel mapping approach for WENO schemes through the 
  use of an approximate constant mapping function which is 
  constructed by employing an approximation of the classic signum 
  function. The new approximate constant mapping function is 
  designed to meet the overall criteria for a proper mapping 
  function required in the design of the WENO-PM6 scheme. The 
  WENO-PM6 scheme was proposed to overcome the potential loss of 
  accuracy of the WENO-M scheme which was developed to recover the 
  optimal convergence order of the WENO-JS scheme at critical 
  points. Our new mapped WENO scheme, denoted as WENO-ACM, maintains 
  almost all advantages of the WENO-PM6 scheme, including 
  low dissipation and high resolution, while decreases the number of 
  mathematical operations remarkably in every mapping process 
  leading to a significant improvement of efficiency. The 
  convergence rates of the WENO-ACM scheme have been shown through 
  one-dimensional linear advection equation with various initial 
  conditions. Numerical results of one-dimensional Euler equations 
  for the Riemann problems, the Mach 3 shock-density wave 
  interaction and the Woodward-Colella interacting blastwaves are
  improved in comparison with the results obtained by the WENO-JS,
  WENO-M and WENO-PM6 schemes. Numerical experiments with
  two-dimensional problems as the 2D Riemann problem, the 
  shock-vortex interaction, the 2D explosion problem, the double 
  Mach reflection and the forward-facing step problem modeled via 
  the two dimensional Euler equations have been conducted to 
  demonstrate the high resolution and the effectiveness of the 
  WENO-ACM scheme. The WENO-ACM scheme provides significantly better 
  resolution than the WENO-M scheme and slightly better resolution 
  than the WENO-PM6 scheme, and compared to the WENO-M and WENO-PM6 
  schemes, the extra computational cost is reduced by more than 
  $83\%$ and $93\%$, respectively.

\end{abstract}



\end{frontmatter}

\section{Introduction}
\label{secIntroduction}
Essentially non-oscillatory (ENO) schemes \cite{ENO1986, ENO1987V24,
ENO1987JCP71, ENO1987JCP83, ENO-Shu1988, ENO-Shu1989} and weighted
ENO (WENO) schemes \cite{WENO-LiuXD, WENO-JS} have been developed 
quite successfully to solve the hyperbolic conservation laws, which 
may develop discontinuities in its solutions even if the initial 
conditions are smooth. The goal of this paper is to propose an 
improved version of the fifth-order WENO scheme for the hyperbolic 
conservation laws taking the form
\begin{equation*}
\begin{array}{ll}
\mathbf{u}_{t} + \displaystyle\sum\limits_{\alpha = 1}^{d}
\mathbf{f}_{\alpha}(\mathbf{u})_{x_{\alpha}} = 0, 
& x_{\alpha} \in \mathbb{R}, t > 0.
\end{array}
\label{governingEquation}
\end{equation*}
Here the function $\mathbf{u} = (u_{1}, u_{2}, \cdots, u_{m})^{
\mathrm{T}}$ is an $m$-dimensional vector of conserved variables,
and flux $\mathbf{f}_{\alpha}(\mathbf{u})$ is a vector-valued 
function of $m$ components with $x_{\alpha}$ and $t$ variables. 

Liu et al. \cite{WENO-LiuXD} developed the first version of WENO 
schemes which convert an $r$th-order ENO scheme \cite{ENO1986,
ENO1987V24,ENO1987JCP71, ENO1987JCP83, ENO-Shu1988, ENO-Shu1989} 
into an $(r+1)$ th-order WENO scheme by using a convex combination 
of all candidate substencils instead of just one as in the original 
ENO scheme. Later, Jiang and Shu \cite{WENO-JS} proposed the 
classic WENO-JS schemes with an improvement that an $r$th-order 
ENO scheme can be converted into a $(2r-1)$th-order WENO scheme by
introducing a new definition of the smoothness indicator used to 
measure the smoothness of the numerical solution on a substencil. 
Then, the weighting method presented in \cite{WENO-LiuXD} and the 
smoothness indicators designed in \cite{WENO-JS} eventually became a 
standard, and the WENO-JS schemes especially the fifth-order one 
\cite{WENO-JS} developed into one of the most popular high-order 
methods \cite{WENO-PPM5}. In recent decades, many successful works 
have been done to raise some issues about WENO schemes 
\cite{WENO-M, WENO-IM, WENO-PM, WENO-Z, WENO-Z01, WENO-Z02, 
WENO-Accuracy, WENO-GlobalAccuracy,WENO-ZS}.

It was clearly pointed out by Henrick et al. \cite{WENO-M} that,
in general, the
fifth-order WENO-JS scheme is only third-order or even less accurate
at critical points of order $n_{\mathrm{cp}} = 1$ in smooth regions, 
where $n_{\mathrm{cp}}$ denotes the order of the critical point; 
e.g., $n_{\mathrm{cp}} = 1$ corresponds to $f'=0, f'' \neq 0$ and 
$n_{\mathrm{cp}} = 2$ corresponds to $f'=0, f''=0, f''' \neq 0$, 
etc. To overcome this problem, Henrick et al. \cite{WENO-M}
introduced a carefully designed mapping function leading to the first
mapped WENO scheme named WENO-M. Compared to the WENO-JS scheme 
\cite{WENO-JS}, the WENO-M scheme is able to recover the optimal 
convergence order near critical points in smooth regions and 
generate more accurate solutions. Another significant contribution of
the work by Henrick et al. \cite{WENO-M} is that they derived a
strong sufficient condition on the weights of substencils for WENO
schemes to achieve optimal convergence orders and this condition has 
become the primary criterion in the design of all other mapped WENO 
schemes \cite{WENO-PM,WENO-IM,WENO-PPM5,WENO-RM260,WENO-RM-Vevek2018,
WENO-AIM,WENO-MAIMi}. Recently, Feng et al. \cite{WENO-PM} found
that the mapping function of the WENO-M scheme \cite{WENO-M} may 
amplify the effect from the non-smooth substencils and thus cause a 
potential loss of accuracy near discontinuities. To address this
issue, they proposed two additional requirements to the mapping 
function $g(\omega)$, that is, $g'(0) = 0$ and $g'(1) = 0$. Then, 
the WENO-PM$k$ scheme has been devised \cite{WENO-PM} by employing a 
piecewise polynomial mapping function which satisfies the two 
additional requirements above and the original criteria in 
\cite{WENO-M}, and the WENO-PM$6$ scheme was recommended, which is 
able to generate more accurate numerical solutions \cite{WENO-PM}
near the discontinuities than the classic WENO-JS scheme and the 
WENO-M scheme. Furthermore, the two additional requirements were 
considered to be very important to decrease the effect from the 
non-smooth substencils \cite{WENO-RM260} and were used in the 
construction of the WENO-RM($mn0$) scheme \cite{WENO-RM260}. 
Similarly, requirements $g'(0)=1$ and $g'(1)=1$, which are also used 
to decrease the effect from the non-smooth substencils, were employed
when the WENO-PPM$n$ ($n=4,5,6$) \cite{WENO-PPM5}, WENO-RM($k,m,s$) 
\cite{WENO-RM-Vevek2018}, WENO-AIM($k,m,c$) \cite{WENO-AIM} and 
WENO-MAIM$i$ \cite{WENO-MAIMi} schemes were constructed. From a 
different perspective, Borges et al. \cite{WENO-Z} proposed another 
version of the fifth-order WENO scheme, which is called the WENO-Z 
scheme. It drives the weights to the optimal values faster than the 
WENO-M scheme by employing a global higher order reference value for 
the smoothness indicators \cite{WENO-eta}. Later, the WENO-Z scheme 
was successfully extended to higher orders by Castro et al. 
\cite{WENO-Z01}. The recommended fifth-order WENO-Z scheme with the 
power parameter $p = 1$ gives less dissipation than the WENO-JS 
scheme but its convergence order is fourth-order at the first-order 
critical points \cite{article_WENO-MS_JSC2019} and will drop to 
second-order at higher order critical points \cite{WENO-eta}. 
Recently, Don et al. \cite{WENO-Accuracy,WENO-Z-optimalAccuracy} 
demonstrated that the optimal order of accuracy for the WENO-Z 
scheme can be obtained regardless of any order critical points by 
setting the parameter $\epsilon$, which is used to prevent the 
denominator becoming zero, as a function of the spacial step 
$\Delta x$. However, it was noted \cite{WENO-PDM} that adjusting 
$\epsilon$ to recover accuracy is actually an implicit switch 
between the WENO scheme and the upstream central scheme, and its 
effect is dependent on the specific problem.

Various kinds of existing mapped WENO schemes, e.g., WENO-M, 
WENO-PM6, WENO-RM260, WENO-PPM5, WENO-RM($k,m,s$), WENO-AIM
($k,m,c$), WENO-MAIM$i$, et al., can improve the performances of the 
classic WENO-JS scheme in some ways like achieving optimal 
convergence orders near critical points in smooth regions, having 
low dissipation and generating more accurate solutions near the 
discontinuities. However, as some specified complicated mapping 
procedures must be performed, the chief drawback of these existing
mapped schemes is that the computational cost increases 
significantly. Taking the WENO-M scheme that is one of the 
concerned WENO schemes in this paper as an example, as reported in
\cite{WENO-Z}, its extra computational cost is $20\%$ to $30\%$ 
compared to the WENO-JS scheme when calculating two dimensional 
Euler equations. For the WENO-PM6 scheme, another concerned WENO
scheme in this paper, it was noted that \cite{WENO-AIM} the 
piecewise nature requires logic operations to be implemented during 
the mapping process, which makes it harder to vectorize the
operations and increases the computational cost. Hong et al. 
indicated that \cite{WENO-PDM}, as the form of the mapping function 
designed in the WENO-PM6 scheme becomes more complicated than the 
WENO-M scheme, its extra computational cost will have a further 
increase. Actually, our extensive numerical tests show that the 
extra computational cost of the WENO-PM6 scheme will increase by 
more than $54\%$ compared to the WENO-JS scheme, and we will give 
the results in subsection \ref{subsec:CostTest} in detail.

In this paper, an approximate constant mapping function is designed 
at first. This new mapping function satisfies the original criteria
proposed in \cite{WENO-M} so that the new corresponding WENO scheme,
abbreviated as WENO-ACM, is able to achieve the optimal convergence
orders near critical points in smooth regions. And also, the new
mapping function maintains $g'(0) = 0$ and $g'(1) = 0$ so that it can
decrease the effect from the non-smooth substencils as the WENO-PM6
scheme does. Thus, the WENO-ACM scheme is able to yield low
dissipation and high resolution results comparable to that of the
WENO-PM6 scheme. The greatest benefit is that the new mapping 
function almost uses only one assignment operation to implement the 
mapping process, instead of evaluating the mapping functions 
involving multiple multiplication and division as other existing 
mapped WENO schemes do. Therefore, the cost of the WENO-ACM scheme 
is very low. Numerical experiments with various benchmark problems 
modeled via the two dimensional Euler equations are conducted to 
demonstrate that the WENO-ACM scheme generates significantly better 
resolution than the WENO-M scheme and slightly better resolution 
than the WENO-PM6 scheme, while the extra computational cost is 
reduced by more than $83\%$ compared to the WENO-M scheme and 
reduced by more than $93\%$ compared to the WENO-PM6 scheme.

The remainder of this paper is organized as follows. In Section 2, we
give a brief description of the finite volume methodology and the
procedures of the WENO-JS \cite{WENO-JS}, WENO-M \cite{WENO-M} and
WENO-PM$k$ \cite{WENO-PM} schemes to clarify our major concern. In
Section 3, we introduce the details on how we construct the new 
mapped WENO scheme with approximate constant mapping and then provide
the parametric study of the new mapping function and the convergence
property of the new mapped WENO scheme. In Section 4, some numerical
experiments are presented to compare the performances of different
WENO schemes, and the computational cost comparisons are also shown 
in this section. Finally, the conclusions are given in Section 5.


\section{Review of finite volume WENO schemes}
\label{secMappedWENO}

\subsection{Finite volume methodology}
Consider the following one-dimensional scalar hyperbolic 
conservation law
\begin{equation}
\partial_{t} u + \partial_{x} f(u) = 0,
\label{eq:1D-hyperbolicLaw}
\end{equation}
which is to be solved on the domain $x \in [x_{l}, x_{r}]$ for 
$t \geq 0$ with the initial condition $u(x,0) = u_{0}(x)$. 
Throughout this paper, we assume that the computational domain is 
discretized into uniform cells $I_{j} = [x_{j-1/2}, x_{j+1/2}], 
j = 1,\cdots,N$ with width $\Delta x = (x_{r} - x_{l}) / N$. The 
cell center of $I_{j}$ is denoted by $x_{j} = x_{l} + (j - 1/2)
\Delta x$ and its cell boundaries are denoted by $x_{j \pm 1/2} = 
x_{j} \pm \Delta x/2$. The cell average $\bar{u}_{j}$ of $I_{j}$ is
defined by
\begin{equation}
\bar{u}_{j} = \dfrac{1}{\Delta x}\int_{x_{j-1/2}}^{x_{j+1/2}}u(\xi,t)
\mathrm{d}\xi.
\label{eq:definition:cellAverage}
\end{equation}
By integrating Eq.(\ref{eq:1D-hyperbolicLaw}) over $I_{j}$ and
employing some simple mathematical manipulations, we can approximate
Eq.(\ref{eq:1D-hyperbolicLaw}) by the following finite volume
conservative formulation
\begin{equation}
\dfrac{\mathrm{d}\bar{u}_{j}(t)}{\mathrm{d}t} \approx -\dfrac{1}
{\Delta x}\bigg( \hat{f}(u_{j+1/2}^{-},u_{j+1/2}^{+}) - 
\hat{f}(u_{j-1/2}^{-},u_{j-1/2}^{+}) \bigg).
\label{eq:discretizedFunction}
\end{equation}
In Eq.(\ref{eq:discretizedFunction}), $\bar{u}_{j}(t)$ is the
numerical approximation to $\bar{u}_{j}$ defined in 
Eq.(\ref{eq:definition:cellAverage}), and the numerical flux
$\hat{f}(u^{-},u^{+})$ where $u^{-}$ and $u^{+}$ refer to the
left-sided and right-sided limits of $u$ is a replacement of the
physical flux function $f(u)$. For the possible presence of 
discontinuities, $u^{-}$ and $u^{+}$ are usually not equal. For 
hyperbolic laws, the numerical flux $\hat{f}(u^{-},u^{+})$ is a 
monotone function and it is consistent with the physical flux, i.e.
$\hat{f}(u,u) = f(u)$. In this paper, the global Lax-Friedrichs flux
$\hat{f}(a,b) = \frac{1}{2}[f(a) + f(b) -\alpha(b - a)]$ is chosen, 
where $\alpha = \max_{u} \lvert f'(u) \rvert$ is a constant and the
maximum is taken over the whole range of $u$. In 
Eq.(\ref{eq:discretizedFunction}), $u_{j\pm1/2}^{\pm}$ can be 
computed by the technique of reconstruction, like some WENO
reconstructions which are described in the following subsections. For
the hyperbolic conservation laws system, a local characteristic 
decomposition is used in the reconstruction, and \cite{WENO-JS} is
referred to for more details. Two commonly used classes of finite
volume WENO schemes in two dimensional Cartesian meshes are studied
in detail in \cite{FVMaccuracyProofs03}, and the one denoted as
class A is taken in this paper.

\subsection{WENO-JS}
We recall the reconstruction process of the fifth-order WENO-JS
scheme \cite{WENO-JS}, which has successfully been extended to 
higher order ones \cite{WENO-Balsara,veryHighOrderWENO}. We describe 
only the procedure of the left-biased reconstruction 
$u_{j + 1/2}^{-}$ as the right-biased one $u_{j + 1/2}^{+}$ can 
easily be obtained by mirror symmetry with respect to the location 
$x_{j + 1/2}$ of that for $u_{j + 1/2}^{-}$. For simplicity of 
notation, we do not use the subscript ``-'' in the following content.

Explicitly, for the five-point stencil $S^{5}=\{I_{j-2}, I_{j-1}, 
I_{j}, I_{j+1}, I_{j+2}\}$, the third-order approximations of 
$u(x_{j+1/2},t)$ associated with three left-biased substencils $S_{s}
=\{I_{j+s-2}, I_{j+s-1}, I_{j+s}\}$, $s = 0,1,2$ are as follows
\begin{equation}
\begin{array}{l}
\begin{aligned}
&u_{j+1/2}^{0} = \dfrac{1}{6}(2\bar{u}_{j-2} - 7\bar{u}_{j-1}
+ 11\bar{u}_{j}), \\
&u_{j+1/2}^{1} = \dfrac{1}{6}(-\bar{u}_{j-1} + 5\bar{u}_{j}
+ 2\bar{u}_{j+1}), \\
&u_{j+1/2}^{2} = \dfrac{1}{6}(2\bar{u}_{j} + 5\bar{u}_{j+1}
- 2\bar{u}_{j+2}).
\end{aligned}
\end{array}
\label{eq:approx_ENO}
\end{equation}
The fifth-order approximation of global stencil $S^{5}$ is built via
the following convex combination of the three third-order
approximations in Eq.(\ref{eq:approx_ENO})
\begin{equation*}
u_{j + 1/2} = \sum\limits_{s = 0}^{2}\omega_{s}u_{j + 1/2}^{s},
\end{equation*}
where $\omega_{s}$ are nonlinear weights. In the classic WENO-JS
scheme, the nonlinear weights are calculated by
\begin{equation} 
\omega_{s}^{\mathrm{JS}} = \dfrac{\alpha_{s}^{\mathrm{JS}}}{\sum_{l =
 0}^{2} \alpha_{l}^{\mathrm{JS}}}, \alpha_{s}^{\mathrm{JS}} = \dfrac{
 d_{s}}{(\epsilon + \beta_{s})^{2}},
\label{eq:weights:WENO-JS}
\end{equation} 
where $\epsilon$ is a small positive number introduced to prevent
the denominator being zero and it was taken to be 
$10^{-6}$ in the original WENO-JS scheme, and $d_{0}=0.1, d_{1}=0.6,
d_{2} = 0.3$ are ideal weights of $\omega_{s}$ satisfying
$\sum\limits_{s = 0}^{2} d_{s} u^{s}_{j + 1/2} =
u(x_{j + 1/2}, t) + O(\Delta x^{5})$ in smooth regions. The
parameters $\beta_{s}$ named smoothness indicators are defined as 
follows \cite{WENO-JS}
\begin{equation*}
\begin{array}{l}
\begin{aligned}
\beta_{0} &= \dfrac{13}{12}\big(\bar{u}_{j - 2} - 2\bar{u}_{j - 1} + 
\bar{u}_{j} \big)^{2} + \dfrac{1}{4}\big( \bar{u}_{j - 2} - 4\bar{u}_
{j - 1} + 3\bar{u}_{j} \big)^{2}, \\
\beta_{1} &= \dfrac{13}{12}\big(\bar{u}_{j - 1} - 2\bar{u}_{j} + \bar
{u}_{j + 1} \big)^{2} + \dfrac{1}{4}\big( \bar{u}_{j - 1} - \bar{u}_{
j + 1} \big)^{2}, \\
\beta_{2} &= \dfrac{13}{12}\big(\bar{u}_{j} - 2\bar{u}_{j + 1} + \bar
{u}_{j + 2} \big)^{2} + \dfrac{1}{4}\big( 3\bar{u}_{j} - 4\bar{u}_{j 
+ 1} + \bar{u}_{j + 2} \big)^{2}.
\end{aligned}
\end{array}
\end{equation*}

In smooth regions without critical points, the classic WENO-JS 
scheme is able to achieve fifth-order of accuracy. However, at 
critical points where the first derivative vanishes 
but the third derivative does not simultaneously, it loses accuracy 
and its order of accuracy decreases to third-order or even less. 
More details can be found in \cite{WENO-M}.

\subsection{WENO-M}
\label{subsecWENO-M}
In order to overcome the problem that the WENO-JS scheme loses 
accuracy at critical points, Henrick et al. \cite{WENO-M} proposed 
the WENO-M scheme by constructing a mapping function of the 
nonlinear weights $\omega$ given by
\begin{equation}
\big( g^{\mathrm{M}} \big)_{s}(\omega) = \dfrac{ \omega \big( d_{s} +
d_{s}^2 - 3d_{s}\omega + \omega^{2} \big) }{ d_{s}^{2} + (1 - 2d_
{s})\omega }, \quad \quad \omega \in [0, 1],
\label{mappingFunctionWENO-M}
\end{equation}
Clearly, $\big( g^{\mathrm{M}}\big)_{s}(\omega)$ is a monotonically
increasing function in $[0, 1]$ with finite slopes, and it satisfies
the following properties.
\begin{lemma} 
The mapping function defined by Eq.(\ref{mappingFunctionWENO-M})
satisfies: \\

C1. $0 \leq \big( g^{\mathrm{M}} \big)_{s}(\omega) \leq 1, \big( g^{
\mathrm{M}} \big)_{s}(0) = 0, \big( g^{\mathrm{M}} \big)_{s}(d_{s}) =
d_{s}, \big( g^{\mathrm{M}} \big)_{s}(1) = 1$;

C2. $\big( g^{\mathrm{M}} \big)_{s}'(d_{s}) = \big(g^{\mathrm{M}}
\big)_{s}''(d_{s}) = 0$. 
\label{lemmaWENO-Mproperties}
\end{lemma}

By employing Eq.(\ref{eq:weights:WENO-JS}) and Eq.(\ref
{mappingFunctionWENO-M}), one can obtain the mapped weights as 
follows
\begin{equation*}
\omega_{s}^{\mathrm{M}} = \dfrac{\alpha _{s}^{\mathrm{M}}}{\sum_{l = 
0}^{2} \alpha _{l}^{\mathrm{M}}}, \alpha_{s}^{\mathrm{M}} = \big( g^{
\mathrm{M}} \big)_{s}(\omega^{\mathrm{JS}}_{s}),
\end{equation*}

It has been analyzed and proved in detail in \cite{WENO-M} that the
WENO-M scheme is able to achieve the optimal order of accuracy in
smooth regions even near the first-order critical point.

\subsection{WENO-PM6}
\label{subsecWENO-PM}
Feng et al. \cite{WENO-PM} found that the mapping operation of the 
WENO-M scheme will cause the potential loss of accuracy near the 
discontinuities or the parts with sharp gradients. To overcome 
this drawback, they add two requirements, that is, $g'_{s}(0) = 0$ 
and $g'_{s}(1) = 0$, to the original criteria (see Lemma 
\ref{lemmaWENO-Mproperties}) by Henrick et al. \cite{WENO-M}. And
then a new mapping function is defined by the following piecewise 
polynomial function
\begin{equation}
\big( g^{\mathrm{PM}} \big)_{s}(\omega) = C_{1}(\omega - d_{s})^{k+1}
(\omega + C_{2}) + d_{s}, \quad \quad k \geq 2,
\label{mappingFunctionWENO-PM}
\end{equation}
where $C_{1},C_{2}$ are constants with specified parameters $k$ and 
$d_{s}$, and they are calculated by
\begin{equation*}
\left\{
\begin{array}{ll}
\begin{aligned}
&C_{1}=(-1)^{k}\dfrac{k+1}{d_{s}^{k+1}}, C_{2}=\dfrac{d_{s}}{k+1}, 
& \text{if} \quad 0 \leq \omega \leq d_{s}, \\
&C_{1}=-\dfrac{k+1}{(1-d_{s})^{k+1}}, C_{2}=\dfrac{d_{s} - (k + 2)}
{k+1}, & \text{if} \quad d_{s} < \omega \leq 1.
\end{aligned}
\end{array}
\right.
\end{equation*}

\begin{lemma} 
The mapping function $\big( g^{\mathrm{PM}} \big)_{s}(\omega)$
defined by Eq.(\ref{mappingFunctionWENO-PM}) satisfies: \\

C1. $\big( g^{\mathrm{PM}} \big)'_{s}(\omega) \geq 0, \omega \in [0, 
1]$;

C2. $\big( g^{\mathrm{PM}} \big)_{s}(0) = 0, \big( g^{\mathrm{PM}} 
\big)_{s}(d_{s}) = d_{s}, \big( g^{\mathrm{PM}} \big)_{s}(1) = 1$;

C3. $\big( g^{\mathrm{PM}} \big)'_{s}(d_{s}) = \cdots = \big( g^{
\mathrm{PM}} \big)^{(k)}_{s}(d_{s}) = 0$;

C4. $\big( g^{\mathrm{PM}} \big)'_{s}(0) = \big( g^{\mathrm{PM}} 
\big)'_{s}(1) = 0$.
\label{lemmaWENO-PMproperties}
\end{lemma}

Similarly, by employing Eq.(\ref{eq:weights:WENO-JS}) and 
Eq.(\ref{mappingFunctionWENO-PM}) where the parameter $k$ is taken 
to be $6$ as recommended in \cite{WENO-PM}, one can obtain the 
mapped weights of the WENO-PM6 scheme as follows
\begin{equation*}
\omega_{s}^{\mathrm{PM}6}=\dfrac{\alpha _{s}^{\mathrm{PM}6}}{\sum_{l 
= 0}^{2} \alpha _{l}^{\mathrm{PM}6}}, \alpha_{s}^{\mathrm{PM}6} = 
\big(g^{\mathrm{PM}6} \big)_{s}(\omega^{\mathrm{JS}}_{s}).
\end{equation*}

It has been verified that the WENO-PM6 scheme is able to achieve the
optimal order of accuracy as the WENO-M scheme does at critical
points. In addition, the resolution of the WENO-PM6 scheme is 
significantly higher than the WENO-JS scheme and the WENO-M scheme, 
especially for a long output time. One can see \cite{WENO-PM} for 
more details.

\subsection{Time discretization}
\label{subsecTime-discretization}
Following the method of lines (MOL) approach, the Partial 
Differential Equation (PDE) Eq.(\ref{eq:1D-hyperbolicLaw}) can be 
turned into an Ordinary Differential Equation (ODE) system of the 
form
\begin{equation}
\dfrac{\mathrm{d}\bar{u}_{j}(t)}{\mathrm{d} t} = \mathcal{L}(u_{j}),
\label{ODEs}
\end{equation}
where
\begin{equation*}
\mathcal{L}(u_{j}):=-\dfrac{1}{\Delta x}\bigg( \hat{f}(u_{j+1/2}^{-},
u_{j+1/2}^{+}) - \hat{f}(u_{j-1/2}^{-},u_{j-1/2}^{+}) \bigg). 
\end{equation*}
Then, the WENO schemes can be applied to obtain $\mathcal{L}(u_{j})$
. 

In all the numerical experiments in this paper, the ODE system 
Eq.(\ref{ODEs}) is solved using the following explicit, third-order, 
Strong Stability Preserving (SSP) Runge-Kutta method 
\cite{ENO-Shu1988,SSPRK1998,SSPRK2001}
\begin{equation*}
\begin{array}{l}
\begin{aligned}
&u^{(1)} = u^{n} + \Delta t \mathcal{L}(u^{n}), \\
&u^{(2)} = \dfrac{3}{4} u^{n} + \dfrac{1}{4} u^{(1)} + \dfrac{1}{4} 
\Delta t \mathcal{L}(u^{(1)}), \\
&u^{n + 1} = \dfrac{1}{3} u^{n} +\dfrac{2}{3}u^{(2)}
+ \dfrac{2}{3} \Delta t \mathcal{L}(u^{(2)}),
\end{aligned}
\end{array}
\end{equation*}
where $u^{(1)}$ and $u^{(2)}$ are the intermediate stages, $u^{n}$ 
is the value of $u$ at time level $t^{n} = n\Delta t$, and $\Delta t$
is the time step satisfying some proper CFL condition.


\section{The new mapped WENO scheme with approximate constant mapping
}\label{secWENO-ACM}
\subsection{Design and properties of the approximate constant mapping
function}
\subsubsection{The new mapping function}
Let $\mathrm{sgm}(x, \delta, A)$ denote the signum-like function, 
taking the form
\begin{equation*}
\mathrm{sgm}\big( x, \delta, A \big) = \left\{ 
\begin{array}{ll}
\begin{aligned}
&\dfrac{x}{|x|}, & |x| \geq \delta, \\ 
&\dfrac{x}{\Big(A\big( \delta ^2 - x^2 \big)\Big)^{k + 3} + \ |x|}, 
& |x| < \delta,
\end{aligned}
\end{array} \right.
\label{sgFunction}
\end{equation*} 
where $k \in \mathbb{N}^{+}$, $\delta > 0$ and $\delta 
\rightarrow 0$. The positive parameter $A$ is a scale 
transformation factor introduced to adjust the shape of the mapping 
function below. As mentioned in \cite{WENO-MAIMi}, we can easily 
verify that $\mathrm{sgm}\big( x, \delta, A \big)$ is monotone 
increasing. 

Then we can construct a global monotonically increasing mapping 
function, denoted as $\big( g^{\mathrm{ACM}} \big)_{s}(\omega)$, by 
directly splicing two signum-like functions as
\begin{equation}
\big( g^{\mathrm{ACM}} \big)_{s}(\omega) = \left\{
\begin{array}{ll}
\begin{aligned}
&\dfrac{d_{s}}{2}\mathrm{sgm}(\omega - \mathrm{CFS}_{s}, \delta_{s}, 
A) + \dfrac{d_{s}}{2}, & \omega \leq d_{s}, \\
&\dfrac{1-d_{s}}{2}\mathrm{sgm}(\omega -\overline{\mathrm{CFS}}_{s}, 
\delta_{s}, A) + \dfrac{1 + d_{s}}{2}, & \omega > d_{s},
\end{aligned}
\end{array}
\right.
\label{eq:mappingFunctionACM}
\end{equation}
where the Control Factor of Smoothness $\mathrm{CFS}_{s}$ is the
same as that in \cite{WENO-MAIMi} satisfying $\mathrm{CFS}_{s} \in
(0,d_{s})$, and $\overline{\mathrm{CFS}}_{s}=1- \frac{1-d_{s}}{d_{s}}
\times \mathrm{CFS}_{s}$ with $\overline{\mathrm{CFS}}_{s} \in 
(d_{s}, 1)$. In addition, the splicing condition $\mathrm{CFS}_{s} + 
\delta_{s} < d_{s} < \overline{\mathrm{CFS}}_{s} - \delta_{s}$ and 
the requirements $\mathrm{CFS}_{s} - \delta_{s} > 0$ and 
$\overline{\mathrm{CFS}}_{s} + \delta_{s} < 1$ need to be satisfied. 
Therefore, the value of the parameter $\delta_{s}$ is limited by 
$\delta_{s} < \min\Big\{\mathrm{CFS}_{s}, d_{s} - \mathrm{CFS}_{s}, 
(1 - d_{s})\Big(1 - \frac{\mathrm{CFS}_{s}}{d_{s}}\Big)
 , \frac{1-d_{s}}{d_{s}}\mathrm{CFS}_{s}\Big\}$. The effects of 
 parameters $\mathrm{CFS}_{s}, k, A$ and $\delta_{s}$ on the mapping 
 function $\big(g^{\mathrm{ACM}}\big)_{s}(\omega)$ will be discussed 
 in the following subsection.

\begin{remark}
The splicing condition $\mathrm{CFS}_{s} + \delta_{s} < d_{s} < \overline{\mathrm{CFS}}_{s} - \delta_{s}$ is used to guarantee $\big( g^{\mathrm{ACM}}\big)_{s}(d_{s}) = d_{s}$ and $\big(g^{\mathrm{ACM}}\big)'_{s}(d_{s}) = \big(g^{\mathrm{ACM}}
\big)''_{s}(d_{s}) = \cdots = 0$. Similarly, the requirements $\mathrm{CFS}_{s} - \delta_{s} > 0$ and $\overline{\mathrm{CFS}}_{s} + \delta_{s} < 1$ are imposed to ensure that $\big( g^{\mathrm{ACM}}\big)_{s}(\omega)$ satisfies the properties at the boundaries $\omega = 0$ and $\omega = 1$, that is, $\big(g^{\mathrm{ACM}}\big)_{s}(0) = \big(g^{\mathrm{ACM}}\big)'_{s}(0^{+}) = \big(g^{\mathrm{ACM}}\big)'_{s}(1^{-}) = 0$ and $\big(g^{\mathrm{ACM}}\big)_{s}(1) = 1$.
\label{rmk:splicing_condition}
\end{remark}	

As a summary, we state the trivial theorem without proof in the 
following.

\begin{theorem}
With appropriate parameters $\mathrm{CFS}_{s}$, $A$ and $\delta_{s}$,
the mapping function $\big(g^{\mathrm{ACM}}\big)_{s}(\omega)$ 
defined by Eq.(\ref{eq:mappingFunctionACM}) satisfies the following 
properties: \\

C1. $\big(g^{\mathrm{ACM}}\big)_{s}(0)=0$, $\big( 
g^{\mathrm{ACM}}\big)_{s}(d_{s}) = d_{s}$, $\big(g^{\mathrm{ACM}}
\big)_{s}(1) = 1$;

C2. $\big(g^{\mathrm{ACM}}\big)'_{s}(\omega) \geq 0$, $\omega \in 
(0,1)$;

C3. $\big(g^{\mathrm{ACM}}\big)'_{s}(d_{s}) = \big(g^{\mathrm{ACM}}
\big)''_{s}(d_{s}) = \cdots = 0$;

C4. $\big(g^{\mathrm{ACM}}\big)'_{s}(0^{+}) = \big(g^{\mathrm{ACM}}
\big)'_{s}(1^{-}) = 0$.  
\label{theorem:mappingFunction:WENO-ACM}
\end{theorem}

We can observe the properties in Theorem 
\ref{theorem:mappingFunction:WENO-ACM} intuitively from 
Fig.\ref{fig:parametricStudy}. Now, we give the 
approximate-constant-mapped WENO scheme, denoted as WENO-ACM, with 
the mapped weights
\begin{equation}
\omega_{s}^{\mathrm{ACM}} = \dfrac{\alpha _{s}^{\mathrm{ACM}}}{
\sum_{l = 0}^{2} \alpha _{l}^{\mathrm{ACM}}}, \alpha_{s}^{\mathrm{
ACM}}=\big( g^{\mathrm{ACM}} \big)_{s}(\omega^{\mathrm{JS}}_{s}).
\label{eq:nonlinearWeightsACM}
\end{equation}

\begin{remark}
By submitting the nonlinear weights $\omega_{s}^{\mathrm{JS}}$ of 
the $(2r-1)$th-order WENO-JS scheme (in \cite{WENO-Balsara} and
\cite{veryHighOrderWENO} the results of $r = 2,\cdots,9$ were
given) into Eq.(\ref{eq:nonlinearWeightsACM}), the above fifth-order
WENO-ACM scheme can be easily extended to $(2r-1)$th-order ones.
\end{remark}

\subsubsection{Parametric study}\label{subsec:par_study}
\begin{figure}[ht]
\centering
\includegraphics[height=0.39\textwidth]
{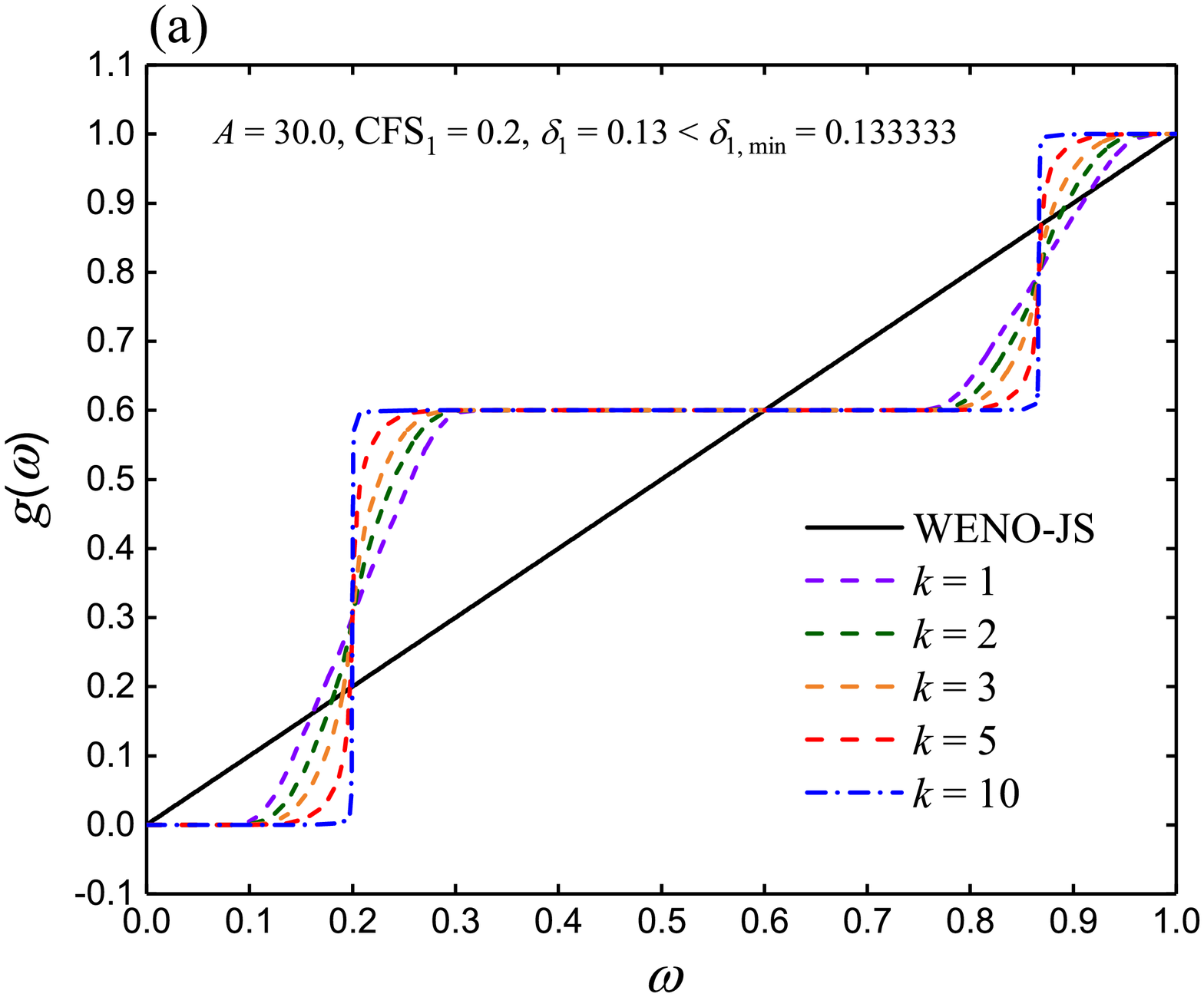}\quad
\includegraphics[height=0.39\textwidth]
{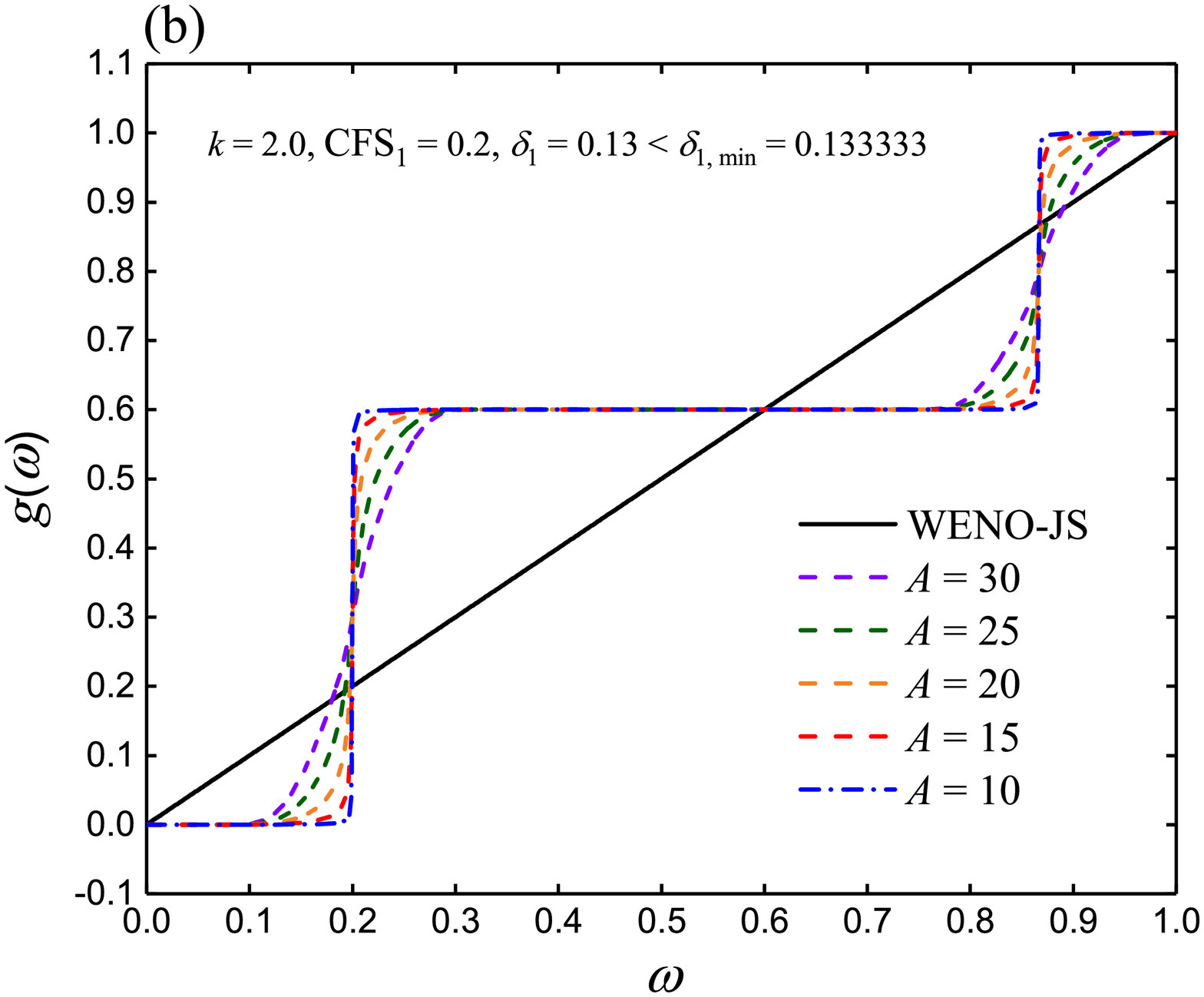} \\
\includegraphics[height=0.39\textwidth]
{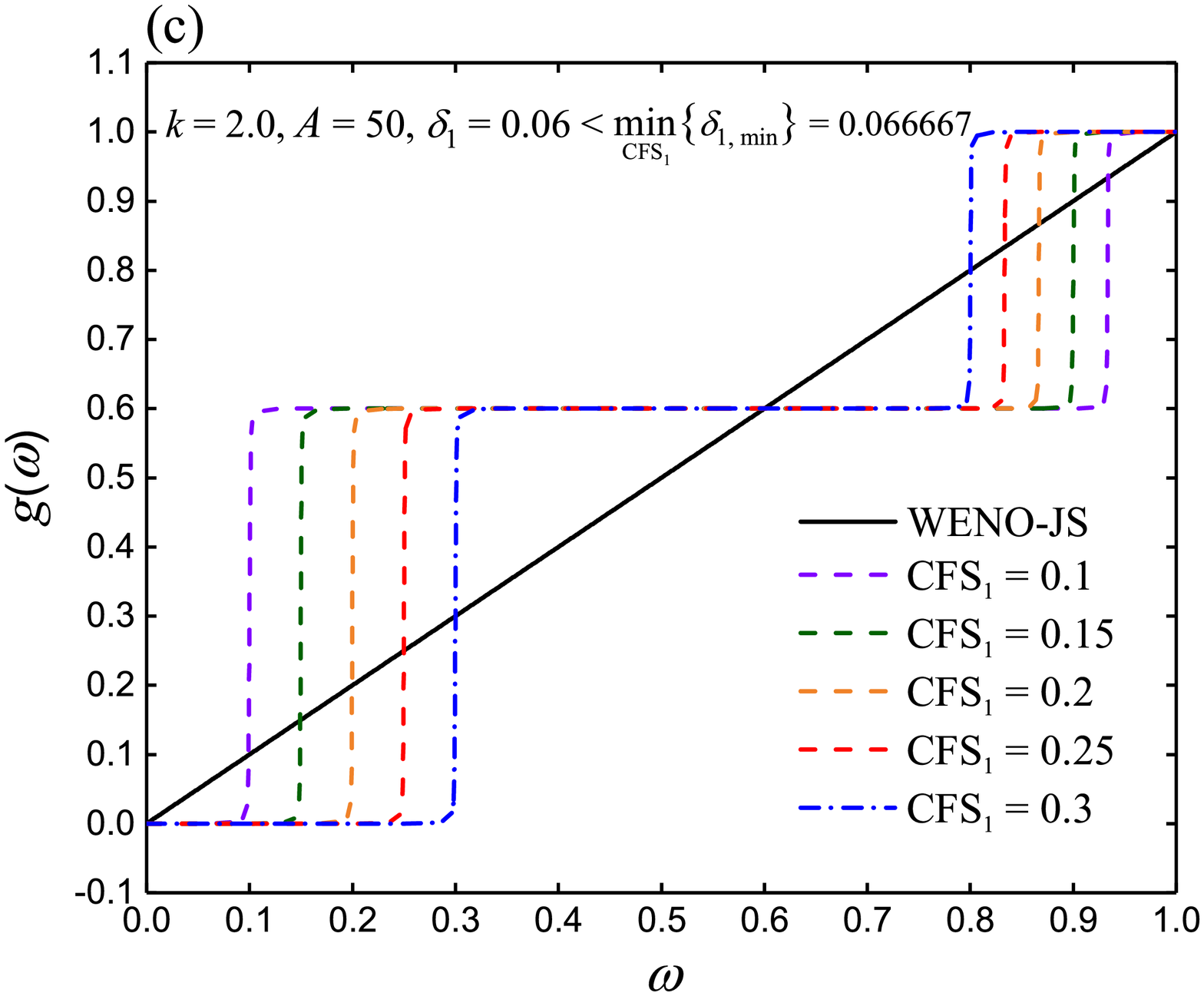}\quad
\includegraphics[height=0.39\textwidth]
{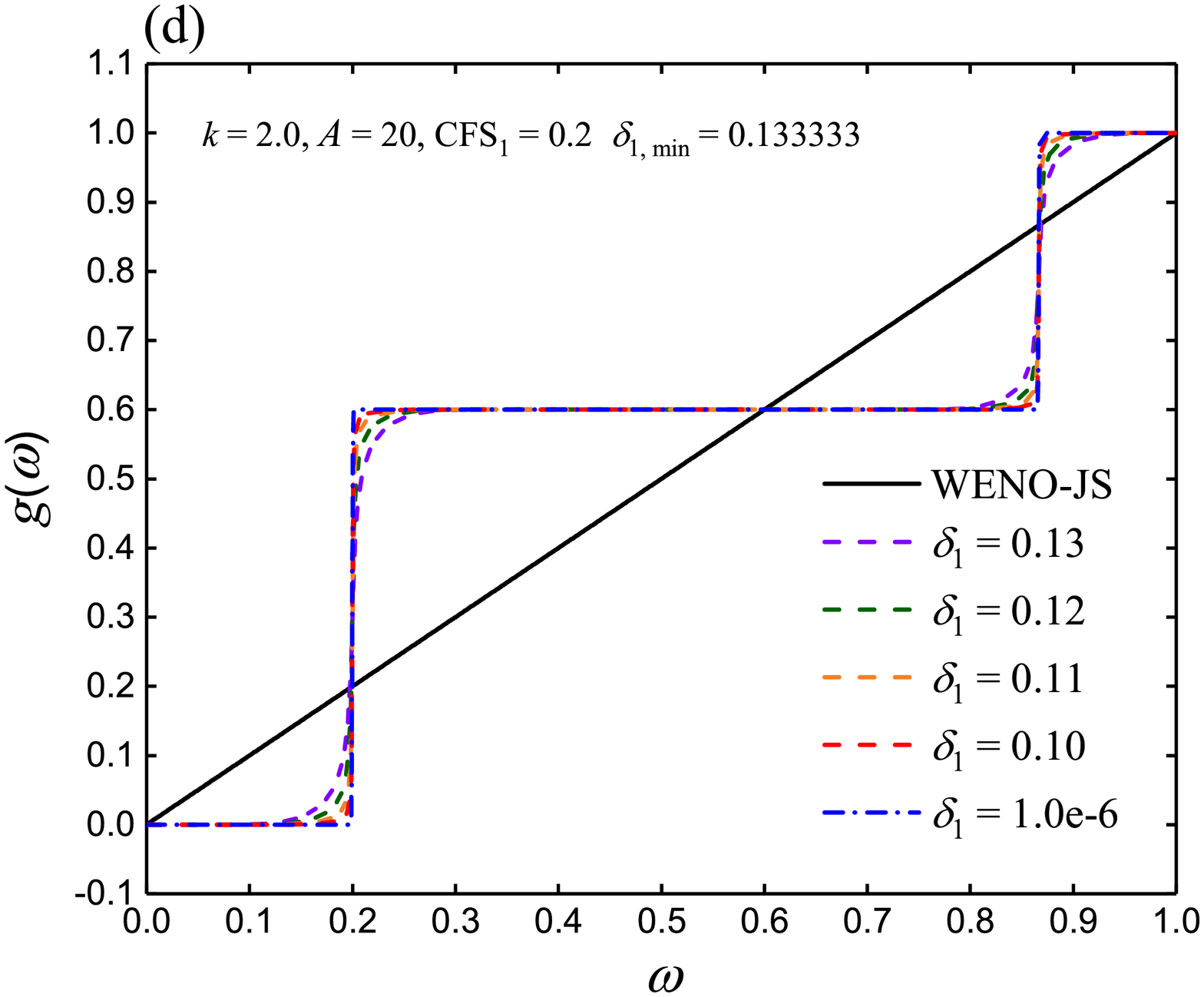} \\
\caption{The new mapping function $\big( g^{\mathrm{ACM}} \big)_{1}
(\omega)$, and effects of varying parameters $k, A, \mathrm{CFS}_{1}$
and $\delta_{1}$ for $d_{1} = 0.6$.}
\label{fig:parametricStudy}
\end{figure}

The curves of $\big(g^{\mathrm{ACM}}\big)_{s}(\omega)$ varying with 
$\omega$ and the effects of the parameters $k, A, \mathrm{CFS}_{s}$
and $\delta_{s}$ are shown in Fig. \ref{fig:parametricStudy}, taking 
the case of $d_{1} = 0.6$ as an example. 

In Fig. \ref{fig:parametricStudy}, we can see the following 
properties: (1) for given $A, \mathrm{CFS}_{1}$ and $\delta_{1}$, 
increasing $k$ will widen the optimal weight interval (standing for
the interval about $\omega = d_{s}$ over which the mapping process
attempts to use the corresponding optimal weight, see 
\cite{WENO-MAIMi}), but narrow the transition intervals (standing
 for the intervals about $\omega = \mathrm{CFS}_{s}$ or 
$\omega = \overline{\mathrm{CFS}}_{s}$ over which the mapping 
results satisfying $0< g(\omega) < d_{s}$ or $d_{s} < g(\omega) < 1$
respectively); (2) for given $k, \mathrm{CFS}_{1}$ and $\delta_{1}$, 
descreasing $A$ will widen the optimal weight interval and narrow the
transition intervals; (3) for given $k, A$ and $\delta_{1}$, 
increasing $\mathrm{CFS}_{1}$ will narrow the optimal weight interval
and maintain the width of the transition intervals unchanged; (4)
for given $k, A$ and $\mathrm{CFS}_{1}$, decreasing $\delta_{1}$ 
will widen the optimal weight interval and narrow the transition 
intervals.

For smooth problems, a wider optimal weight interval will bring the 
scheme closer to the corresponding linear upwind scheme leading to 
lower dissipation and higher resolution, and it is apparent from the
success of fifth-order WENO-IM($2,0.1$) scheme \cite{WENO-IM} that 
widening the optimal weight interval results in better performance.
However, an excessive optimal weight interval may 
lead to possible over-amplification of the contribution from a 
non-smooth stencil that creats a serious problem when a shock or 
nonlinear interaction between two shocks, like in the blastwave 
problem \cite{interactingBlastWaves-Woodward-Colella}, appears in 
the solution. It would generate numerical oscillations, even produce 
negative density and pressure, due to the mapping. This would be 
likely to happen in various mapped WENO schemes as their mappings 
push the nonlinear weights over the optimal weight intervals to be 
the ideal weights. Taking the WENO-ACM scheme as an example, its 
mapping functions compress all the large and small weights 
$\omega_{s} \in (\mathrm{CFS}_{s}+\delta_{s}, 
\overline{\mathrm{CFS}}_{s} - \delta_{s})$ closer together toward 
the ideal weights with less numerical dissipation but higher risk of 
generating numerical oscillations. It is easy to verify from 
Eq.(\ref{eq:mappingFunctionACM}) that the maximum optimal weight 
interval of the mapping function $\big(g^{\mathrm{ACM}}\big)_{s}(
\omega)$ is determined by $\mathrm{CFS}_{s}$ and we can also find 
this intuitively from Fig. \ref{fig:parametricStudy}. Therefore, in 
the WENO-ACM scheme, an optimal $\mathrm{CFS}_{s}$ should be desired 
that helps to obtain solutions with less numerical dissipation 
leading to higher resolution and 
prevent the scheme from generating numerical oscillations in the 
meantime. In other words, we can treat $\mathrm{CFS}_{s}$ as a 
tunable parameter: a larger $\mathrm{CFS}_{s}$ mekes the performance 
of the WENO-ACM scheme get closer to that of the WENO-M or even 
WENO-JS scheme, and a smaller $\mathrm{CFS}_{s}$ mekes the 
performance of the WENO-ACM scheme get closer to that of the linear 
upwind scheme. It is difficult and unsolved here yet to determine 
the optimal $\mathrm{CFS}_{s}$ theoretically, but after extensive 
numerical tests, we find that $\mathrm{CFS}_{s} = d_{s}/10$ should 
be a good choice. We will give a more detailed discussion about this 
through the blastwave problem as shown in Example \ref{Euler4} in 
Subsection \ref{subsec:1DEuler}.

Clearly, the new mapping method uses only one assignment operation 
when $\omega$ is out of the transition intervals. However, when 
$\omega$ is on the transition intervals, the new mapping method uses 
multiple multiplication and division as the existing mapped WENO 
schemes (like the WENO-M and WENO-PM6 schemes) do in their mapping 
processes on the whole interval of $\omega \in [0,1]$. Thus, 
narrower transition intervals will introduce fewer mathematical 
operations during the mapping process of the WENO-ACM scheme so that 
the CPU time will decrease significantly. Furthermore, our tests 
have shown that narrower transition intervals will not bring any 
adverse effects on the resolution and convergence rate of accuracy 
of the WENO-ACM scheme. Actually, the WENO-ACM scheme still performs 
very well even if the transition intervals decrease to near zero, 
and this will be verified in the calculation results of the numerical
experiments in Section \ref{secNumericalExperiments}.

\subsection{The rate of convergence}
Before giving Theorem \ref{theorem:convergenceRates} to show the
convergence rate property of the WENO-ACM scheme, we state the
following two necessary lemmas explicitly that can be found in the
statement of page 456 to 457 in \cite{WENO-IM} and page 565 in
\cite{WENO-M}.

\begin{lemma}
The sufficient condition for the $(2r - 1)$th-order WENO scheme to 
achieve the optimal order of accuracy is 
\begin{equation*}
\omega_{s} - d_{s} = O\Big((\Delta x)^{r}\Big), \quad r = 2,\cdots,9, s = 0,
\cdots,r-1.
\end{equation*}
\label{lemma:sufficient_condition}
\end{lemma}

\begin{lemma}
For $n_{\mathrm{cp}} \leq r - 1$, the weights
$\omega_{s}^{\mathrm{JS}}$ in the $(2r - 1)$th-order WENO-JS scheme 
satisfy
\begin{equation*}
\omega_{s}^{\mathrm{JS}} - d_{s} = O\Big((\Delta x)^{r - 1 - 
n_{\mathrm{cp}}}\Big), \quad r = 2, 3, \cdots, 9,
\label{Eq1_lemma_ncpWENO-JS}
\end{equation*} 
then, the convergence order is
\begin{equation*}
r_{\mathrm{c}} = \left\{
\begin{array}{ll}
2r - 1, & \mathrm{if} \quad n_{\mathrm{cp}} = 0, \\
2r - 2 - n_{\mathrm{cp}}, & \mathrm{if} \quad n_{\mathrm{cp}} = 1, 2
, \cdots, r - 1.
\end{array}\right.
\end{equation*} 
\label{convergence:WENO-JS}
\end{lemma}

\begin{theorem}
For $n_{\mathrm{cp}} < r - 1$, the $(2r-1)$th-order WENO-ACM scheme 
can achieve the optimal convergence rate of accuracy if the new
mapping function $\big( g^{\mathrm{ACM}} \big)_{s}(\omega)$ is
applied to weights of the $(2r-1)$th-order WENO-JS scheme.
\label{theorem:convergenceRates}
\end{theorem}

One can easily prove Theorem \ref{theorem:convergenceRates} by 
employing the Taylor series analysis and using Theorem 
\ref{theorem:mappingFunction:WENO-ACM}, Lemma 
\ref{lemma:sufficient_condition} and Lemma \ref{convergence:WENO-JS},
and the detailed proof process is almost identical to the one in 
\cite{WENO-M}.


\section{Numerical experiments}
\label{secNumericalExperiments}
In this section, we compare the calculation results of the WENO-ACM 
scheme with those of several typical WENO schemes, i.e., the 
classic WENO-JS scheme, the original WENO-M scheme, and the WENO-PM6 
scheme that generates less dissipation and better resolution.
The value of $\epsilon$ is chosen to be $10^{-40}$ 
in all considered schemes. 
Several problems with different initial and boundary conditions are 
used for the comparison, such as the one-dimensional linear 
advection equation, one- and two- dimensional Euler equations of 
compressible gas dynamics. The one-dimensional linear advection 
equation is computed for the accuracy test, and the numerical 
experiments of the Euler equations are conducted to demonstrate the 
performance of the WENO-ACM scheme in solving hyperbolic systems. 
Meanwhile, the computational costs of the WENO-JS, WENO-M, WENO-PM6 
and WENO-ACM schemes in simulating the benchmark problems of 
two-dimensional Euler equations are compared by the CPU timing per 
Runge-Kutta step.

In all the numerical experiments below, the global Lax-Friedrichs
numerical flux is employed, and the parameters in the WENO-ACM scheme
are chosen to be $k = 2, A = 20, \delta_{s} = 1.0\mathrm{e-}6, 
\mathrm{CFS}_{s} = 0.1d_{s}$. 

\subsection{One-dimensional linear advection equation}
\begin{example}
\bf{(Accuracy test without critical points \cite{WENO-IM})} 
\rm{We solve the one-dimensional linear advection equation $u_{t} + 
u_{x} = 0$ with the periodic boundary conditions and the following 
initial condition} 
\label{LAE1}
\end{example} 
\begin{equation}
u(x, 0) = \sin ( \pi x ). 
\label{eq:LAE:IC1}
\end{equation}
It is noted that although the initial 
condition in Eq.(\ref{eq:LAE:IC1}) has two first-order critical 
points, their first and third derivatives vanish simultaneously. To 
ensure that the error for the overall scheme is a measure of the 
spatial convergence only, and note that we consider only the 
fifth-order methods here, we set the CFL number to be 
$(\Delta x)^{2/3}$. The $L_{1}, L_{2}, L_{\infty}$ norms of the 
error are calculated by comparing the numerical solution 
$(u_{h})_{j}$ with the exact solution $u_{j}^{\mathrm{exact}}$
according to 
\begin{equation*}
\displaystyle
\begin{array}{l}
L_{1} = h \cdot \displaystyle\sum\limits_{j} \big\lvert u_{j}^{
\mathrm{exact}} - (u_{h})_{j} \big\rvert, \\
L_{2} = \sqrt{h \cdot \displaystyle\sum\limits_{j} (u_{j}^{
\mathrm{exact}} - (u_{h})_{j})^{2}}, \\
L_{\infty} = \displaystyle\max_{j} \big\lvert u_{j}^{
\mathrm{exact}} - (u_{h})_{j} \big\rvert,
\end{array}
\end{equation*}
where $h = \Delta x$ is the uniform spatial step size.

Table \ref{table_LAE1} shows the $L_{1}, L_{2}, L_{\infty}$ errors
and convergence orders of various considered WENO schemes for 
Example 1 at output time $t = 2.0$. The results of the three rows 
are $L_{1}$-, $L_{2}$- and $L_{\infty}$- norm errors and orders in 
turn (similarly hereinafter). All the schemes can achieve the 
optimal convergence orders. In terms of accuracy, the WENO-M, 
WENO-PM6 and WENO-ACM schemes provide more accurate numerical 
solutions than the results of the WENO-JS scheme in general. 
Besides, it is noted that, in terms of the $L_{1}$-norm error, the 
WENO-ACM scheme gives almost equally accurate numerical solutions as
those of the WENO-PM6 scheme which are more accurate than the 
numerical solutions of the WENO-M scheme.

\begin{table}[ht]
\begin{scriptsize}
\centering
\caption{Convergence properties of various considered schemes 
solving $u_{t} + u_{x} = 0$ with initial condition $u(x, 0) = \sin 
(\pi x)$.}
\label{table_LAE1}
\begin{tabular*}{\hsize}
{@{}@{\extracolsep{\fill}}lllllll@{}}
\hline
$h = \Delta x$       
& $0.2$               & $0.1$               & $0.05$        
& $0.025$             & $0.0125$            & $0.00625$ \\
\hline
WENO-JS 
& 6.18628e-02(-)      & 2.96529e-03(4.3821) & 9.27609e-05(4.9985) 
& 2.89265e-06(5.0031) & 9.03392e-08(5.0009) & 2.82330e-09(4.9999)  \\
\space
& 4.72306e-02(-)      & 2.42673e-03(4.2826) & 7.64332e-05(4.9887) 
& 2.33581e-06(5.0322) & 7.19259e-08(5.0213) & 2.23105e-09(5.0107)  \\
\space
& 4.87580e-02(-)      & 2.57899e-03(4.2408) & 9.05453e-05(4.8320) 
& 2.90709e-06(4.9610) & 8.85753e-08(5.0365) & 2.72458e-09(5.0228)  \\
WENO-M 
& 2.01781e-02(-)      & 5.18291e-04(5.2829) & 1.59422e-05(5.0228) 
& 4.98914e-07(4.9979) & 1.56021e-08(4.9990) & 4.99356e-10(4.9977)  \\
\space
& 1.55809e-02(-)      & 4.06148e-04(5.2616) & 1.25236e-05(5.0193) 
& 3.91875e-07(4.9981) & 1.22541e-08(4.9991) & 3.83568e-10(4.9976)  \\
\space
& 1.47767e-02(-)      & 3.94913e-04(5.2256) & 1.24993e-05(4.9816) 
& 3.91808e-07(4.9956) & 1.22538e-08(4.9988) & 3.83541e-10(4.9977)  \\
WENO-PM6 
& 1.74869e-02(-)      & 5.02923e-04(5.1198) & 1.59130e-05(4.9821) 
& 4.98858e-07(4.9954) & 1.56020e-08(4.9988) & 4.88355e-10(4.9977)  \\
\space
& 1.35606e-02(-)      & 3.95215e-04(5.1006) & 1.25010e-05(4.9825) 
& 3.91831e-07(4.9957) & 1.22541e-08(4.9989) & 3.83568e-10(4.9976)  \\
\space
& 1.27577e-02(-)      & 3.94515e-04(5.0151) & 1.24960e-05(4.9805) 
& 3.91795e-07(4.9952) & 1.22538e-08(4.9988) & 3.83543e-10(4.9977)  \\
WENO-ACM 
& 1.52184e-02(-)      & 5.02844e-04(4.9196) & 1.59130e-05(4.9818) 
& 4.98858e-07(4.9954) & 1.56020e-08(4.9988) & 4.88355e-10(4.9977)  \\
\space
& 1.19442e-02(-)      & 3.95138e-04(4.9178) & 1.25010e-05(4.9822) 
& 3.91831e-07(4.9957) & 1.22541e-08(4.9989) & 3.83568e-10(4.9976)  \\
\space
& 1.17569e-02(-)      & 3.94406e-04(4.8977) & 1.24960e-05(4.9801) 
& 3.91795e-07(4.9952) & 1.22538e-08(4.9988) & 3.83543e-10(4.9977)  \\
\hline
\end{tabular*}
\end{scriptsize}
\end{table}

\begin{example}
\bf{(Accuracy test with first-order critical points \cite{WENO-M})} 
\rm{We solve the one-dimensional linear advection equation $u_{t} + 
u_{x} = 0$ with the periodic boundary conditions and the following
initial condition}
\label{LAE2}
\end{example}
\begin{equation}
u(x, 0) = \sin \bigg( \pi x - \dfrac{\sin(\pi x)}{\pi} \bigg). 
\label{eq:LAE:IC2}
\end{equation}
As mentioned earlier, the CFL number is set to be $(\Delta x)^{2/3}$.
It is easy to verify that the particular initial condition 
Eq.(\ref{eq:LAE:IC2}) has two first-order critical points, which both
have a non-vanishing third derivative. 

The $L_{1}, L_{2}, L_{\infty}$ errors and convergence orders of 
various considered WENO schemes for Example \ref{LAE2} at output
time $t = 2.0$ are shown in Table \ref{table_LAE2}. From Table
\ref{table_LAE2}, we can observe that the WENO-M, WENO-PM6 and 
WENO-ACM schemes can retain the optimal orders even in the presence
of critical points. Moreover, in terms of accuracy, the WENO-ACM 
scheme provides the equally accurate results as those of the WENO-M 
and WENO-PM6 schemes, which are much more accurate than the 
solutions of the WENO-JS scheme whose $L_{\infty}$ convergence rate 
of accuracy drops by almost 2 orders leading to an overall accuracy 
loss shown with $L_{1}$ and $L_{2}$ convergence orders.

\begin{table}[ht]
\begin{scriptsize}
\centering
\caption{Convergence properties of various considered schemes 
solving $u_{t} + u_{x} = 0$ with initial condition $u(x, 0) = \sin 
( \pi x - \sin(\pi x)/\pi )$.}
\label{table_LAE2}
\begin{tabular*}{\hsize}
{@{}@{\extracolsep{\fill}}lllllll@{}}
\hline
$h = \Delta x$       
& $0.2$               & $0.1$               & $0.05$        
& $0.025$             & $0.0125$            & $0.00625$ \\
\hline
WENO-JS 
& 1.24488e-01(-)      & 1.01260e-02(3.6199) & 7.22169e-04(3.8096) 
& 3.42286e-05(4.3991) & 1.58510e-06(4.4326) & 7.95517e-08(4.3165)  \\
\space
& 1.09463e-01(-)      & 8.72198e-03(3.6496) & 6.76133e-04(3.6893) 
& 3.63761e-05(4.2162) & 2.29598e-06(3.9858) & 1.68304e-07(3.7700)  \\
\space
& 1.24471e-01(-)      & 1.43499e-02(3.1167) & 1.09663e-03(3.7099) 
& 9.02485e-05(3.6030) & 8.24022e-06(3.4531) & 8.31702e-07(3.3085)  \\
WENO-M 
& 7.53259e-02(-)      & 3.70838e-03(4.3443) & 1.45082e-04(4.6758) 
& 4.80253e-06(4.9169) & 1.52120e-07(4.9805) & 4.77083e-09(4.9948)  \\
\space
& 6.39017e-02(-)      & 3.36224e-03(4.2484) & 1.39007e-04(4.5962) 
& 4.52646e-06(4.9406) & 1.42463e-07(4.9897) & 4.45822e-09(4.9980)  \\
\space
& 7.49250e-02(-)      & 5.43666e-03(3.7847) & 2.18799e-04(4.6350) 
& 6.81451e-06(5.0049) & 2.14545e-07(4.9893) & 6.71080e-09(4.9987)  \\
WENO-PM6 
& 9.51313e-02(-)      & 4.82173e-03(4.3023) & 1.55428e-04(4.9552) 
& 4.87327e-06(4.9952) & 1.52750e-07(4.9956) & 4.77729e-09(4.9988)  \\
\space
& 7.83600e-02(-)      & 4.29510e-03(4.1894) & 1.43841e-04(4.9001) 
& 4.54036e-06(4.9855) & 1.42488e-07(4.9939) & 4.45807e-09(4.9983)  \\
\space
& 9.32356e-02(-)      & 5.91037e-03(3.9796) & 2.09540e-04(4.8180) 
& 6.83270e-06(4.9386) & 2.14532e-07(4.9932) & 6.71079e-09(4.9986)  \\
WENO-ACM 
& 8.75629e-02(-)      & 4.39527e-03(4.3163) & 1.52219e-04(4.8517) 
& 4.86436e-06(4.9678) & 1.52735e-07(4.9931) & 4.77728e-09(4.9987)  \\
\space
& 6.98131e-02(-)      & 4.02909e-03(4.1150) & 1.42172e-04(4.8247) 
& 4.53770e-06(4.9695) & 1.42486e-07(4.9931) & 4.45807e-09(4.9983)  \\
\space
& 7.91292e-02(-)      & 5.89045e-03(3.7478) & 2.09893e-04(4.8107) 
& 6.83017e-06(4.9416) & 2.14533e-07(4.9926) & 6.71079e-09(4.9986)  \\
\hline
\end{tabular*}
\end{scriptsize}
\end{table}

\begin{example}
\bf{(Accuracy test with discontinuous initial condition)} 
\rm{We solve the one-dimensional linear advection equation 
$u_{t} + u_{x} = 0$ with the following initial condition}
\label{LAE3}
\end{example} 
\begin{equation}
\begin{array}{l}
u(x, 0) = \left\{
\begin{array}{ll}
\dfrac{1}{6}\big[ G(x, \beta, z - \hat{\delta}) + 4G(x, \beta, z) + G
(x, \beta, z + \hat{\delta}) \big], & x \in [-0.8, -0.6], \\
1, & x \in [-0.4, -0.2], \\
1 - \big\lvert 10(x - 0.1) \big\rvert, & x \in [0.0, 0.2], \\
\dfrac{1}{6}\big[ F(x, \alpha, a - \hat{\delta}) + 4F(x, \alpha, a) +
 F(x, \alpha, a + \hat{\delta}) \big], & x \in [0.4, 0.6], \\
0, & \mathrm{otherwise},
\end{array}\right. 
\end{array}
\label{eq:LAE:IC3}
\end{equation}
where $G(x, \beta, z) = \mathrm{e}^{-\beta (x - z)^{2}}, F(x, \alpha
, a) = \sqrt{\max \big(1 - \alpha ^{2}(x - a)^{2}, 0 \big)}$, and 
the constants are $z = -0.7, \hat{\delta} = 0.005, \beta = \dfrac{
\log 2}{36\hat{\delta} ^{2}}, a = 0.5$ and $\alpha = 10$. The 
periodic boundary conditions are used in the two directions and the 
CFL number is set to be $0.1$. This example is obtained from \cite
{WENO-JS} and it consists of a Gaussian, a square wave, a sharp 
triangle and a semi-ellipse.

Table \ref{table_LAE3} shows the $L_{1}, L_{2}, L_{\infty}$ errors 
and convergence orders of various considered WENO schemes for 
Example \ref{LAE3} at output times $t = 2$ and $t = 2000$. At output 
time $t = 2$, we find: (1) for all schemes, the $L_{1}$ and $L_{2}$ 
orders are approximately $1.0$ and $0.4$ to $0.5$, respectively, and 
the $L_{\infty}$ orders are all negative; (2) in terms of accuracy,
the WENO-ACM scheme provides the most accurate results closely 
followed by the WENO-PM6 and WENO-M schemes, which are more 
accurate than that of the WENO-JS scheme. At output time $t = 2000$, 
we find: (1) for the WENO-JS and WENO-M schemes, the $L_{1}$, $L_{2}$
orders decrease to very small values and even become negative; (2) 
however, for the WENO-ACM and WENO-PM6 schemes, their $L_{1}$ orders 
are maintained at approximately $1.0$, and their $L_{2}$ orders 
increase to approximately $0.6$ to $0.7$; (3) for all schemes, the
$L_{\infty}$ orders are very small and even become negative. 
Overall, for both short and long output times, the WENO-ACM scheme
performs as well as the WENO-PM6 scheme in calculating this problem 
that includes various discontinuities.

\begin{table}[ht]
\begin{scriptsize}
\centering
\caption{Convergence properties of various considered schemes 
solving $u_{t} + u_{x} = 0$ with initial condition 
Eq.(\ref{eq:LAE:IC3}).}
\label{table_LAE3}
\begin{tabular*}{\hsize}
{@{}@{\extracolsep{\fill}}lllllll@{}}
\hline
\space  &\multicolumn{3}{l}{t = 2}  &\multicolumn{3}{l}{t = 2000}  \\
\cline{2-4}  \cline{5-7}
$h = \Delta x$
& $0.01$         & $0.005$             & $0.0025$ 
& $0.01$         & $0.005$             & $0.0025$   \\
\hline
WENO-JS     
& 6.30497e-02(-) & 2.81654e-02(1.2103)  & 1.41364e-02(0.9945)  
& 6.12899e-01(-) & 5.99215e-01(0.0326)  & 5.50158e-01(0.1232)   \\
\space
& 1.08621e-01(-) & 7.71111e-02(0.4943)  & 5.69922e-02(0.4362)  
& 5.08726e-01(-) & 5.01160e-01(0.0216)  & 4.67585e-01(0.1000)   \\
\space
& 4.09733e-01(-) & 4.19594e-01(-0.0343) & 4.28463e-01(-0.0302)  
& 7.99265e-01(-) & 8.20493e-01(-0.0378) & 8.14650e-01(0.0103)   \\
WENO-M     
& 4.77201e-02(-) & 2.23407e-02(1.0949)  & 1.11758e-02(0.9993)  
& 3.81597e-01(-) & 3.25323e-01(0.2302)  & 3.48528e-01(-0.0994)  \\
\space
& 9.53073e-02(-) & 6.91333e-02(0.4632)  & 5.09232e-02(0.4411)  
& 3.59205e-01(-) & 3.12970e-01(0.1988)  & 3.24373e-01(-0.0516)  \\
\space
& 3.94243e-01(-) & 4.05856e-01(-0.0419) & 4.16937e-01(-0.0389)  
& 6.89414e-01(-) & 6.75473e-01(0.0295)  & 6.25645e-01(0.1106)   \\
WENO-PM6     
& 4.66681e-02(-) & 2.13883e-02(1.1256)  & 1.06477e-02(1.0063)  
& 2.17323e-01(-) & 1.05197e-01(1.0467)  & 4.47030e-02(1.2347)   \\
\space
& 9.45566e-02(-) & 6.82948e-02(0.4694)  & 5.03724e-02(0.4391)  
& 2.28655e-01(-) & 1.47518e-01(0.6323)  & 9.34250e-02(0.6590)   \\
\space
& 3.96866e-01(-) & 4.06118e-01(-0.0332) & 4.15277e-01(-0.0322)  
& 5.63042e-01(-) & 5.04977e-01(0.1570)  & 4.71368e-01(0.0994)   \\
WENO-ACM     
& 4.45059e-02(-) & 2.03633e-02(1.1280)  & 1.02139e-02(0.9954)  
& 2.21313e-01(-) & 1.06583e-01(1.0541)  & 4.76305e-02(1.1620)   \\
\space
& 9.24356e-02(-) & 6.69718e-02(0.4649)  & 4.95672e-02(0.4342)  
& 2.28433e-01(-) & 1.46401e-01(0.6418)  & 9.40930e-02(0.6378)   \\
\space
& 3.92505e-01(-) & 4.03456e-01(-0.0397) & 4.13217e-01(-0.0345)  
& 5.36234e-01(-) & 5.03925e-01(0.0897)  & 5.15924e-01(-0.0339)   \\
\hline
\end{tabular*}
\end{scriptsize}
\end{table}

\begin{example}
\bf{(High resolution performance test with high-order critical
points)} 
\rm{We solve the one-dimensional linear advection equation 
$u_{t} + u_{x} = 0$ with the periodic boundary conditions and
the following initial condition \cite{WENO-PM}} 
\label{LAE4}
\end{example} 
\begin{equation}
u(x, 0) = \sin^{9} ( \pi x ) 
\label{eq:LAE:IC4}
\end{equation}
Again, the CFL number is set to be $(\Delta x)^{2/3}$. It is easy to
verify that the initial condition in Eq.(\ref{eq:LAE:IC4}) has
high-order critical points. 

Table \ref{table_LAE4} shows the $L_{1}, L_{2}, L_{\infty}$ errors 
of various considered WENO schemes for Example \ref{LAE4} at several
output times with a uniform mesh size of $\Delta x = 1/200$.
Clearly, at short output times, the WENO-ACM scheme achieves similar 
results as those of the WENO-M and WENO-PM6 schemes. However, at 
long output times, the numerical solutions computed by the WENO-M
scheme are far less accurate than that of the WENO-PM6 scheme, while 
the solutions of the WENO-ACM scheme are still as accurate as 
results of the WENO-PM6 scheme. Another observation is that the 
WENO-M, WENO-PM6 and WENO-ACM schemes all provide more accurate
numerical solutions than the WENO-JS scheme.

\begin{figure}[ht]
\centering
  \includegraphics[height=0.39\textwidth]
  {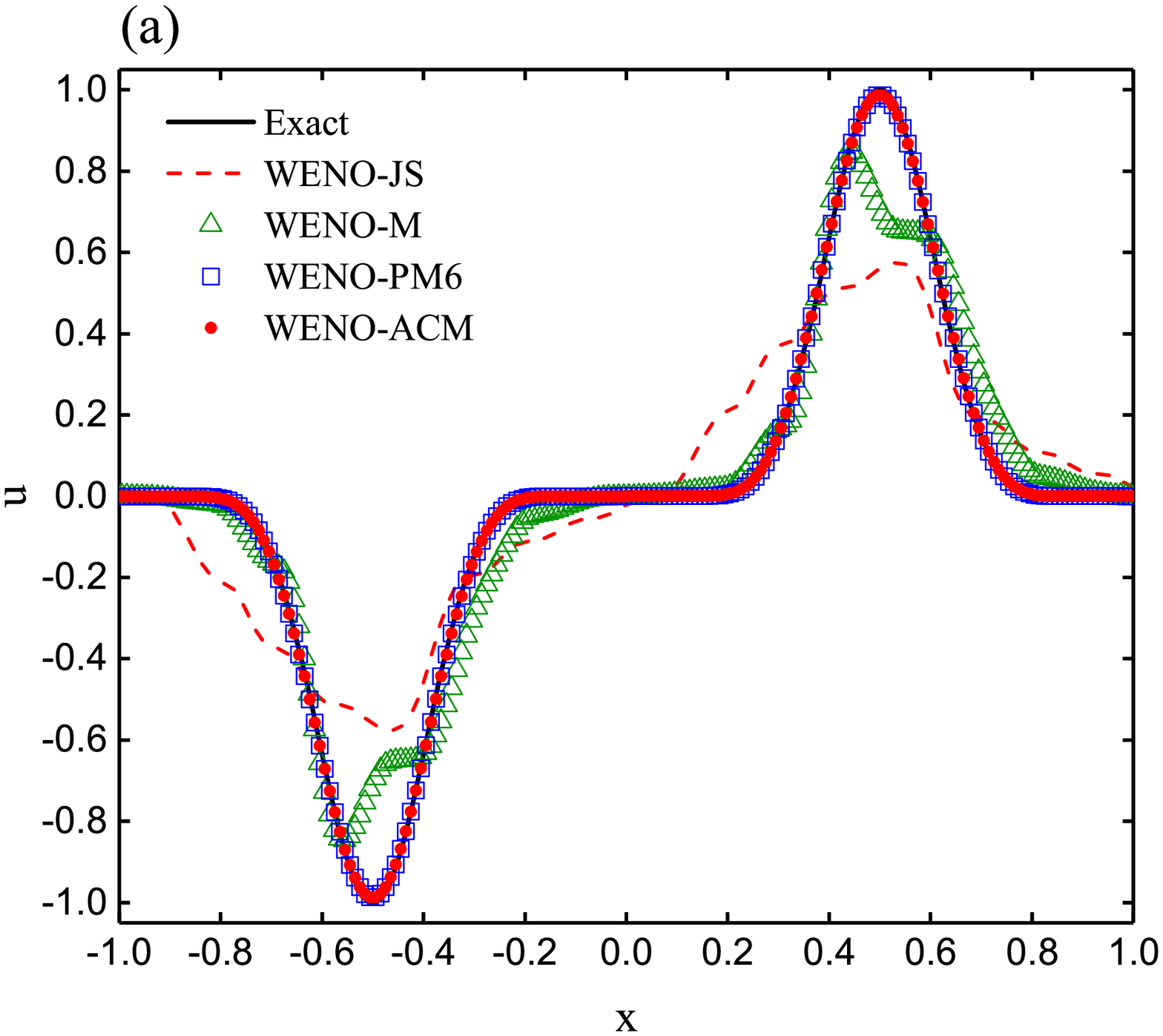}
  \includegraphics[height=0.39\textwidth]
  {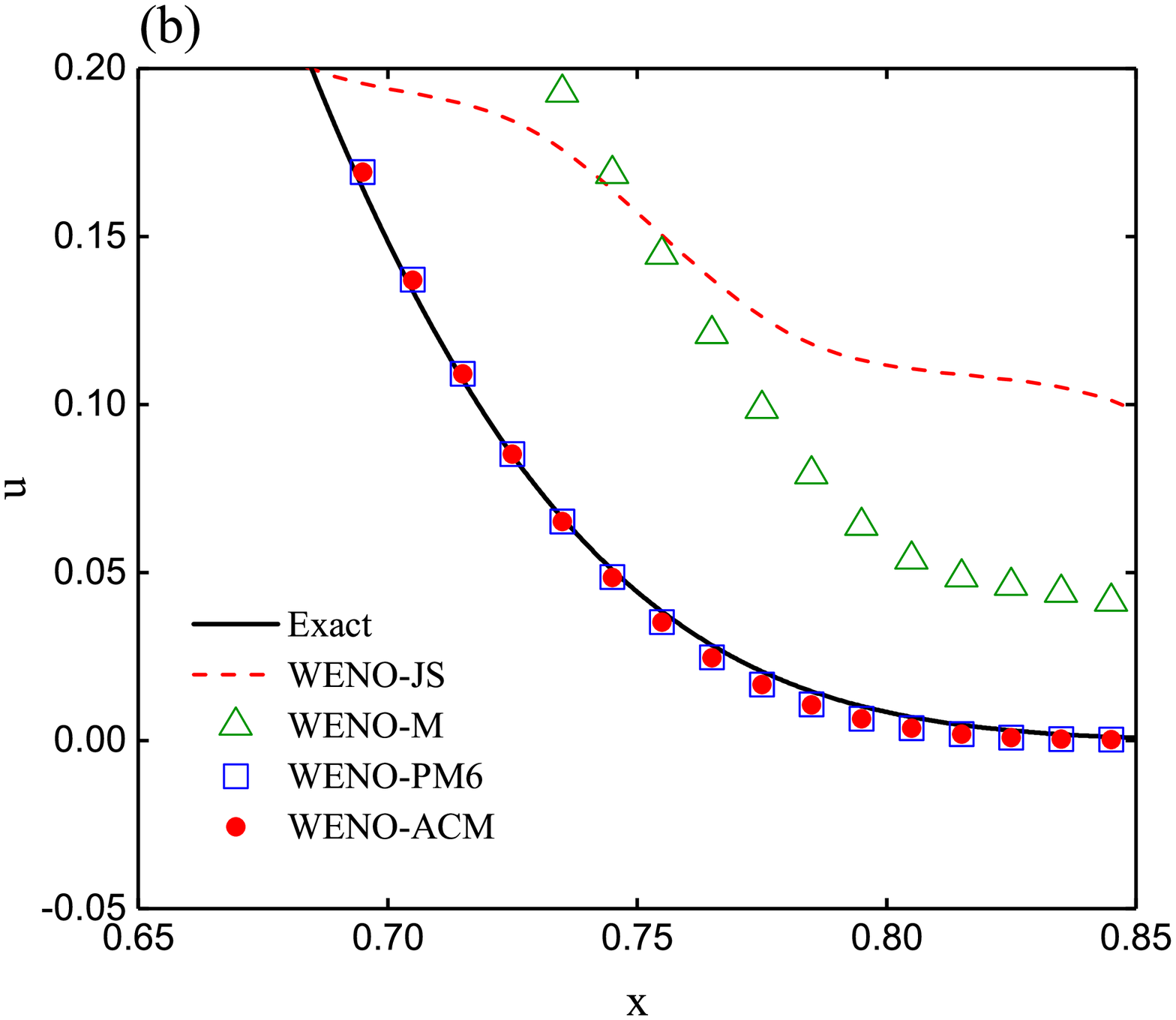}
\caption{Performance of the WENO-JS, WENO-M, WENO-PM6 and WENO-ACM
schemes for Example \ref{LAE4} at output time $t = 1000$ with a 
uniform mesh size of $\Delta x= 1/200$.}
\label{fig:ex:LAE4}
\end{figure}

Fig. \ref{fig:ex:LAE4} shows the performance of the WENO-JS, 
WENO-M, WENO-PM6 and WENO-ACM schemes for Example \ref{LAE4} at 
output time $t = 1000$ with a uniform mesh size of $\Delta x= 1/200$.
Clearly the WENO-ACM and WENO-PM6 schemes give the highest 
resolution followed by the WENO-M scheme whose resolution decreases 
significantly, and the WENO-JS scheme shows the lowest resolution.

\begin{table}[ht]
\begin{scriptsize}
\centering
\caption{Performance of various considered schemes solving $u_{t} + 
u_{x} = 0$ with $u(x, 0) = \sin^{9} (\pi x), \Delta x= 1/200$.}
\label{table_LAE4}
\begin{tabular*}{\hsize}
{@{}@{\extracolsep{\fill}}lllll@{}}
\hline
Scheme      & WENO-JS     & WENO-M      & WENO-PM6    & WENO-ACM \\
\hline
$t = 1$     & 3.87826e-05 & 8.84565e-06 & 8.52448e-06 & 8.43356e-06\\
\space      & 3.62689e-05 & 8.31248e-06 & 8.22944e-06 & 8.20366e-06\\
\space      & 6.69118e-05 & 1.38461e-05 & 1.38389e-05 & 1.38389e-05\\
$t = 10$    & 3.86931e-04 & 8.90890e-05 & 8.40259e-05 & 8.42873e-05\\
\space      & 3.52611e-04 & 8.32089e-05 & 8.19676e-05 & 8.19107e-05\\
\space      & 5.36940e-04 & 1.38348e-04 & 1.38205e-04 & 1.38205e-04\\
$t = 30$    & 1.17988e-03 & 2.73430e-04 & 2.51117e-04 & 2.52378e-04\\
\space      & 1.06511e-03 & 2.51737e-04 & 2.45084e-04 & 2.45090e-04\\
\space      & 1.58134e-03 & 4.13887e-04 & 4.13397e-04 & 4.13398e-04\\
$t = 50$    & 2.05488e-03 & 4.81901e-04 & 4.17588e-04 & 4.19825e-04\\
\space      & 1.84782e-03 & 4.39983e-04 & 4.07311e-04 & 4.07429e-04\\
\space      & 2.69500e-03 & 6.87879e-04 & 6.86969e-04 & 6.86983e-04\\
$t = 100$   & 5.42288e-03 & 1.29154e-03 & 8.30374e-04 & 8.35747e-04\\
\space      & 5.17716e-03 & 1.28740e-03 & 8.09152e-04 & 8.09679e-04\\
\space      & 1.20056e-02 & 3.32665e-03 & 1.36410e-03 & 1.36404e-03\\
$t = 200$   & 2.35657e-02 & 5.74021e-03 & 1.63963e-03 & 1.65557e-03\\
\space      & 2.68753e-02 & 7.66721e-03 & 1.59697e-03 & 1.59929e-03\\
\space      & 6.47820e-02 & 2.37125e-02 & 2.68938e-03 & 2.68955e-03\\
$t = 500$   & 1.55650e-01 & 4.89290e-02 & 3.88864e-03 & 3.95849e-03\\
\space      & 1.46859e-01 & 6.23842e-02 & 3.83159e-03 & 3.84802e-03\\
\space      & 2.57663e-01 & 1.78294e-01 & 6.45650e-03 & 6.45564e-03\\
$t = 1000$  & 2.91359e-01 & 1.34933e-01 & 7.17606e-03 & 7.24723e-03\\
\space      & 2.66692e-01 & 1.46524e-01 & 7.19008e-03 & 7.21626e-03\\
\space      & 4.44664e-01 & 3.17199e-01 & 1.21637e-02 & 1.21593e-02\\
\hline
\end{tabular*}
\end{scriptsize}
\end{table}

\subsection{One-dimensional Euler system}\label{subsec:1DEuler}
In this subsection, we calculate the one-dimensional Euler system of 
gas dynamics with different initial and boundary conditions. The 
one-dimensional Euler system is given by the following strong 
conservation form of mass, momentum and energy
\begin{equation}
\begin{array}{ll}
\begin{aligned}
&\dfrac{\partial \rho}{\partial t} + \dfrac{\partial (\rho u)}
{\partial x} = 0, \\
&\dfrac{\partial (\rho u)}{\partial t} + \dfrac{\partial (\rho u^{2}
+ p)}{\partial x} = 0, \\
&\dfrac{\partial E}{\partial t} + \dfrac{\partial (uE + up)}
{\partial x} = 0, \\
\end{aligned}
\end{array}
\label{1DEulerEquations}
\end{equation}
where $\rho, u, p$ and $E$ are the density, velocity, pressure and
total energy, respectively. The Euler system 
Eq.(\ref{1DEulerEquations}) is closed by the equation of state for 
an ideal polytropic gas, which is given by
\begin{equation*}
p = (\gamma - 1)\Big( E - \dfrac{1}{2}\rho u^{2} \Big),
\end{equation*}
where $\gamma$ is the ratio of specific heat, and we use $\gamma=1.4$
in this paper. The finite volume version of the characteristic-wise 
one-dimensional WENO procedure is employed, and we refer to 
\cite{FVMaccuracyProofs03} for details. In all examples of this 
subsection, the CFL number is set to be $0.5$.

\begin{example}
\bf{(Sod's shock tube problem)} 
\rm{We consider the Sod's shock tube problem 
\cite{SodShock-tubeProblem}, specified by the following initial 
condition}
\label{Euler1}
\end{example}
\begin{equation*}
\big( \rho, u, p \big)(x, 0) =\left\{
\begin{array}{ll}
(1.0, 0.0, 1.0), & x \in [0.0, 0.5], \\
(0.125, 0.0, 0.1), & x \in [0.5, 1.0].
\end{array}\right.
 \label{initial_1DEuler1}
\end{equation*}
The transmissive boundary conditions are used in two directions.

\begin{figure}[ht]
\centering
  \includegraphics[height=0.39\textwidth]
  {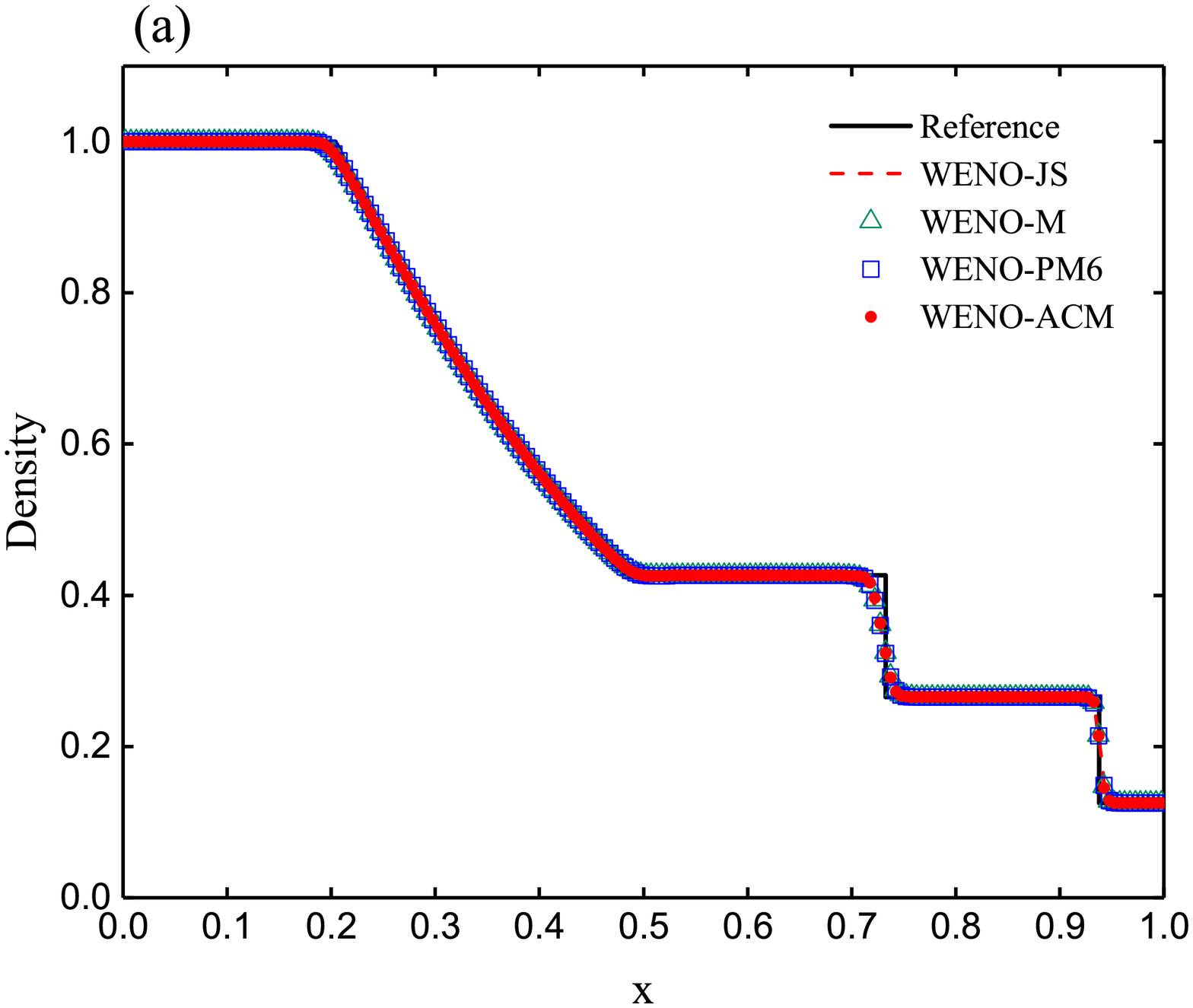}
  \includegraphics[height=0.39\textwidth]
  {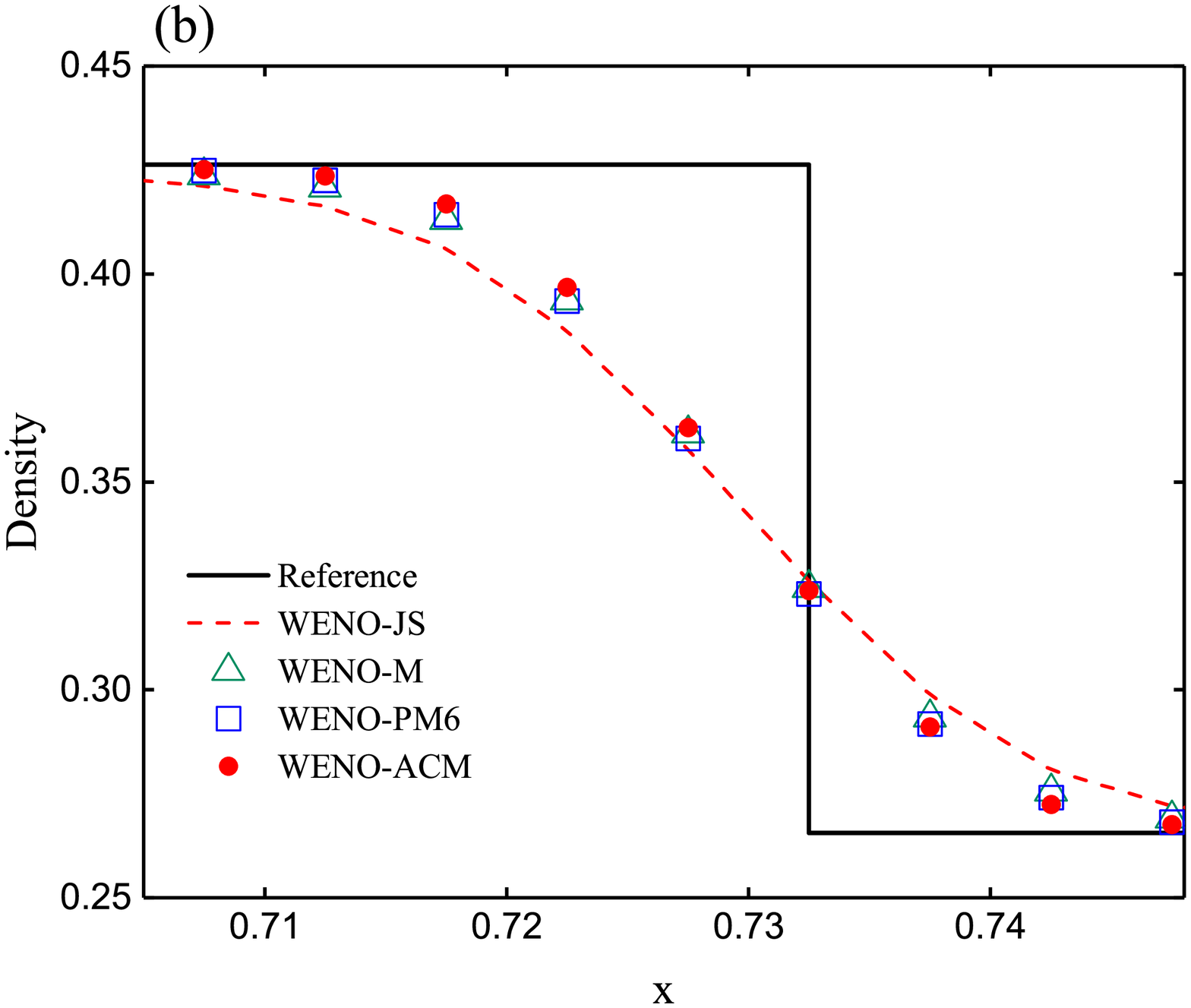}
\caption{The density profiles of the Sod's shock tube problem.}
\label{fig:ex:Sod-shock-tube}
\end{figure}

Fig. \ref{fig:ex:Sod-shock-tube} presents the density profiles
computed by the WENO-JS, WENO-M, WENO-PM6 and WENO-ACM schemes at 
output time $t = 0.25$ with a uniform mesh size of $N = 200$. We 
observe that the WENO-M, WENO-PM6 and WENO-ACM schemes capture 
sharper discontinuity compared to the WENO-JS scheme, and the 
WENO-ACM scheme gives slightly better resolution than the WENO-PM6 
and WENO-M schemes.

\begin{example}
\bf{(Riemann problem of Lax)} 
\rm{The second considered problem modeled by the one-dimensional 
Euler system Eq.(\ref{1DEulerEquations}) is the Lax's problem 
\cite{LaxShock-tubeProblem}, specified by the following initial 
condition}
\label{Euler2}
\end{example}
\begin{equation*}
\big( \rho, u, p \big)(x, 0) =\left\{
\begin{array}{ll}
(0.445, 0.698, 3.528), & x \in [-5, 0], \\
(0.500, 0.000, 0.571), & x \in [0, 5].
\end{array}\right.
 \label{initial_1DEuler2}
\end{equation*}
The transmissive boundary conditions are used at $x = \pm 5$, and
the uniform cell number is chosen to be $N = 200$.

\begin{figure}[ht]
\centering
  \includegraphics[height=0.39\textwidth]
  {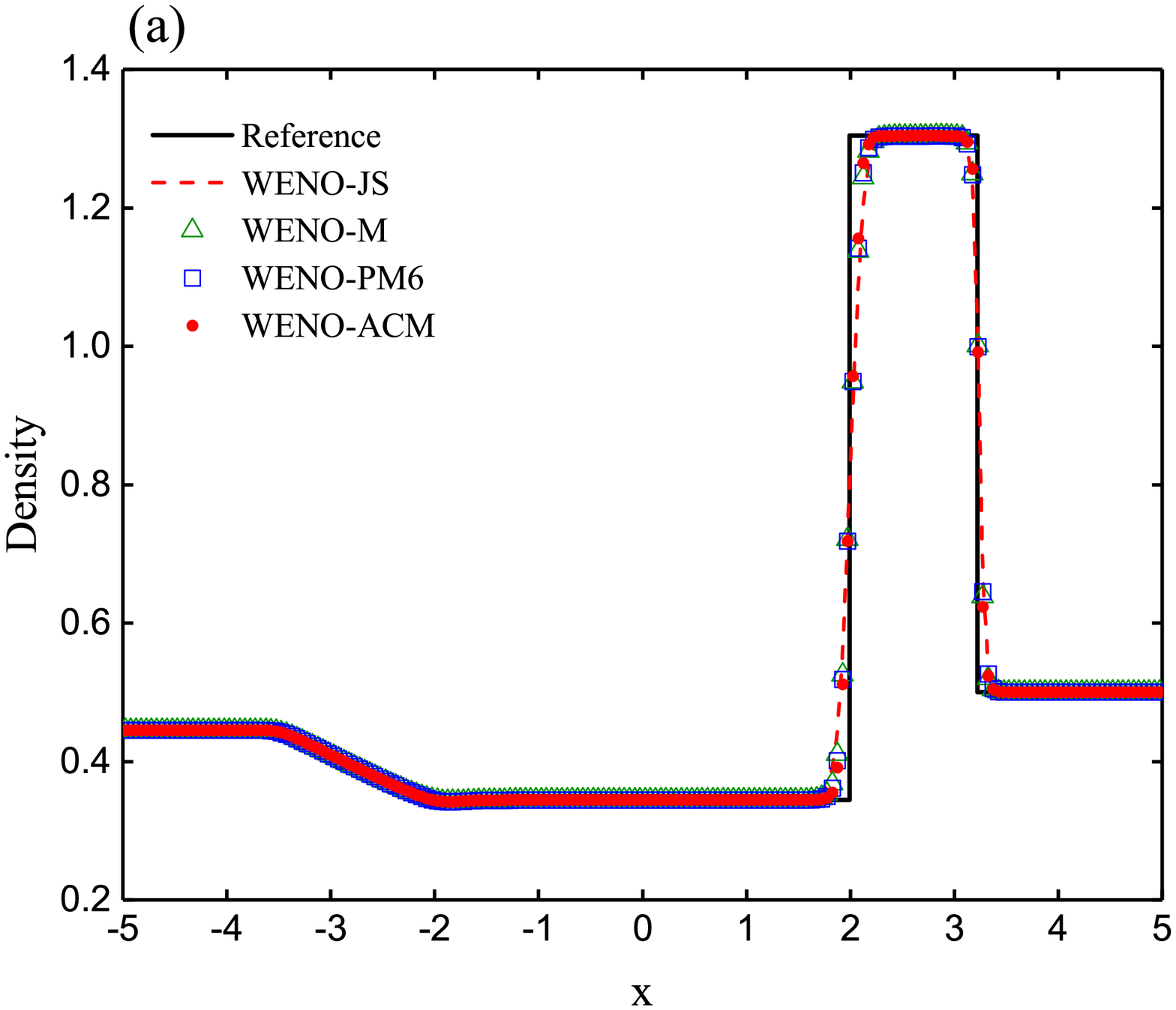}
  \includegraphics[height=0.39\textwidth]
  {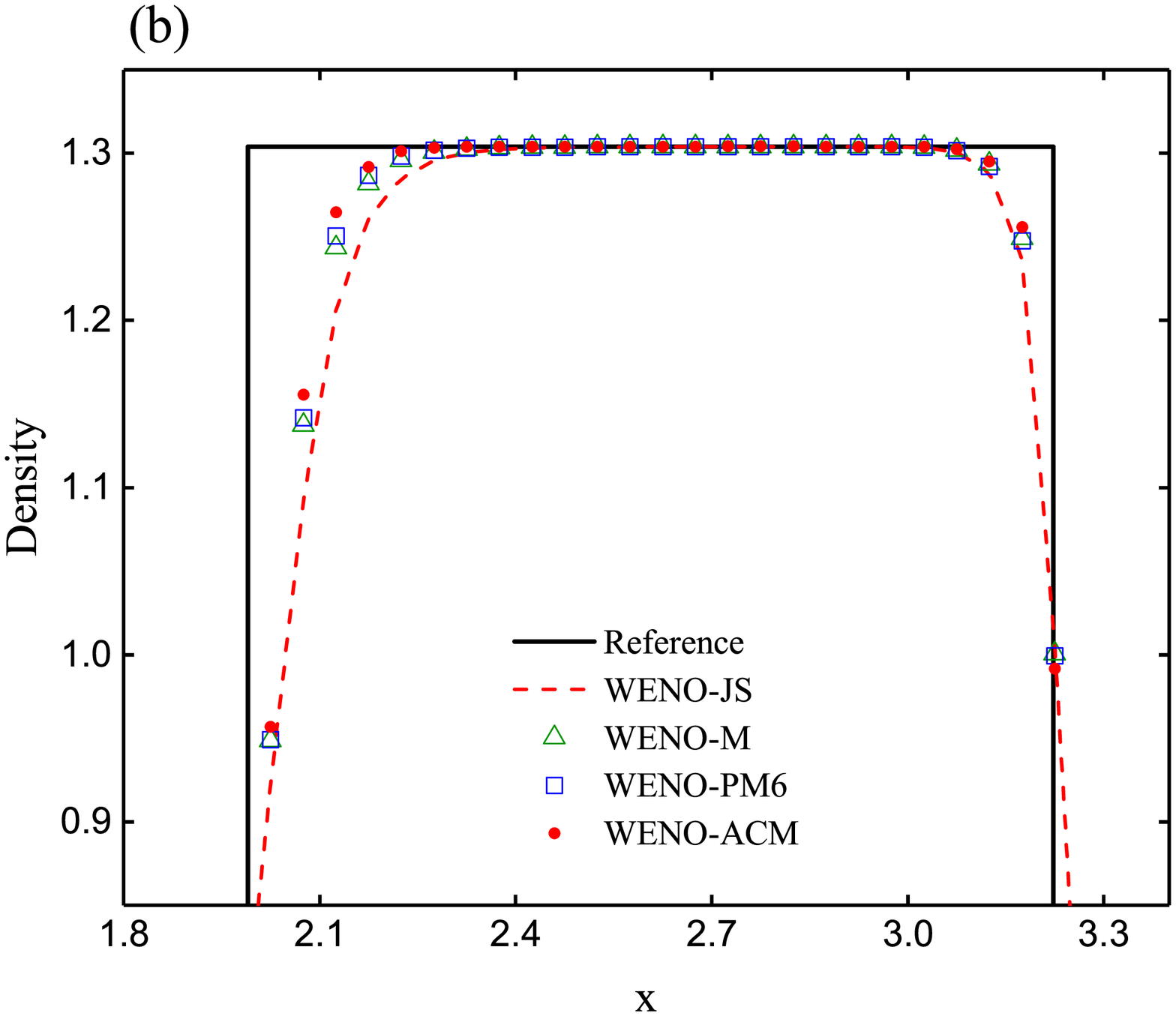}
\caption{The density profiles of the Riemann problem of Lax.}
\label{fig:ex:Lax-shock-tube}
\end{figure}

Fig. \ref{fig:ex:Lax-shock-tube} presents the density profiles
computed by the WENO-JS, WENO-M, WENO-PM6 and WENO-ACM schemes at 
output time $t = 1.3$. It is observed that the WENO-M, WENO-PM6 and 
WENO-ACM schemes have higher resolution than the WENO-JS scheme near 
the discontinuity. Also, if we take a closer look at $x\in (1.8,3.3)$
, it demonstrates that the WENO-ACM scheme performs slightly better 
than the WENO-PM6 and WENO-M schemes.

\begin{example}
\bf{(Mach 3 shock-density wave interaction)} 
\rm{We solve the Mach 3 shock-density wave interaction 
\cite{ENO-Shu1989}, whose solution consists of several shocklets and 
fine-scale structures that are located behind a main right-going 
shock \cite{WENO-Z}. Its boundaries at $x = \pm 5$ are specified by 
zero-gradient boundary condition, and its initial condition is given 
by}
\label{Euler3}
\end{example}
\begin{equation*}
\big( \rho, u, p \big)(x, 0) =\left\{
\begin{array}{ll}
(3.857143, 2.629369, 10.333333), & x \in [-5.0, -4.0], \\
(1.0 + 0.2\sin(5x), 0, 1), & x \in [-4.0, 5.0].
\end{array}\right.
 \label{initial_1DEuler3}
\end{equation*}

\begin{figure}[ht]
\centering
  \includegraphics[height=0.39\textwidth]
  {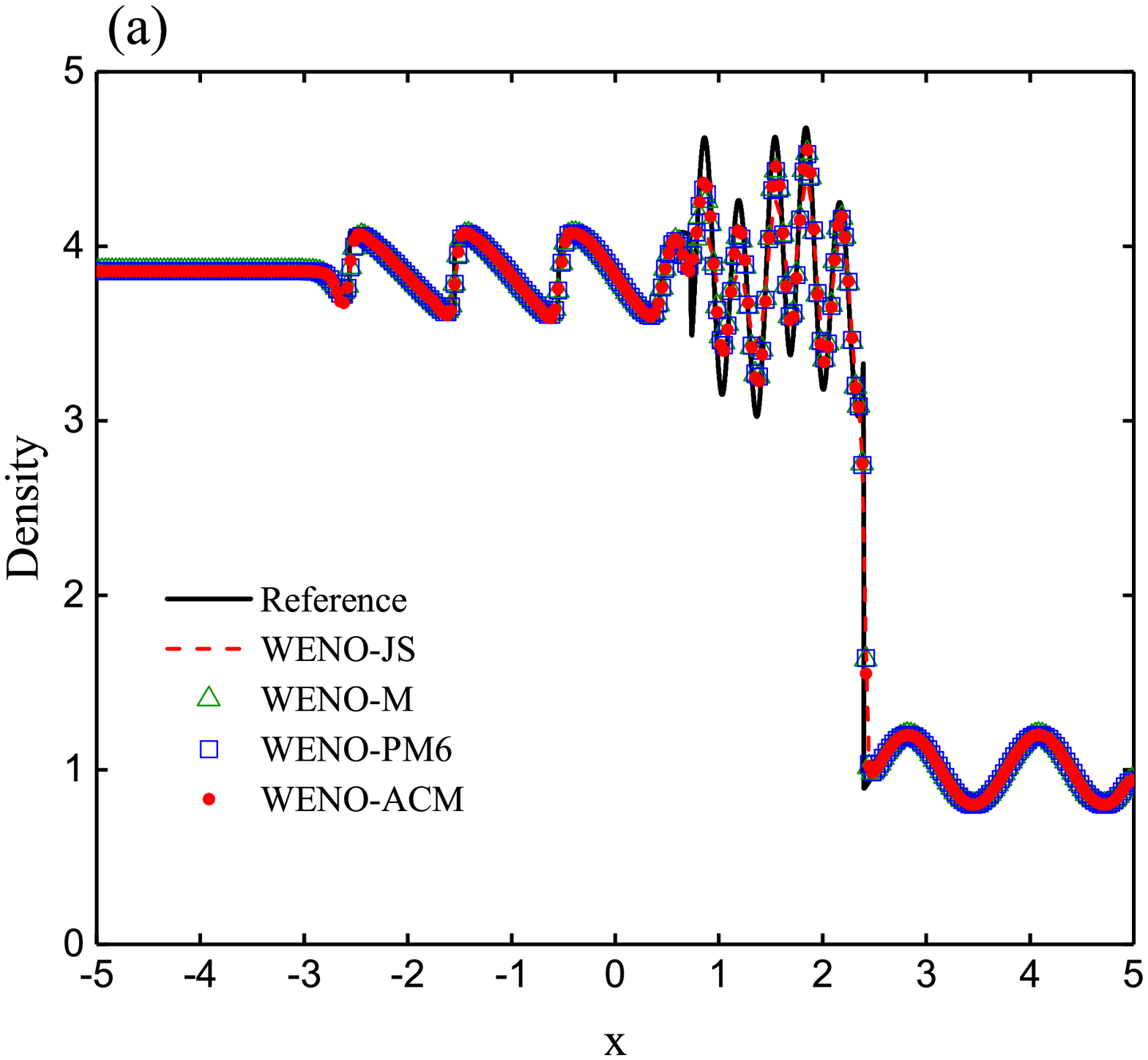}
  \includegraphics[height=0.39\textwidth]
  {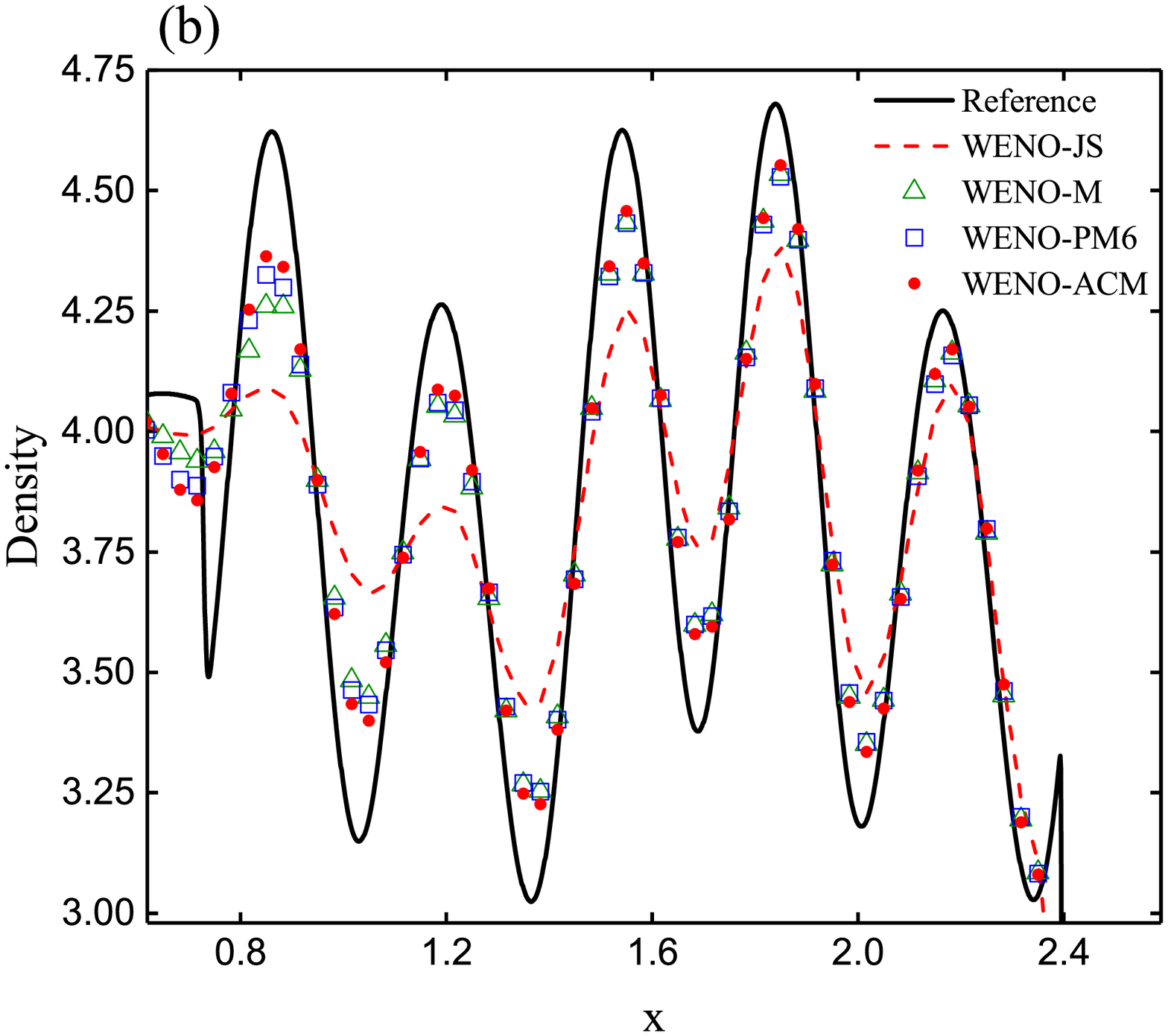}
\caption{The density profiles of the Mach 3 shock-density wave 
interaction.}
\label{fig:ex:Mach3Shock-density}
\end{figure}

Fig. \ref{fig:ex:Mach3Shock-density} gives the comparison on 
density between the WENO-JS, WENO-M, WENO-PM6 and WENO-ACM schemes 
at output time $t=1.8$ with the uniform cell number $N = 300$. The 
solid line is the reference solution using the WENO-JS scheme with 
$N = 10000$. The WENO-M, WENO-PM6 and WENO-ACM schemes capture much 
more fine-scale structures of the solution than the WENO-JS scheme. 
Furthermore, the WENO-ACM scheme shows the best description near 
shocklets and high-frequency waves behind the main shock.

\begin{example}
\bf{(Woodward-Colella interacting blastwaves)} 
\rm{We solve the standard test problem first used by Woodward and
Colella \cite{interactingBlastWaves-Woodward-Colella}. The reflective
boundary conditions are used at $x = 0, 1$, and the initial
condition is set to be}
\label{Euler4}
\end{example}
\begin{equation*}
\big( \rho, u, p \big)(x, 0) =\left\{
\begin{array}{ll}
(1, 0, 1000),   & x \in [0, 0.1),   \\
(1, 0, 0.01), & x \in [0.1, 0.9), \\
(1, 0, 100), & x \in [0.9, 1.0].
\end{array}\right.
 \label{initial_1DEuler4}
\end{equation*}

\begin{figure}[ht]
\centering
  \includegraphics[height=0.39\textwidth]
  {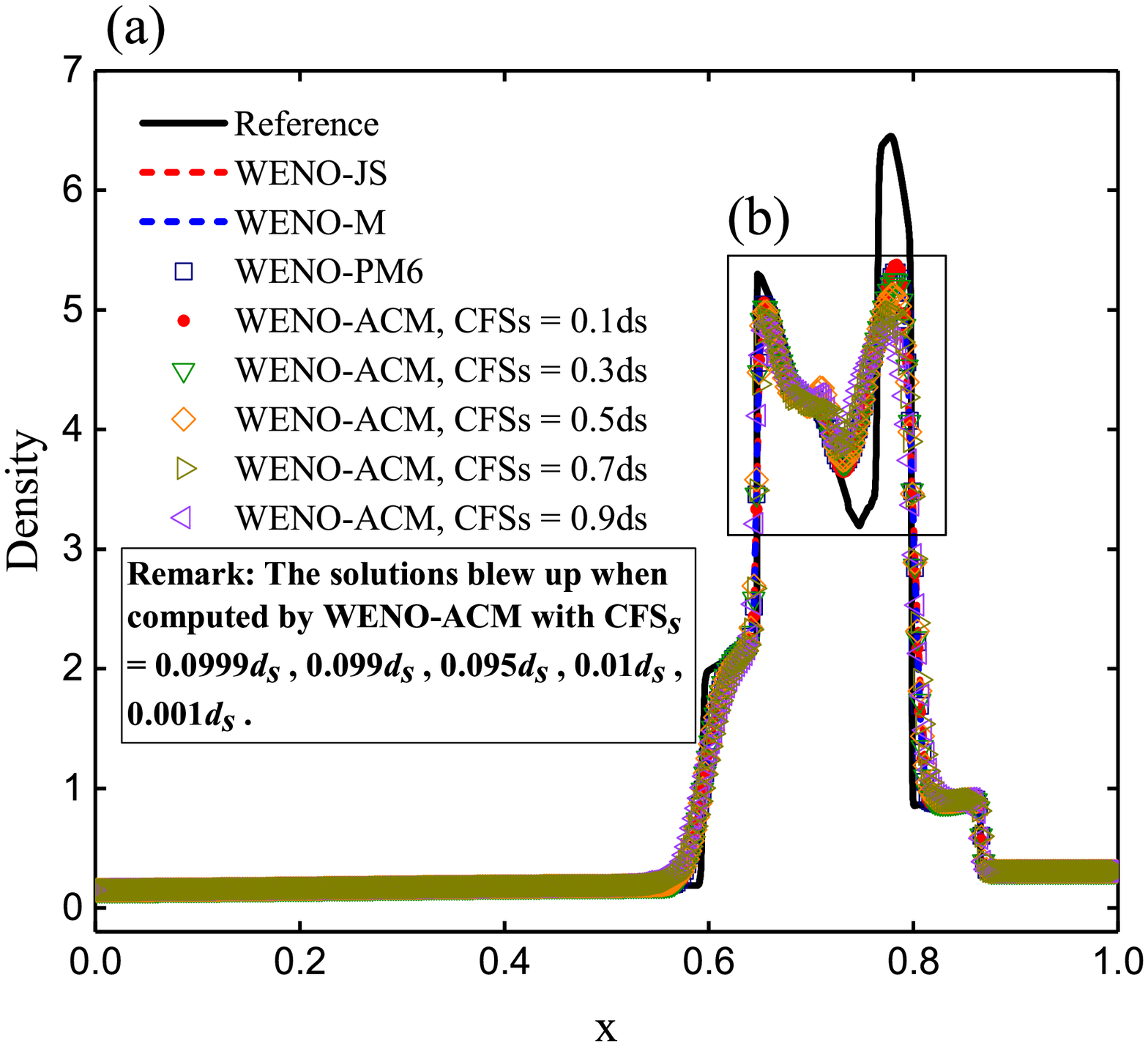}
  \includegraphics[height=0.39\textwidth]
  {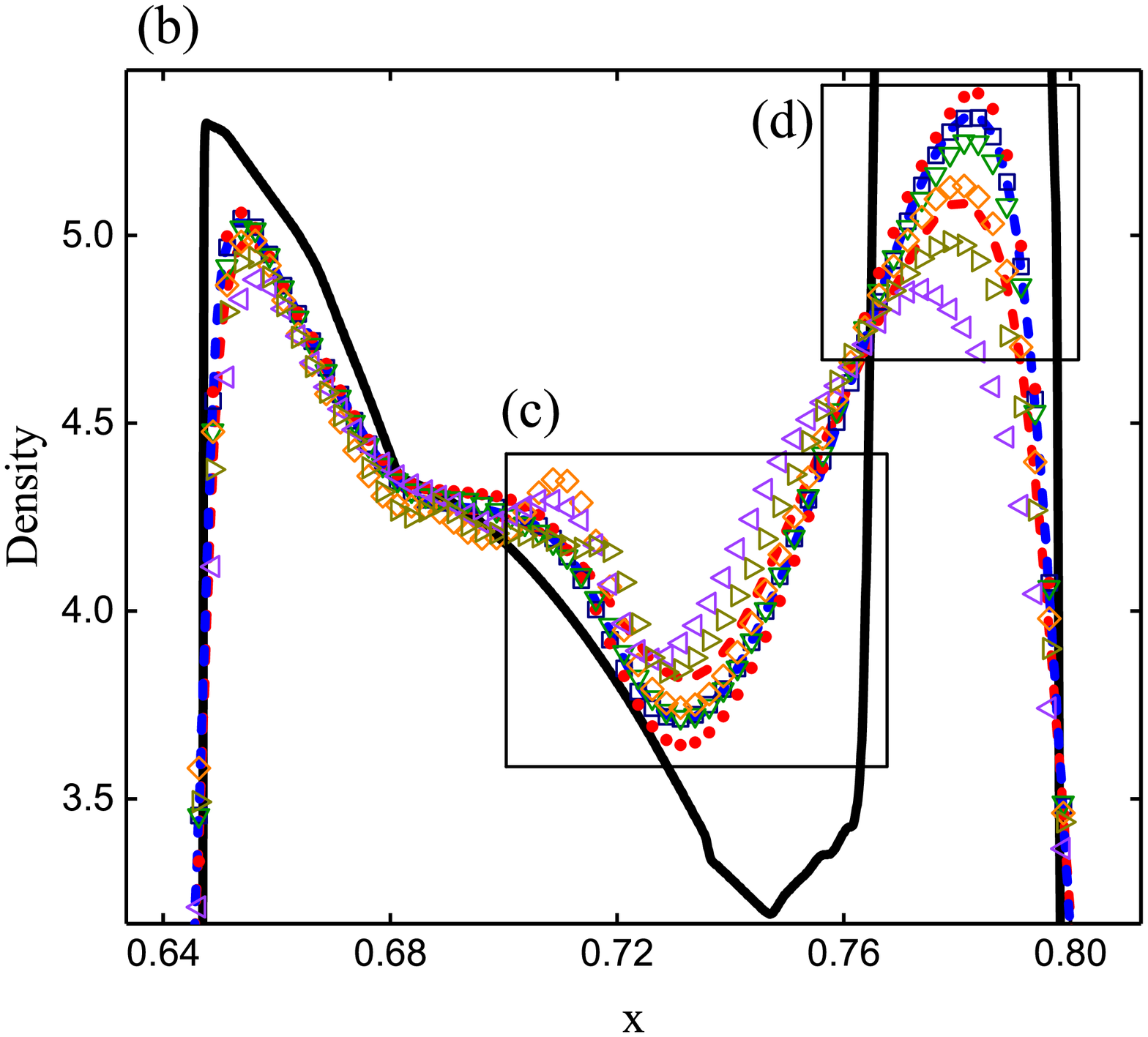}\\
  \includegraphics[height=0.39\textwidth]
  {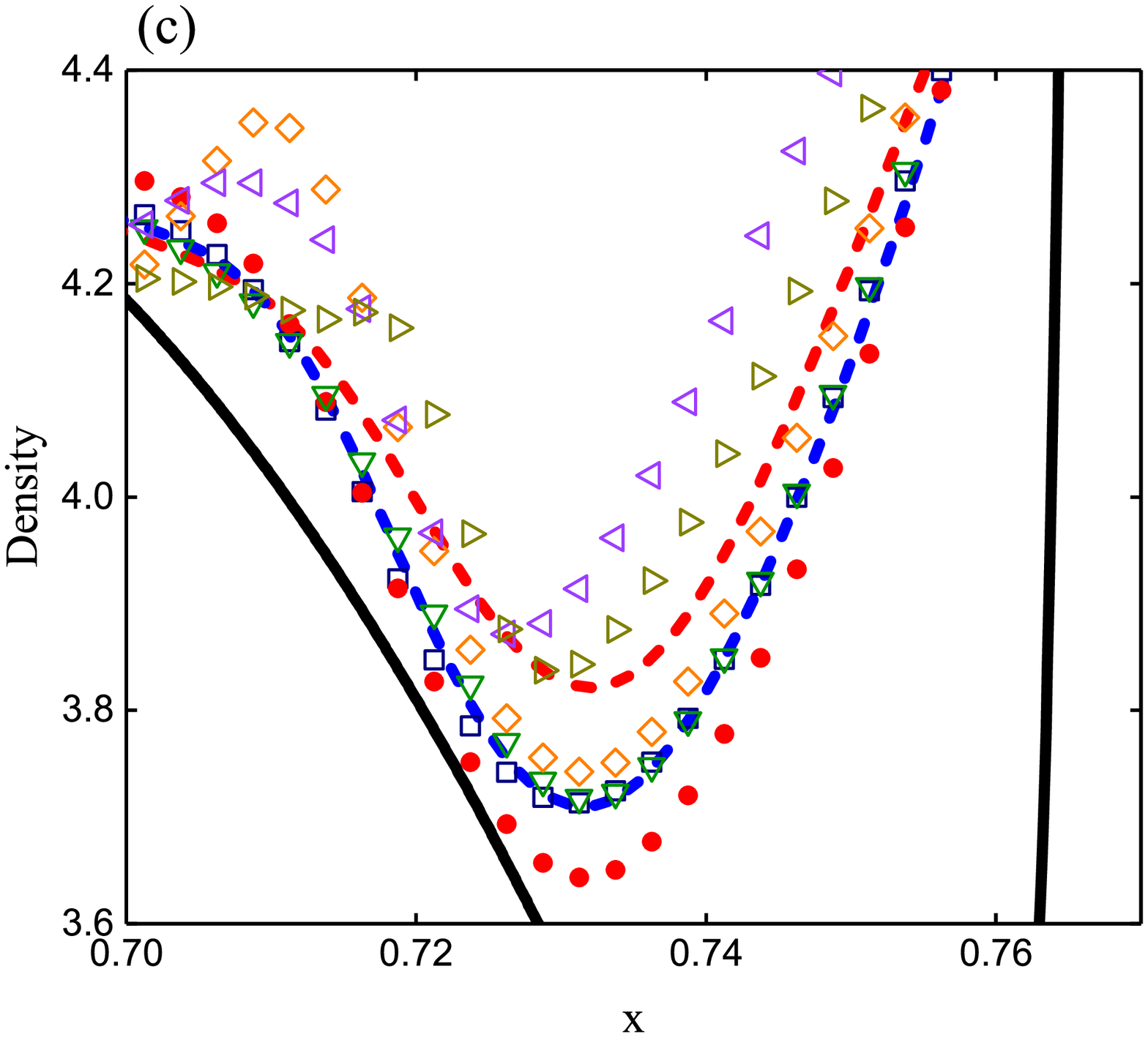}
  \includegraphics[height=0.39\textwidth]
  {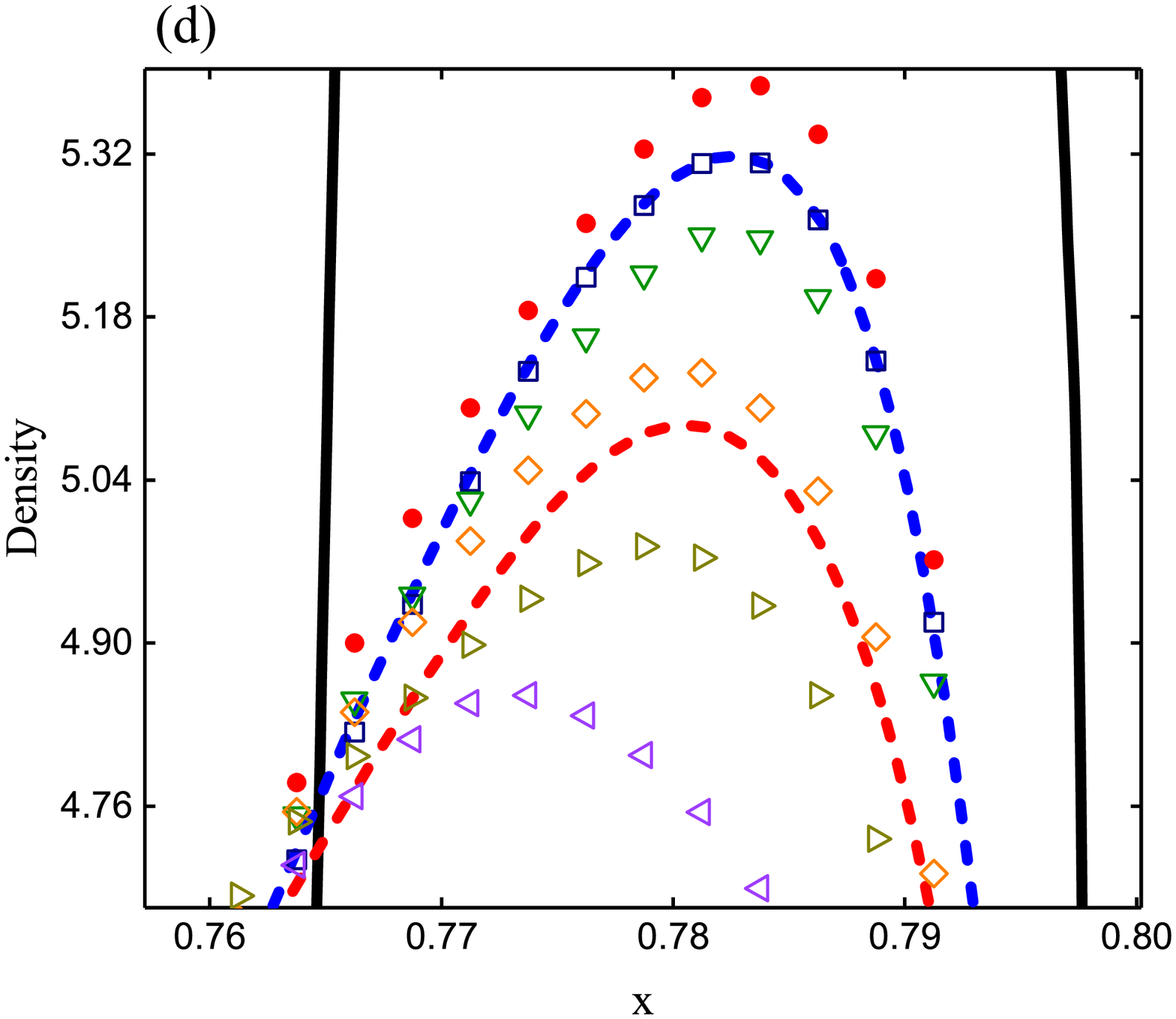}\\  
\caption{The density profiles of the 
Woodward-Colella interacting blastwaves.}
\label{fig:ex:Woodward-Colella-Problem}
\end{figure}

In order to test the robustness of the WENO-ACM scheme as 
discussed in subsection \ref{subsec:par_study}, we test this 
blastwave problem by using the WENO-ACM scheme with more different 
$\mathrm{CFS}_{s}$, that is, $\mathrm{CFS}_{s} = 0.001d_{s}, 
0.01d_{s}, 0.095d_{s}, 0.099d_{s}$, $0.0999d_{s}, 0.1d_{s}, 0.3d_{s}
, 0.5d_{s}, 0.7d_{s}, 0.9d_{s}$, as well as the considered WENO 
schemes used in previous examples. Fig. 
\ref{fig:ex:Woodward-Colella-Problem} shows the density profiles at
output time $t=0.038$ with the uniform cell number $N = 400$. The
reference solution is calculated by using the WENO-JS scheme with
$N = 10000$.

As expected, when $\mathrm{CFS}_{s}$ gets too small, 
like $\mathrm{CFS}_{s} < 0.1d_{s}$, the solutions have blown-up. In 
other words, the effect of the parameter $\mathrm{CFS}_{s}$ in the 
WENO-ACM scheme for solving time-dependent PDEs, especially for the 
robustness of the WENO-ACM scheme, is nonnegligible. In addition, 
for the computing cases when $\mathrm{CFS}_{s} \geq 0.1d_{s}$, the 
solutions have not blown-up. From Fig.
\ref{fig:ex:Woodward-Colella-Problem}, we can see these 
phenomena very clear that: (1) when $\mathrm{CFS}_{s}$ gets larger 
and $\mathrm{CFS}_{s} \geq 0.1d_{s}$, 
the WENO-ACM scheme gives solution with more numerical dissipation 
and lower resolution; (2) when $\mathrm{CFS}_{s} = 0.3d_{s}$, the 
WENO-ACM scheme gives solution close to that of the WENO-PM6 scheme 
whose solution is comparable with that of the WENO-M scheme; (3) 
when $\mathrm{CFS}_{s} = 0.5d_{s}$, the WENO-ACM scheme gives 
solution close to that of the WENO-JS scheme; (4) when 
$\mathrm{CFS}_{s} = 0.7d_{s}, 0.9d_{s}$, the WENO-ACM scheme gives 
solutions with lower resolution than that of the WENO-JS scheme 
whose solution shows lowest resolution among those of the WENO-JS, 
WENO-M and WENO-PM6 schemes; (5) when $\mathrm{CFS}_{s} = 0.1d_{s}$, 
the WENO-ACM scheme gives solution with the highest resolution among 
those of all considered schemes, while the solutions blow up when 
$\mathrm{CFS}_{s} < 0.1d_{s}$.

\subsection{Two-dimensional Euler system}
\label{subsec:examples_2D_Euler}
In this subsection, we calculate the two-dimensional Euler system
of gas dynamics with different initial and boundary conditions. The 
two-dimensional Euler system is given by the following strong conservation form of mass, momentum and energy
\begin{equation}
\begin{array}{ll}
\begin{aligned}
&\dfrac{\partial \rho}{\partial t} + \dfrac{\partial (\rho u)}
{\partial x} + \dfrac{\partial (\rho v)}{\partial y} = 0, \\
&\dfrac{\partial (\rho u)}{\partial t} + \dfrac{\partial (\rho u^{2}
+ p)}{\partial x} + \dfrac{\partial (\rho uv)}{\partial y}= 0, \\
&\dfrac{\partial (\rho v)}{\partial t} + \dfrac{\partial (\rho vu)}
{\partial x} + \dfrac{\partial (\rho v^{2} + p)}{\partial y}= 0, \\
&\dfrac{\partial E}{\partial t} + \dfrac{\partial (uE + up)}
{\partial x} + \dfrac{\partial (vE + vp)}{\partial y} = 0, \\
\end{aligned}
\end{array}
\label{2DEulerEquations}
\end{equation}
where $v$ represents the velocity component in the $y-$direction, and
the other variables are the same as in Eq.(\ref{1DEulerEquations}).
The equation of state for an ideal polytropic gas used to close the 
two-dimensional Euler system Eq.(\ref{2DEulerEquations}) is given by
\begin{equation*}
p = (\gamma - 1)\Big( E - \dfrac{1}{2}\rho(u^{2} + v^{2})\Big).
\end{equation*}
In all numerical examples of this subsection, the CFL number is set 
to be $0.5$.

\begin{example}
\bf{(Shock-vortex interaction)}
\rm{In the shock-vortex interaction problem
\cite{Shock-vortex_interaction-1,Shock-vortex_interaction-2,
Shock-vortex_interaction-3}, a left-moving shock wave interacts with
a right-moving vortex. The initial condition is set over the 
computational domain $[0, 1]\times[0,1]$ by}
\label{ex:shock-vortex}
\end{example}
\begin{equation*}
\big( \rho, u, v, p \big)(x, y, 0) = \left\{
\begin{aligned}
\begin{array}{ll}
(\rho_{\mathrm{L}}, u_{\mathrm{L}}, v_{\mathrm{L}}, p_{\mathrm{L}}),
& x < 0.5, \\
(\rho_{\mathrm{R}}, u_{\mathrm{R}}, v_{\mathrm{R}}, p_{\mathrm{R}}),
& x \geq 0.5, \\
\end{array}
\end{aligned}
\right.
\label{eq:initial_Euer2D:shock-vortex-interaction}
\end{equation*}
where the left state is taken as $(\rho_{\mathrm{L}}, u_{\mathrm{L}}
, v_{\mathrm{L}}, p_{\mathrm{L}}) = (1, \sqrt{\gamma}, 0, 1)$, and
the right state is given as
\begin{equation*}
\begin{array}{l}
p_{\mathrm{R}} = 1.3, \rho_{\mathrm{R}} = \rho_{\mathrm{L}}\bigg( 
\dfrac{\gamma - 1 + (\gamma + 1)p_{\mathrm{R}}}{\gamma + 1 + (\gamma
- 1)p_{\mathrm{R}}} \bigg)\\
u_{\mathrm{R}} = u_{\mathrm{L}}\bigg( \dfrac{1 - p_{\mathrm{R}}}{
\sqrt{\gamma - 1 + p_{\mathrm{R}}(\gamma + 1)}}\bigg), v_{\mathrm{R}}
= 0.
\end{array}
\label{eq:Euler2D:shock-vortex-interaction:rightState}
\end{equation*}
A vortex given by perturbations $(\delta \rho, \delta u, \delta v, 
\delta p)$ is superimposed onto the state when $x < 0.5$, and the
perturbations are defined by
\begin{equation*}
\delta \rho = \dfrac{\rho_{\mathrm{L}}^{2}}{(\gamma - 1)
p_{\mathrm{L}}}\delta T, 
\delta u = \epsilon \dfrac{y - y_{\mathrm{c}}}{r_\mathrm{c}}
\mathrm{e}^{\alpha(1-r^{2})}, 
\delta v = - \epsilon \dfrac{x - x_{\mathrm{c}}}{r_\mathrm{c}}
\mathrm{e}^{\alpha(1-r^{2})}, 
\delta p = \dfrac{\gamma \rho_{\mathrm{L}}^{2}}{(\gamma - 1)
\rho_{\mathrm{L}}}\delta T,
\label{eq:Euler2D:shock-vortex-interactions:Perturbations}
\end{equation*}
where $\epsilon = 0.3, r_{\mathrm{c}} = 0.05, \alpha = 0.204, 
x_{\mathrm{c}} = 0.25, y_{\mathrm{c}} = 0.5$ and
\begin{equation*}
r = \sqrt{\dfrac{(x - x_{\mathrm{c}})^{2} + (y - y_{\mathrm{c}})^{2}}
{r_{\mathrm{c}}^{2}}}, \\
\delta T = - \dfrac{\gamma - 1}{4\alpha \gamma}\epsilon^{2}
\mathrm{e}^{2\alpha (1 - r^{2})}.
\label{eq:Euler2D:shock-vortex-interaction:deltaT}
\end{equation*}

\begin{figure}[ht]
\centering
  \includegraphics[height=0.39\textwidth]
  {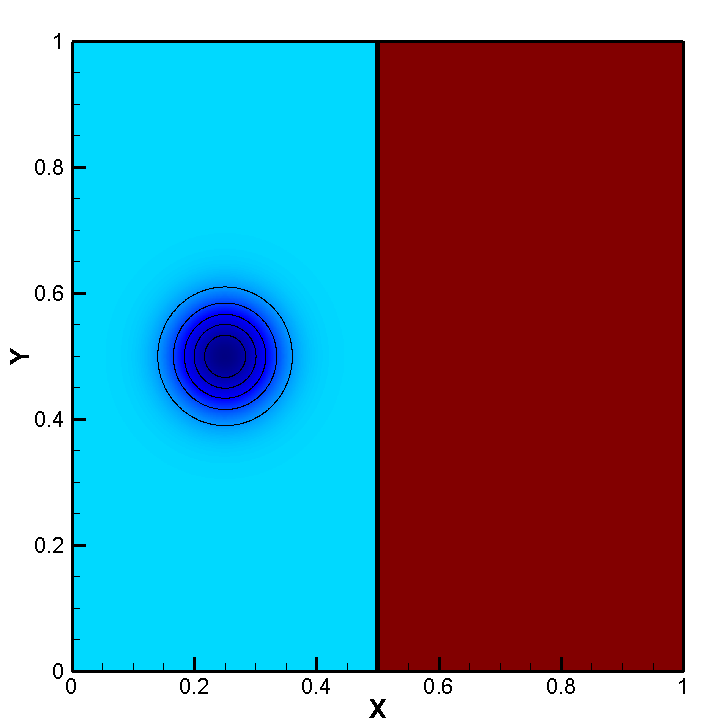}
  \includegraphics[height=0.39\textwidth]
  {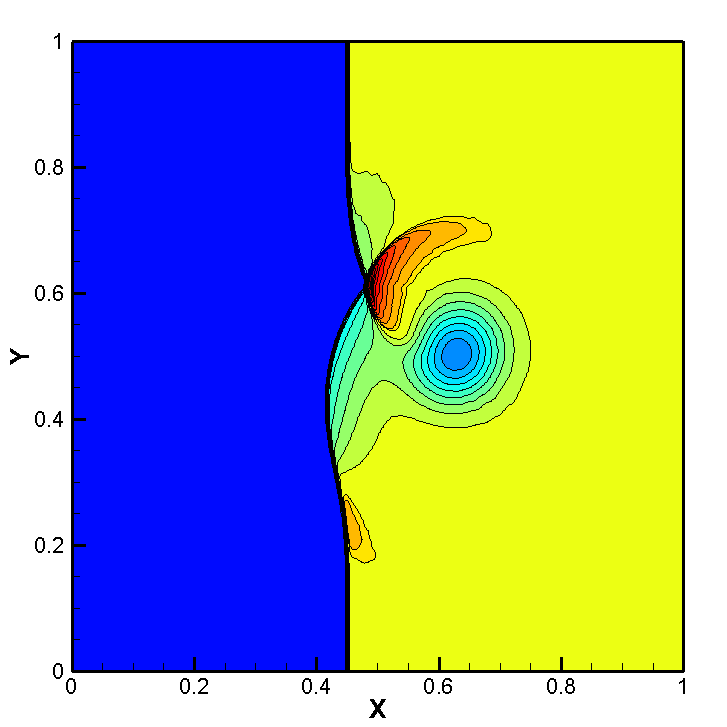}
\caption{Density plots for the shock-vortex interaction problem using
$30$ contour lines with range from $0.9$ to $1.4$, computed using the
WENO-ACM scheme at output time $t = 0.35$ with a uniform mesh size 
of $400\times400$. (left: initial density, right: final density)}
\label{fig:ex:SVI}
\end{figure}

\begin{figure}[ht]
\centering
  \includegraphics[height=0.38\textwidth]
  {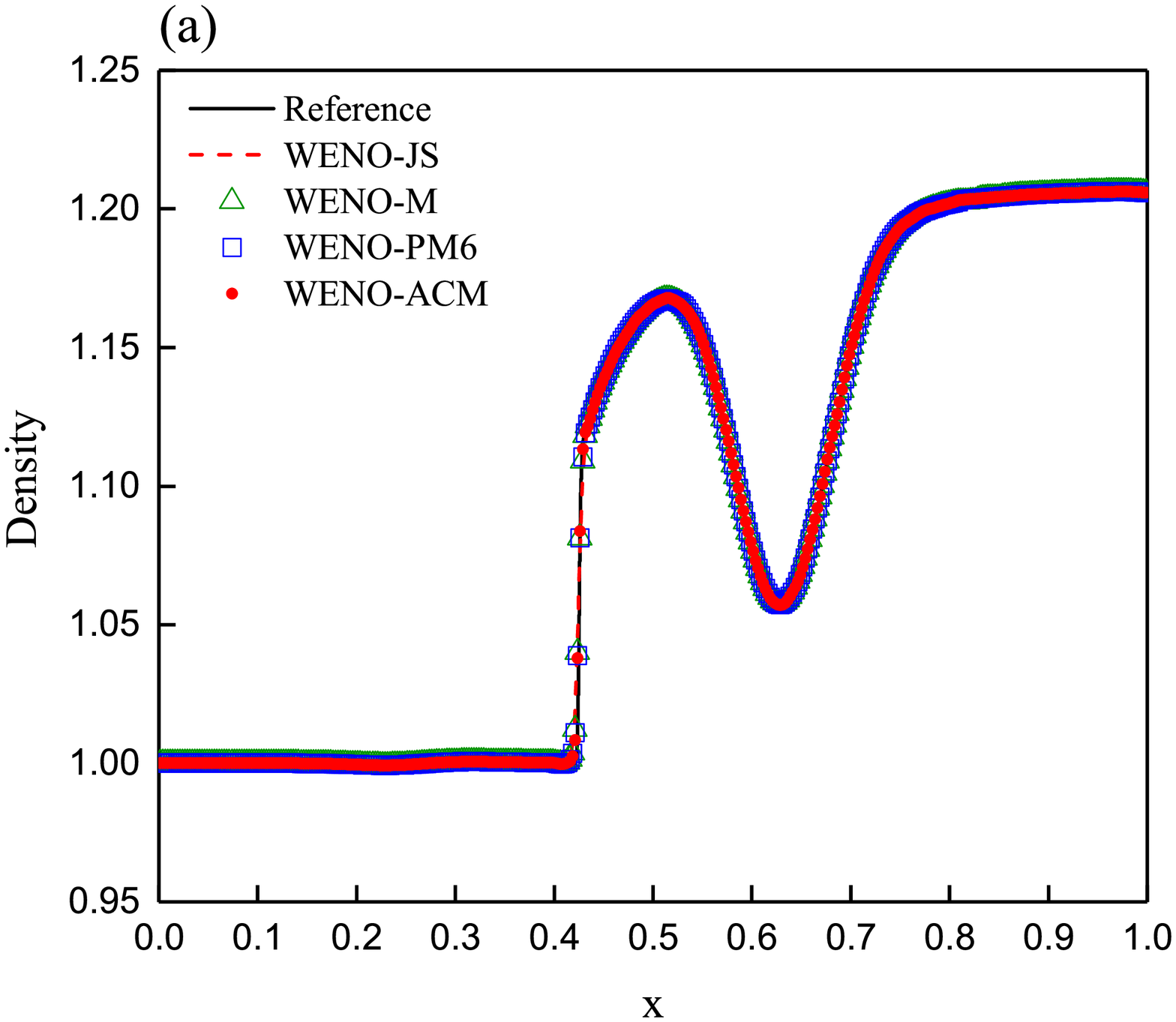}
  \includegraphics[height=0.38\textwidth]
  {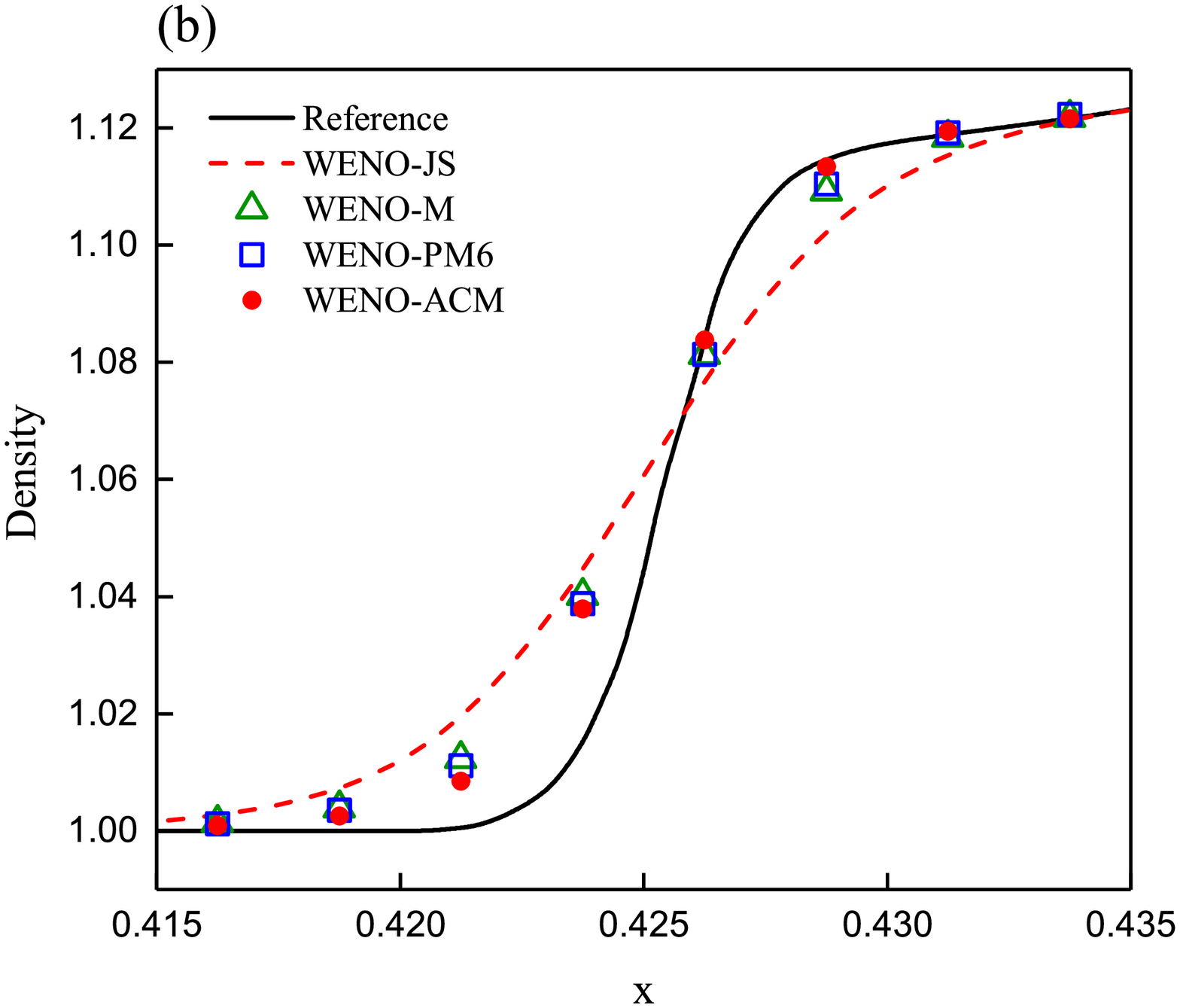}\\
  \includegraphics[height=0.38\textwidth]
  {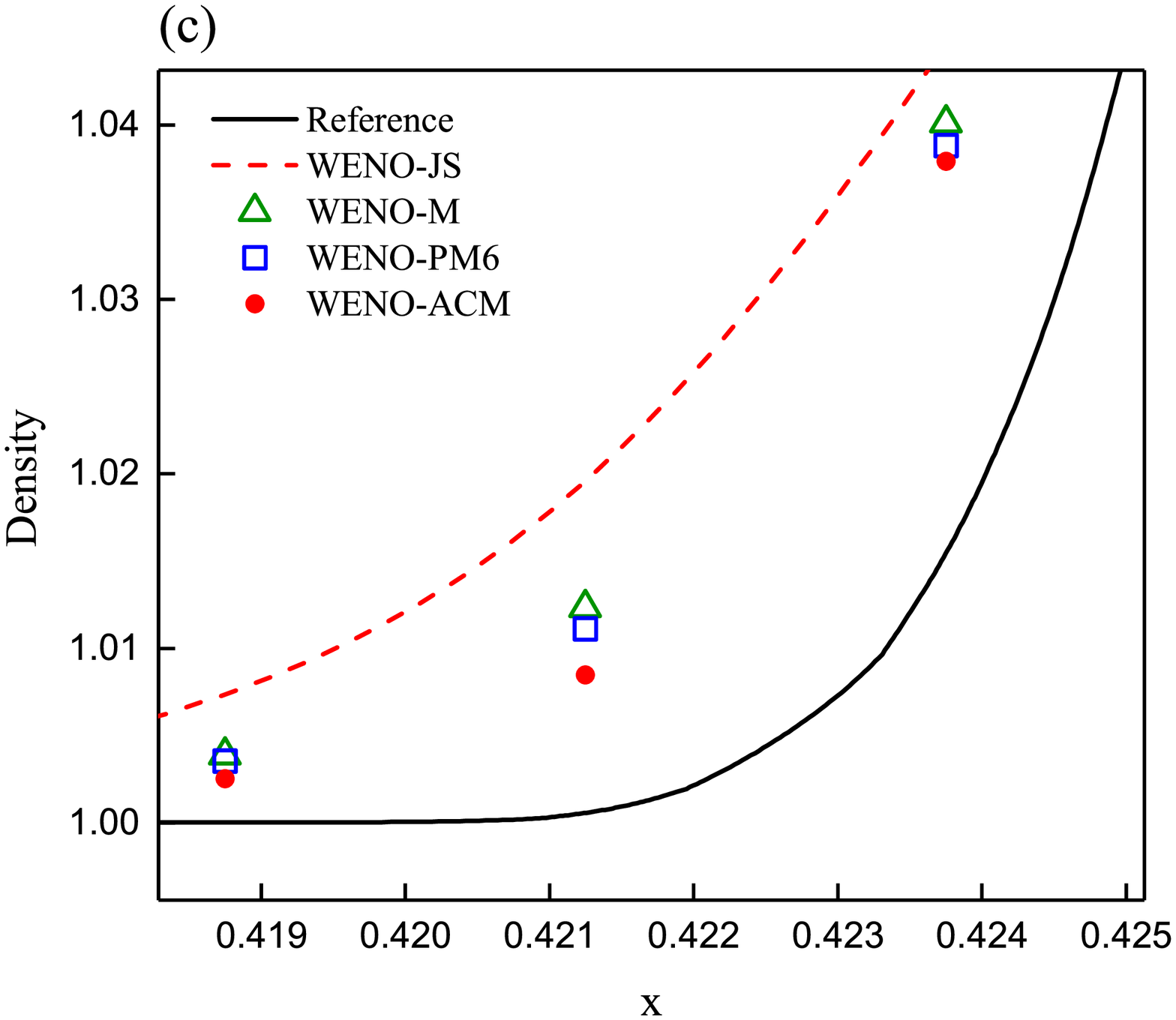}
  \includegraphics[height=0.38\textwidth]
  {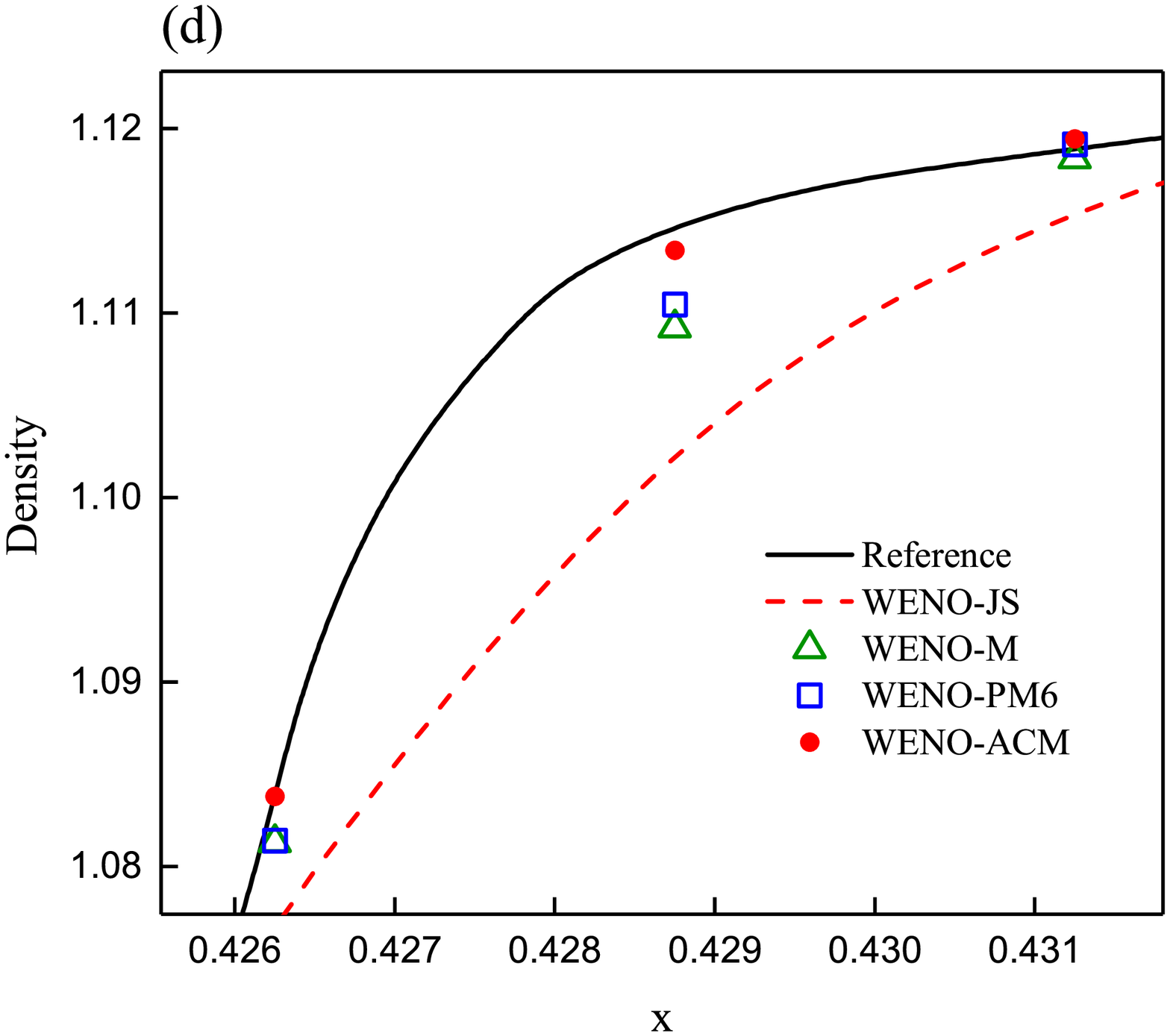}
\caption{The cross-sectional slices of density plot along the plane
$y = 0.5$, computed using the WENO-JS, WENO-M, WENO-PM6 and WENO-ACM
schemes at output time $t = 0.35$ with a uniform mesh size of
$400\times400$.}
\label{fig:ex:SVI-Y05}
\end{figure}

We discretize the computational domain with a uniform mesh size of 
$400 \times 400$. The transmissive boundary conditions are used and 
the output time is taken as $t = 0.35$. The initial and final
positions of the shock and vortex in density profile, computed using
the WENO-ACM scheme, have been shown in Fig. \ref{fig:ex:SVI}. We 
can easily observe that the WENO-ACM scheme performs very well in 
capturing the complex structure of the shock and vortex after the 
interaction. In Fig. \ref{fig:ex:SVI-Y05}, we have presented the 
cross-sectional slices of density profile along the plane $y=0.5$, 
computed by the WENO-JS, WENO-M, WENO-PM6 and WENO-ACM 
schemes at output time $t = 0.35$. The reference solution is 
obtained using the WENO-JS scheme with a uniform mesh size of 
$1000 \times 1000$. As shown in the zoomed-in plots around the 
shock, the WENO-ACM scheme resolves the shock in a non-oscillatory 
manner and provides a better resolution than the other considered 
schemes. The WENO-PM6 and WENO-M schemes perform better than the 
WENO-JS scheme, while the WENO-PM6 scheme performs slightly better 
than the WENO-M scheme.

\begin{example}
\bf{(Explosion problem)} 
\rm{We solve the explosion problem \cite{Toro_RiemannSolvers} 
modeled by the two-dimensional Euler equations 
Eq.(\ref{2DEulerEquations}) on the square domain $[-1, 1] \times 
[-1, 1]$ in the $x-y$ plane. It involves two constant states of flow 
variables separated with a circle of radius $R = 0.4$ centered at
$(0, 0)$. The initial condition is given as}
\label{ex:explosion}
\end{example}
\begin{equation*}
\big( \rho, u, v, p \big)(x, y, 0) = \left\{
\begin{array}{ll}
(1, 0, 0, 1), & \mathrm{if} \sqrt{x^{2} + y^{2}} < 0.4,\\
(0.125, 0, 0, 0.1), & \mathrm{else}.
\end{array}\right.
\label{eq:initial_Euer2D:Explosion}
\end{equation*}
Transmissive boundary conditions are used on all boundaries. The
computational domain is discretized with a uniform mesh size of 
$400\times400$ and the final time is chosen to be $t=0.25$. 

\begin{figure}[ht]
\centering
  \includegraphics[height=0.3\textwidth]
  {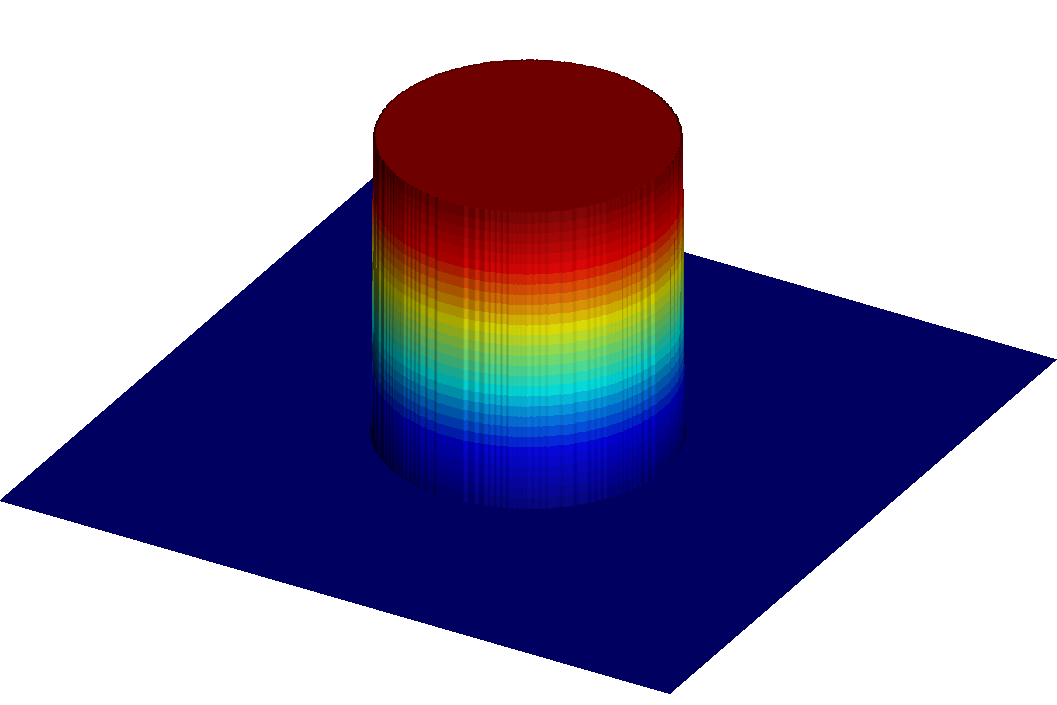}
  \includegraphics[height=0.3\textwidth]
  {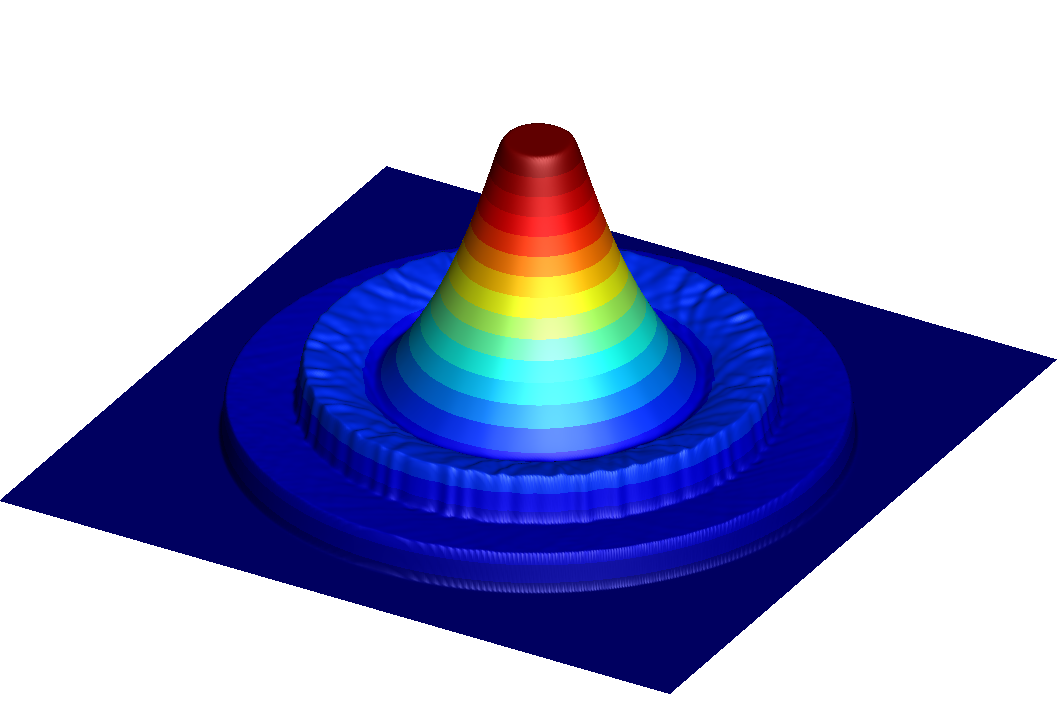}
\caption{Density plots for the explosion problem with range from
$0.125$ to $1.0$, computed using the WENO-ACM scheme at output time
$t = 0.25$ with a uniform mesh size of $400\times400$. (left: initial
density, right: final density)}
\label{fig:ex:Explosion}
\end{figure}

\begin{figure}[ht]
\centering
  \includegraphics[height=0.38\textwidth]
  {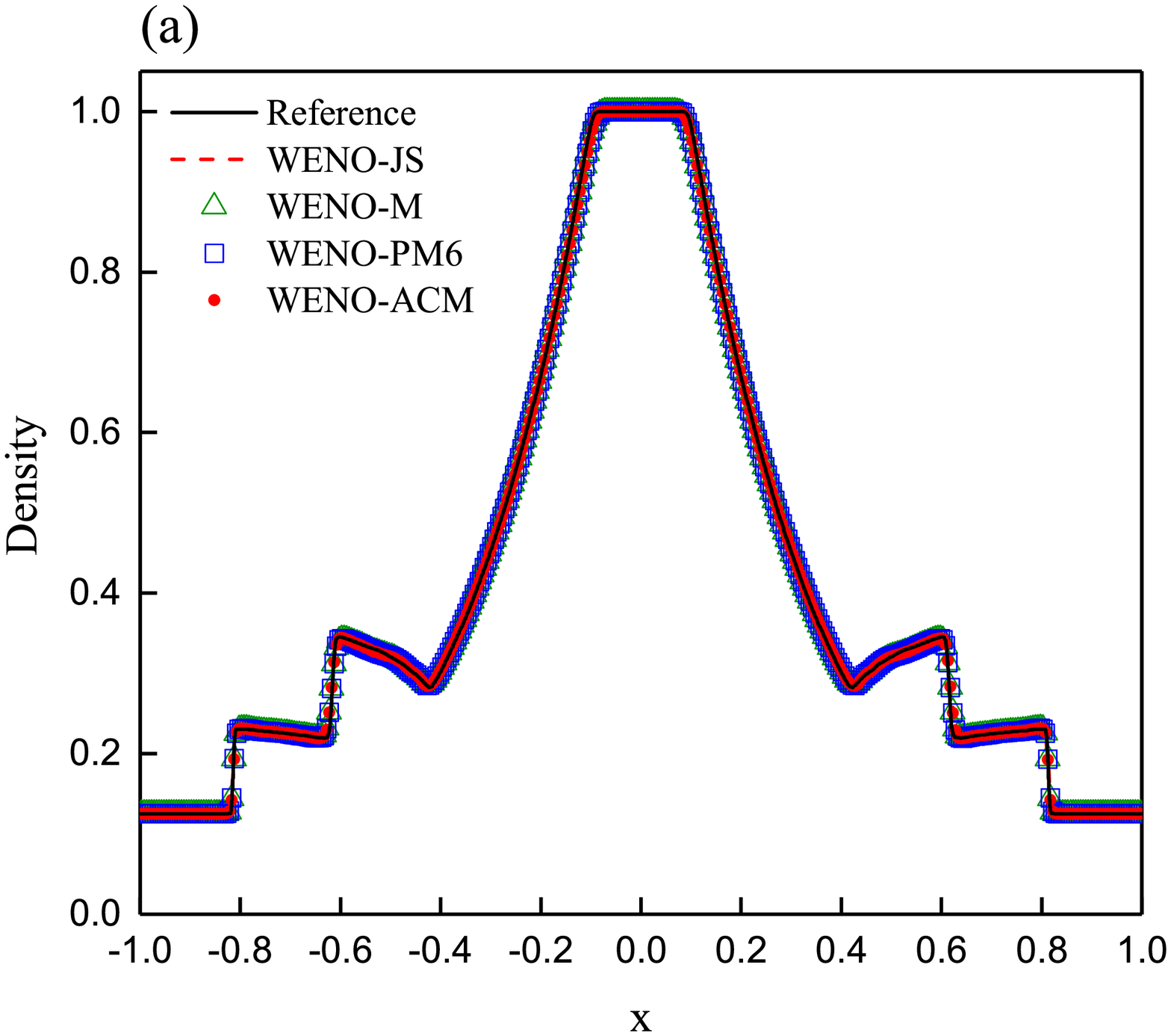}
  \includegraphics[height=0.38\textwidth]
  {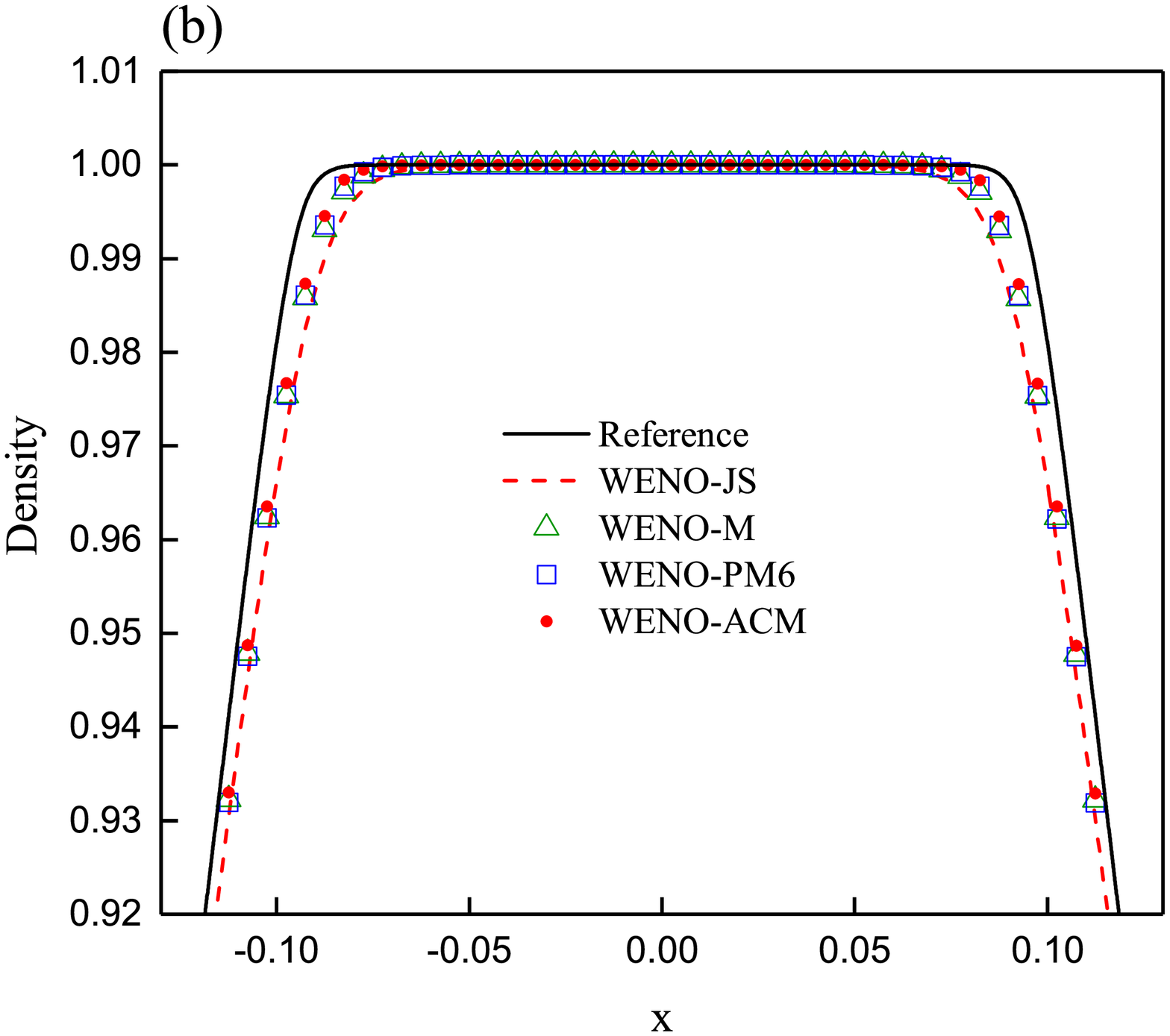}\\
  \includegraphics[height=0.38\textwidth]
  {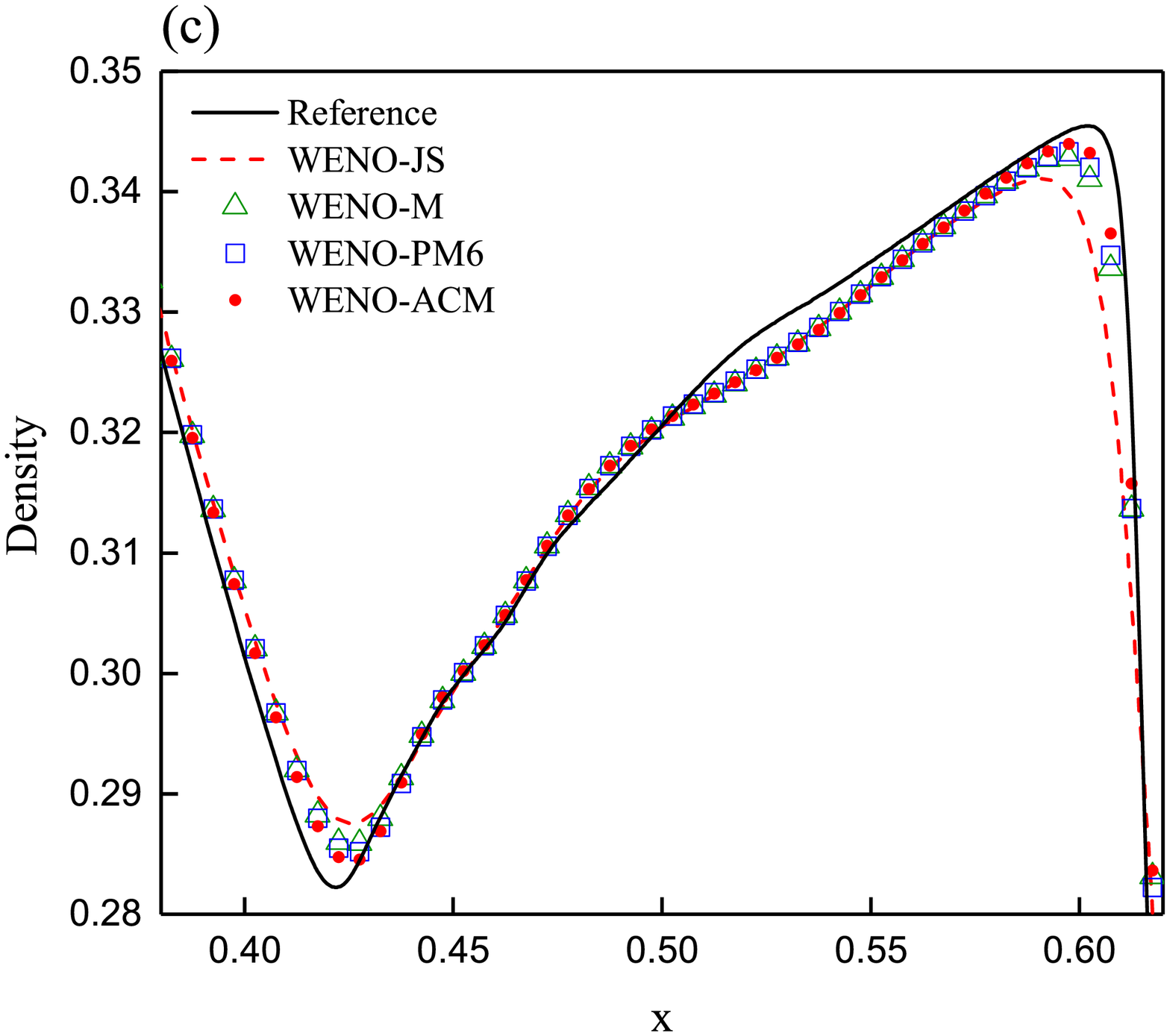}
  \includegraphics[height=0.38\textwidth]
  {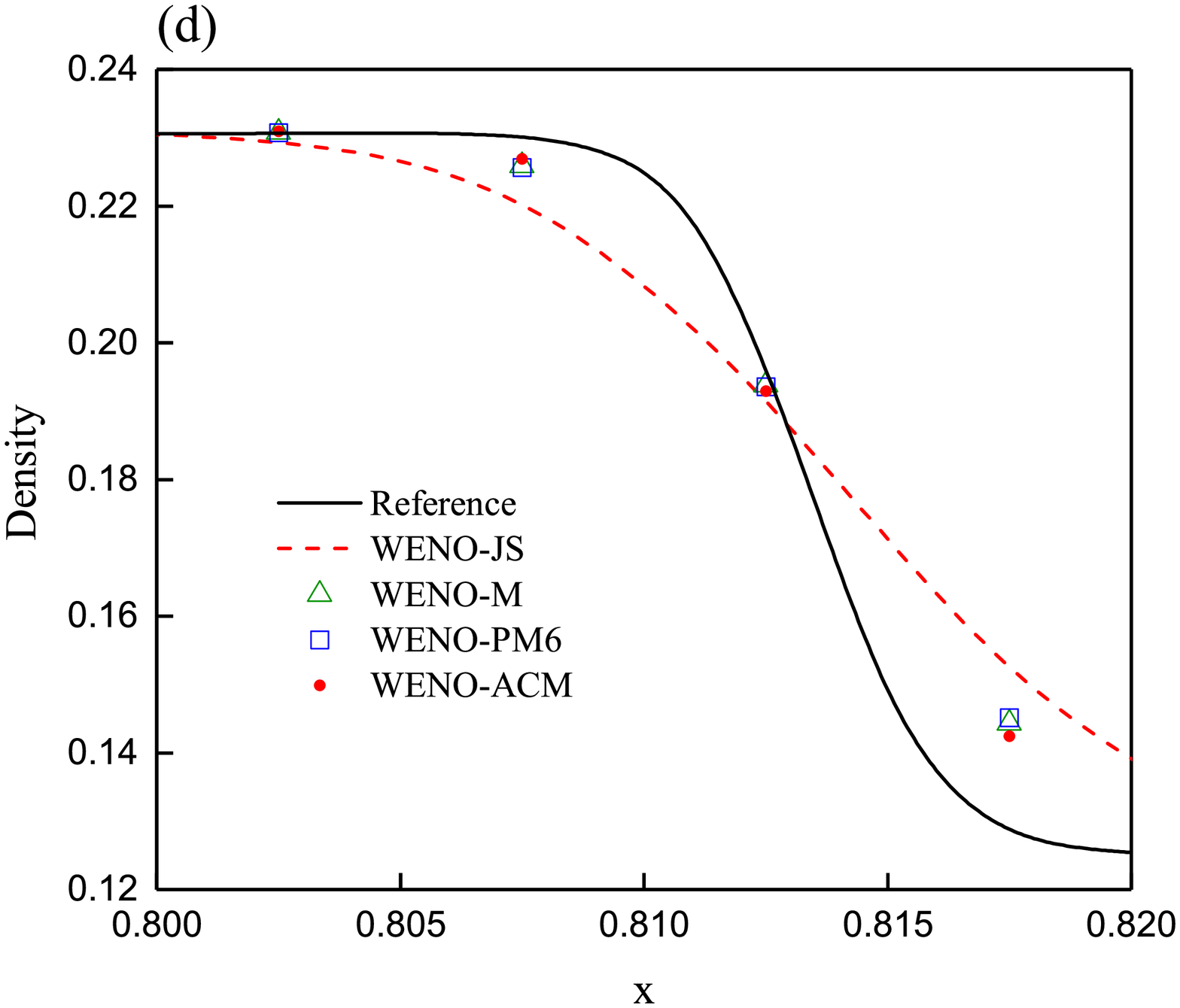}
\caption{The cross-sectional slices of density plot along the plane
$y = 0.0$, computed using the WENO-JS, WENO-M, WENO-PM6 and WENO-ACM
schemes at output time $t = 0.25$ with a uniform mesh size of
$400\times400$.}
\label{fig:ex:Explosion-Y0}
\end{figure}

The density contours of the initial and final states, computed using 
the WENO-ACM scheme, have been shown in Fig. \ref{fig:ex:Explosion}. 
The WENO-ACM scheme is able to capture the structure of the 
explosion problem. In Fig.\ref{fig:ex:Explosion-Y0}, we have 
presented the cross-sectional slices of density profile along the 
plane $y=0.0$, calculated by the WENO-JS, WENO-M, WENO-PM6 and 
WENO-ACM schemes. The reference solution is obtained using the 
WENO-JS scheme with a uniform mesh size of $1000 \times 1000$. We 
can observe that the WENO-M, WENO-PM6 and WENO-ACM schemes provide
significantly better resolutions than the WENO-JS scheme. In Fig. 
\ref{fig:ex:Explosion-Y0} (b) and (d), at $x \in (-0.15, 0.15)$ and
$x \in (0.8, 0.82)$, the WENO-ACM scheme presents a higher 
resolution than the WENO-PM6 and WENO-M schemes whose resolutions 
are comparable. In Fig. \ref{fig:ex:Explosion-Y0} (c), the WENO-ACM 
scheme gives better resolution at $x \in (0.38, 0.62)$ than the 
WENO-PM6 scheme which performs slightly better than the WENO-M 
scheme.

\begin{example}
\bf{(2D Riemann problem)} 
\rm{The calculation of the 2D Riemann problem \cite{Riemann-2D-01,
Riemann2D-02} is done over a unit square domain $[0,1] \times[0,1]$, 
initially involves the constant states of flow variables over each 
quadrant which is got by dividing the computational domain using 
lines $x = x_{0}$ and $y = y_{0}$. There are many different
configurations for the 2D Riemann problem \cite{Riemann2D-02}. In 
this example, the configuration is taken from 
\cite{WENO-PPM5,Riemann2D-02} with the following initial data}
\label{ex:Riemann2D}
\end{example}
\begin{equation*}
\big( \rho, u, v, p \big)(x, y, 0) = \left\{
\begin{aligned}
\begin{array}{ll}
(1.0, 0.0, -0.3, 1.0),       & 0.5 \leq x \leq 1.0, 
                               0.5 \leq y \leq 1.0, \\
(2.0, 0.0, 0.3, 1.0),        & 0.0 \leq x \leq 0.5, 
                               0.5 \leq y \leq 1.0, \\
(1.0625, 0.0, 0.8145, 0.4),  & 0.0 \leq x \leq 0.5, 
                               0.0 \leq y \leq 0.5, \\
(0.5313, 0.0, 0.4276, 0.4),  & 0.5 \leq x \leq 1.0, 
                             0.0 \leq y \leq 0.5. \\
\end{array}
\end{aligned}
\right.
\label{eq:initial_Euer2D:Riemann2D}
\end{equation*}
The transmission boundary conditions are used on all boundaries, and
the numerical solutions are calculated using considered WENO schemes
at output time $t = 0.3$ with a mesh size of $1200 \times 1200$.

\begin{figure}[ht]
\centering
  \includegraphics[height=0.39\textwidth]
  {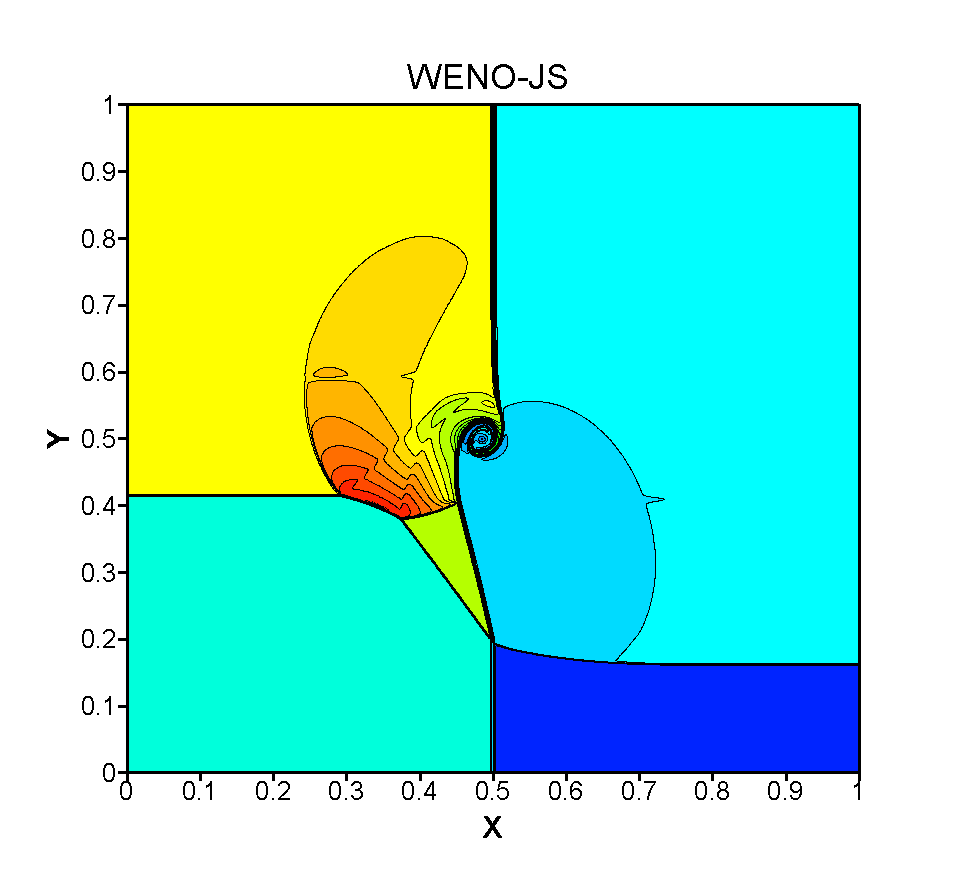}
  \includegraphics[height=0.39\textwidth]
  {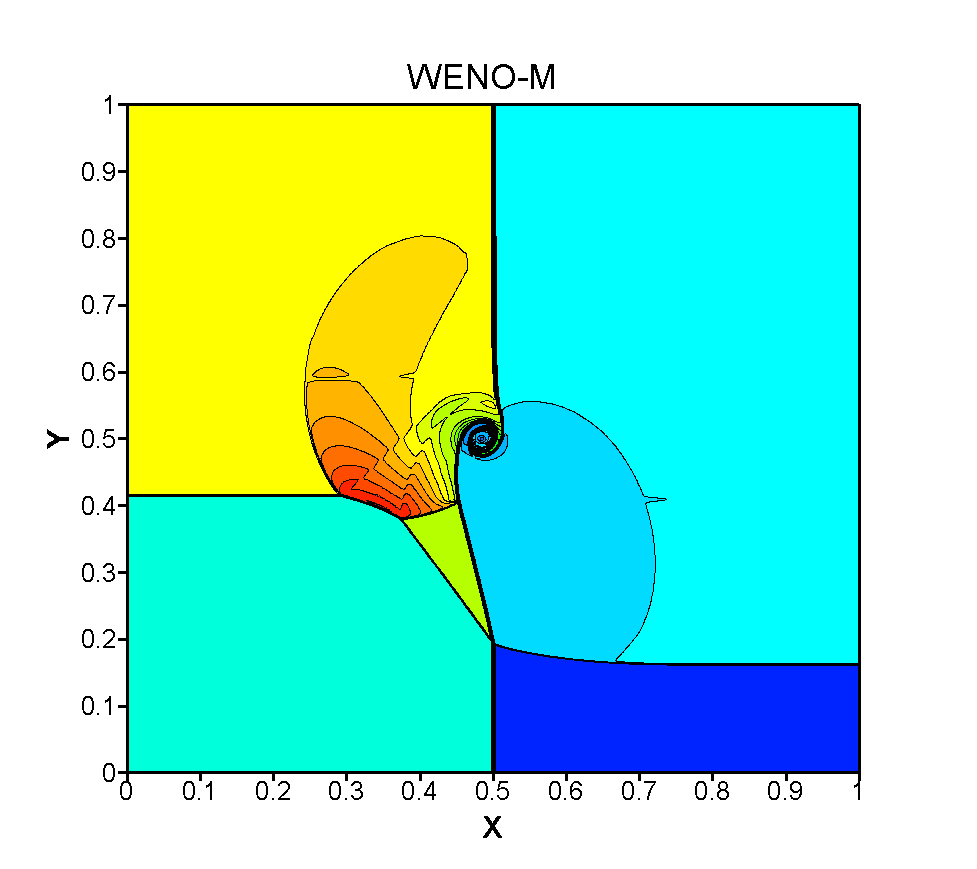}\\
  \includegraphics[height=0.39\textwidth]
  {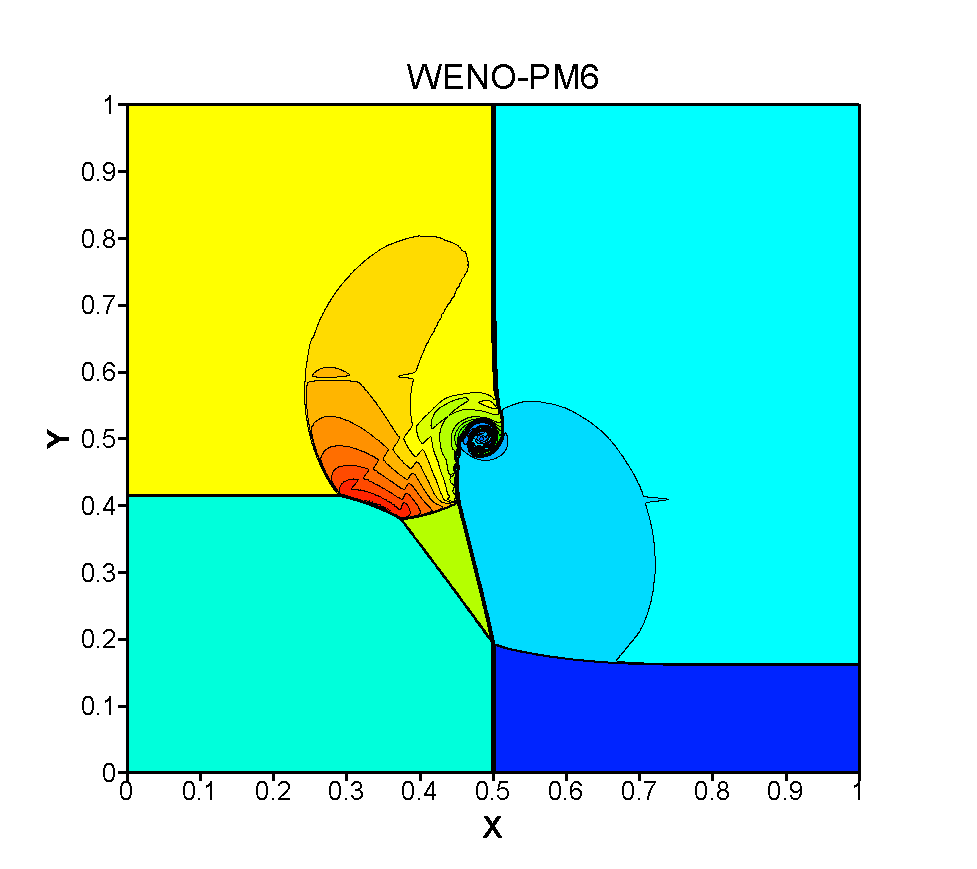}
  \includegraphics[height=0.39\textwidth]
  {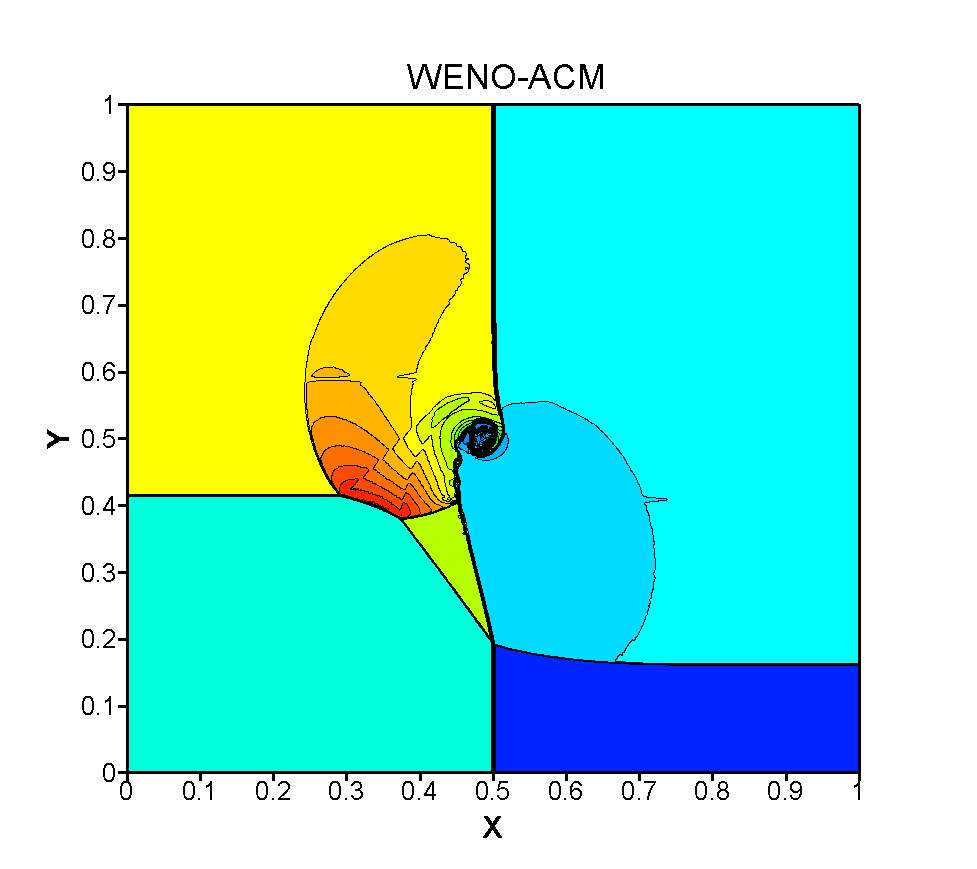}    
\caption{Density plots for the 2D Riemann problem using $30$ contour
lines with range from $0.5$ to $2.4$, computed using the WENO-JS, 
WENO-M, WENO-PM6, WENO-ACM schemes at output time $t = 0.3$ with a 
uniform mesh size of $1200\times1200$.}
\label{fig:ex:Riemann2D}
\end{figure}

\begin{figure}[ht]
\centering
  \includegraphics[height=0.39\textwidth]
  {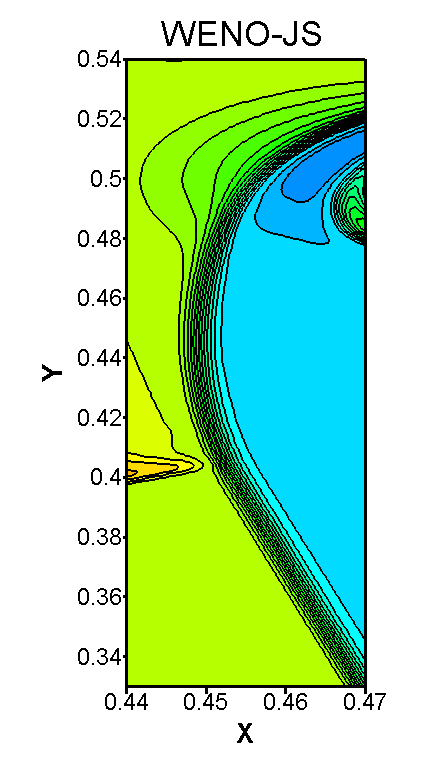}
  \includegraphics[height=0.39\textwidth]
  {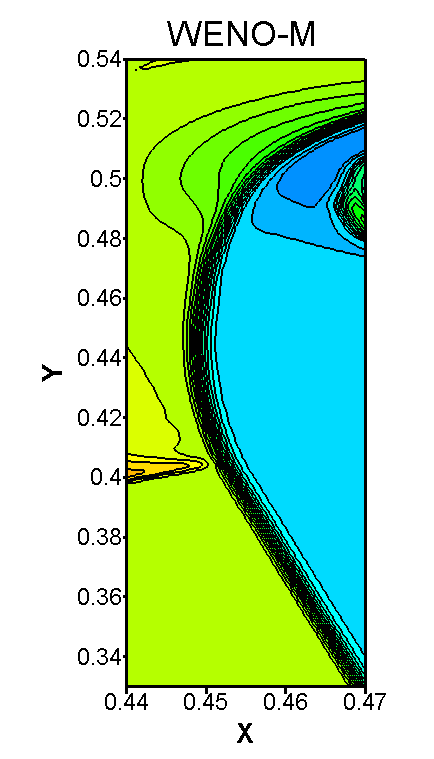}
  \includegraphics[height=0.39\textwidth]
  {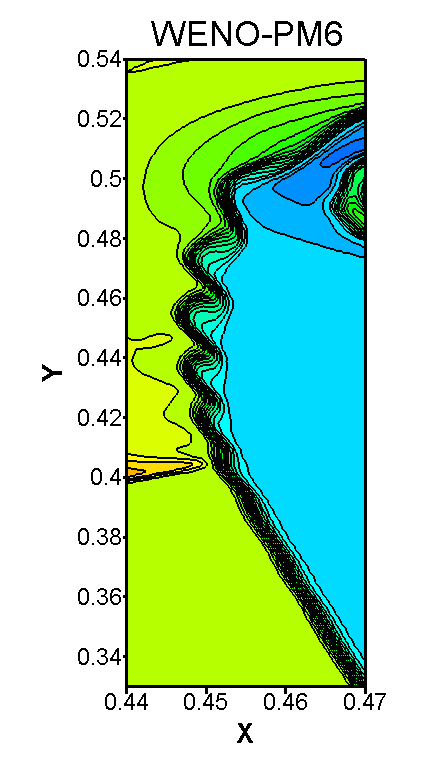}
  \includegraphics[height=0.39\textwidth]
  {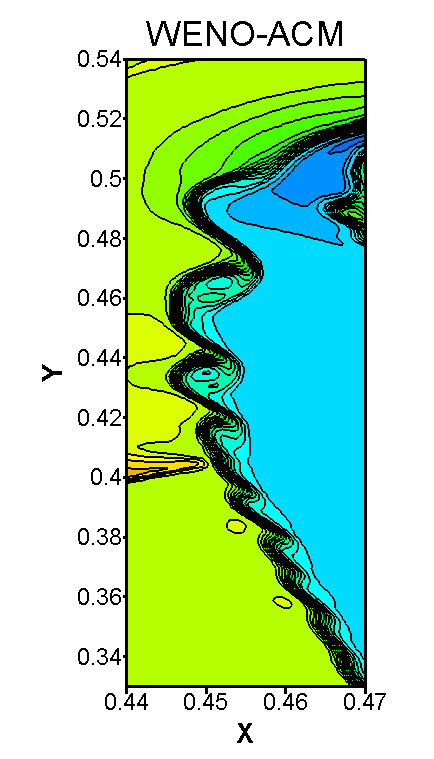}    
\caption{The zoomed-in density plots for the 2D Riemann problem, 
computed using the WENO-JS, WENO-M, WENO-PM6, WENO-ACM schemes at 
output time $t = 0.3$ with a uniform mesh size of $1200\times1200$.}
\label{fig:ex:Riemann2D-ZoomedIn}
\end{figure}

In Fig. \ref{fig:ex:Riemann2D}, we have shown the numerical results 
of density obtained using the WENO-JS, WENO-M, WENO-PM6 and WENO-ACM 
schemes. All considered schemes can capture the main structure of the
solution. However, this example is commonly focused on the 
description of instability of the slip line \cite{WENO-PPM5,
Riemann2D-04}, and we have displayed the close-up view of 
this instability in Fig. \ref{fig:ex:Riemann2D-ZoomedIn}. We can 
observe that the WENO-JS and WENO-M schemes failed to resolve the 
instability of the slip line under current spatial resolution, while 
both the WENO-PM6 and WENO-ACM schemes have evidently resolved this 
instability. The unstable physical structures in the solution of the 
WENO-ACM scheme appear to have larger length scale and bigger wave 
numbers when compared with WENO-PM6.

\begin{example}
\bf{(Double Mach reflection, DMR)}
\rm{Now we apply the considered WENO schemes to the two-dimensional 
double Mach reflection problem 
\cite{interactingBlastWaves-Woodward-Colella}. This problem is an 
important test where a vertical shock wave moves horizontally into a
wedge that is inclined by some angle. The computational domain of 
this problem is $[0,4]\times[0,1]$ and the initial condition is 
given by}
\label{ex:DMR}
\end{example}
\begin{equation*}
\big( \rho, u, v, p \big)(x, y, 0) = \left\{
\begin{aligned}
\begin{array}{ll}
(8.0, 8.25\cos \dfrac{\pi}{6}, -8.25\sin \dfrac{\pi}{6}, 116.5),         & x < x_{0} + \dfrac{y}{\sqrt{3}}, \\
(1.4, 0.0, 0.0, 1.0),    & x \geq x_{0} + \dfrac{y}{\sqrt{3}}, \\
\end{array}
\end{aligned}
\right.
\label{eq:initial_Euer2D:DMR}
\end{equation*}
where $x_{0} = \frac{1}{6}$. The left boundary condition is inflow,
with the post-shock values as stated above, and the outflow boundary 
condition is used on the right boundary. On the bottom boundary,
the reflective boundary condition is applied to the interval $[x_{0},
4]$, while at $(0, x_{0})$, the post-shock values are imposed. The
boundary condition on the upper boundary is implemented as follows
\begin{equation*}
\begin{aligned}
\begin{array}{l}
\big( \rho, u, v, p \big)(x, 1, t) = \left\{
\begin{array}{ll}
(8.0, 8.25\cos \dfrac{\pi}{6}, -8.25\sin \dfrac{\pi}{6}, 116.5), & x \in [0, s(t)), \\
(1.4, 0.0, 0.0, 1.0), & x \in [s(t), 4].
\end{array}
\right.
\end{array}
\end{aligned}
\label{eq:boundary:Euler2D:DMR:UB}
\end{equation*}
where $s(t)$ is the position of the shock wave at time $t$ on the 
upper boundary and given by $s(t) = x_{0}+\frac{1 + 20t}{\sqrt{3}}$. 
The computational domain $[0,3]\times[0,1]$ is discretized with a
uniform mesh size of $2000\times500$ and the output time is chosen 
to be $t=0.2$.

\begin{figure}[ht]
\centering
  \includegraphics[height=0.3\textwidth]
  {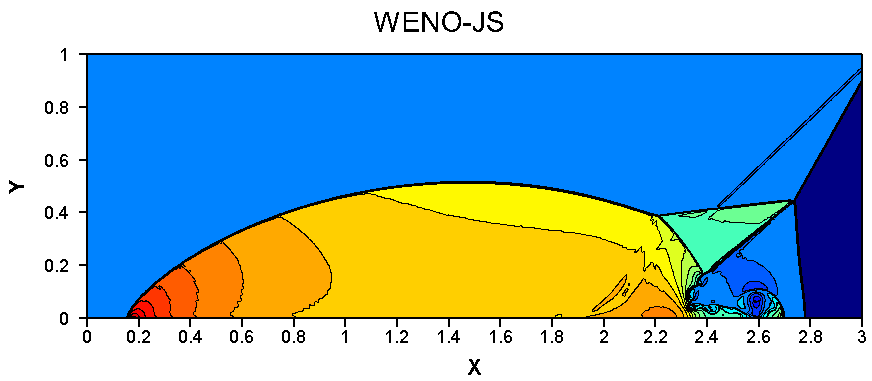}
  \includegraphics[height=0.3\textwidth]
  {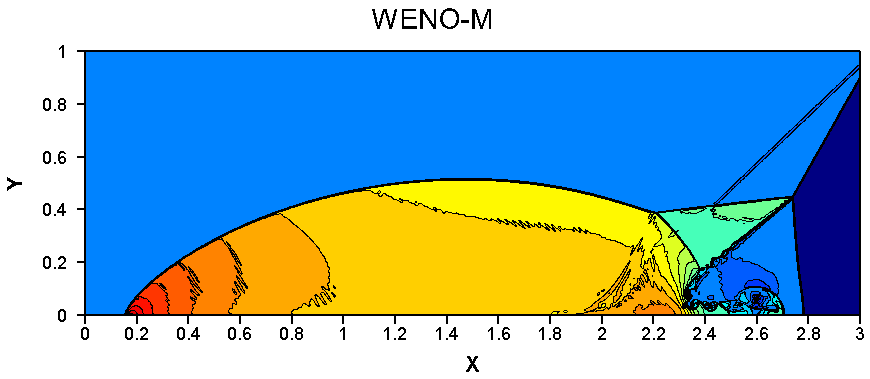}\\
  \includegraphics[height=0.3\textwidth]
  {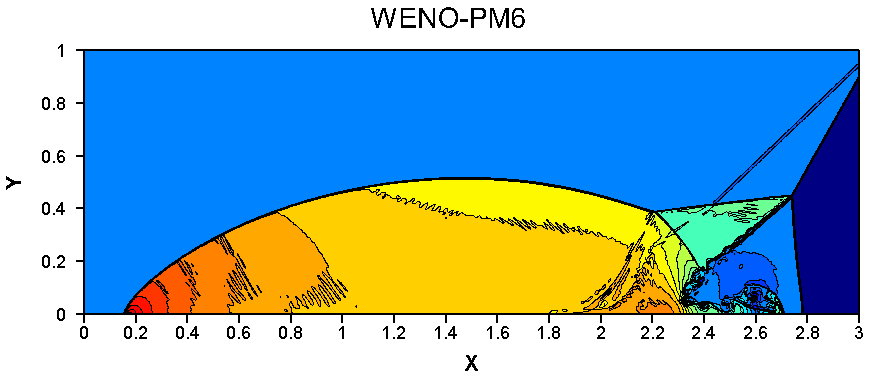}
  \includegraphics[height=0.3\textwidth]
  {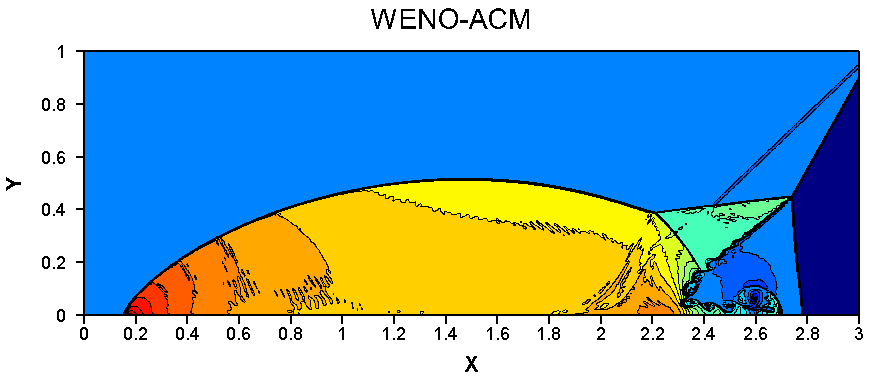}\\
\caption{Density plots for the DMR problem using 30 contour lines
with range from $1.4$ to $25$, computed using the WENO-JS, WENO-M,
WENO-PM6 and WENO-ACM schemes at output time $t = 0.2$ with a 
uniform mesh size of $2000\times500$.}
\label{fig:ex:DMR}
\end{figure}

\begin{figure}[ht]
\centering
  \includegraphics[height=0.32\textwidth]
  {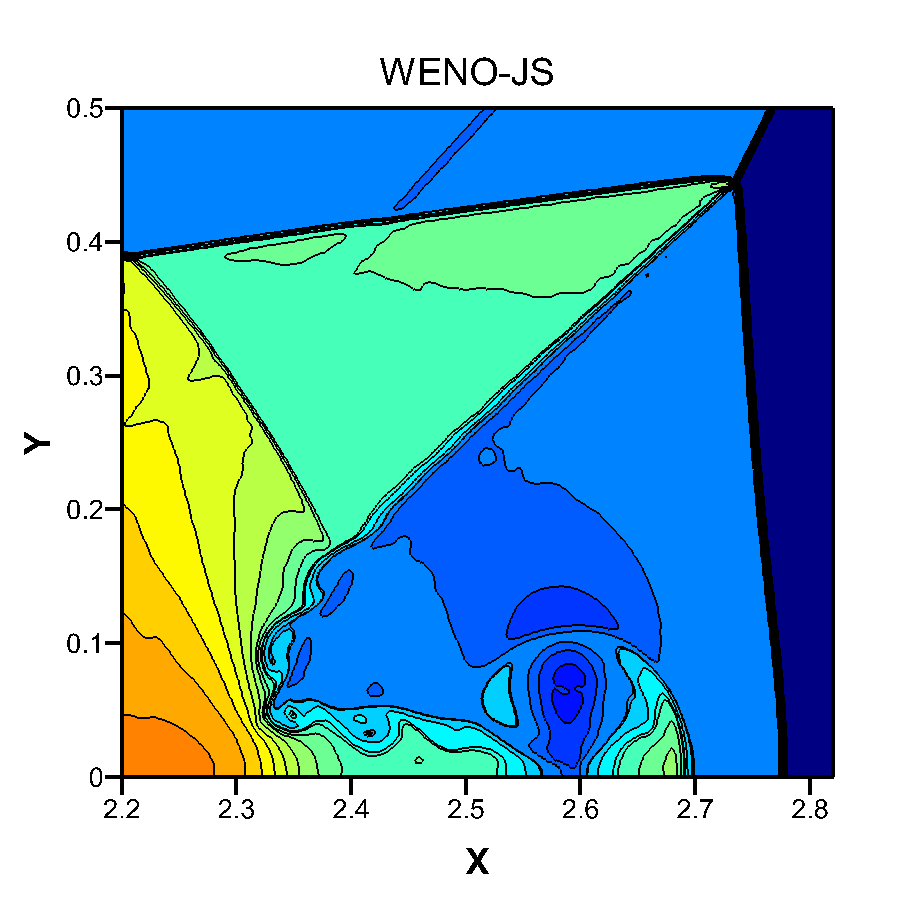}
  \includegraphics[height=0.32\textwidth]
  {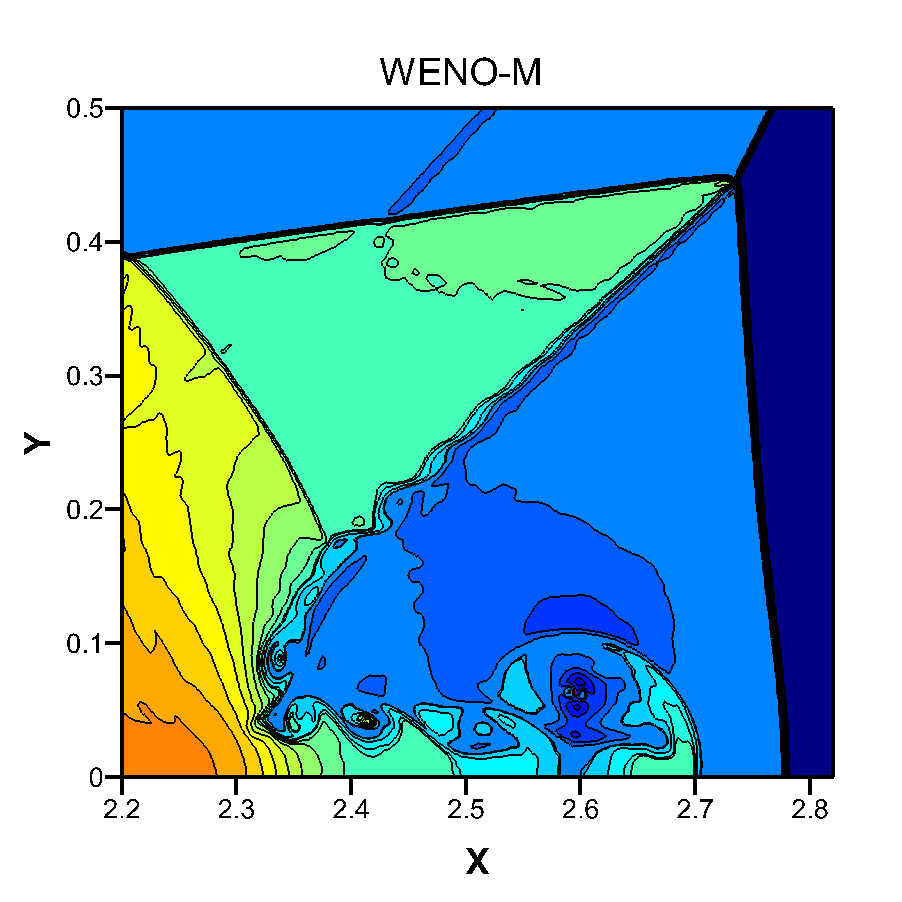}\\
  \includegraphics[height=0.32\textwidth]
  {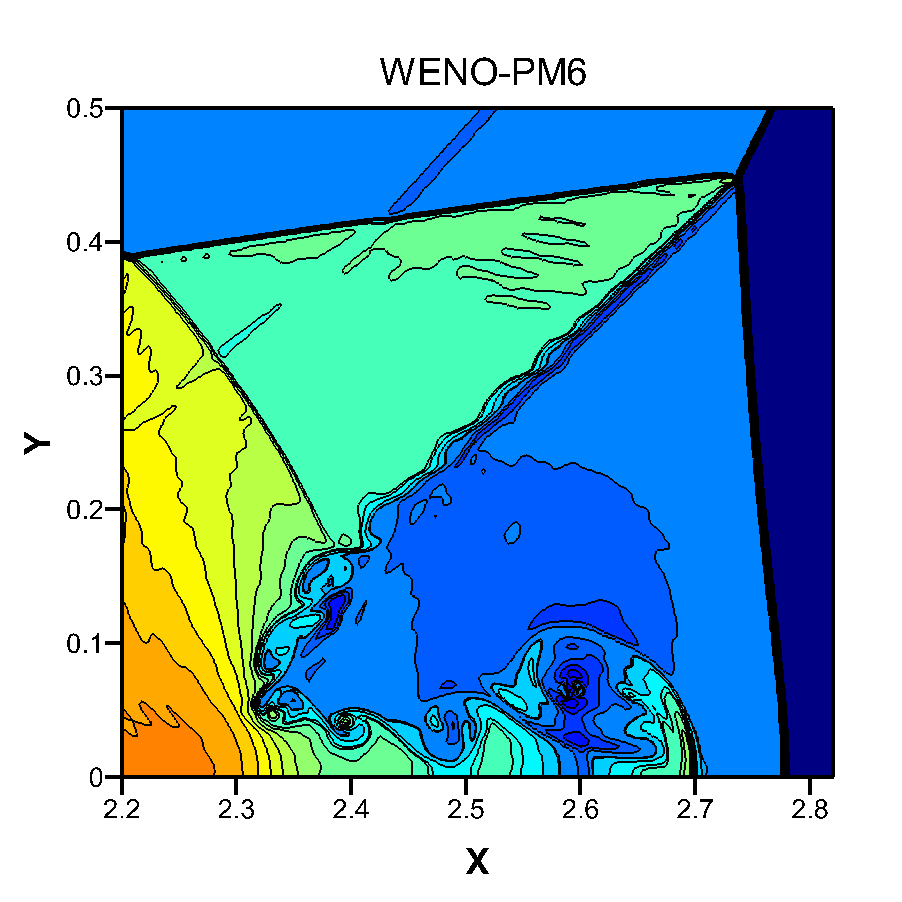}
  \includegraphics[height=0.32\textwidth]
  {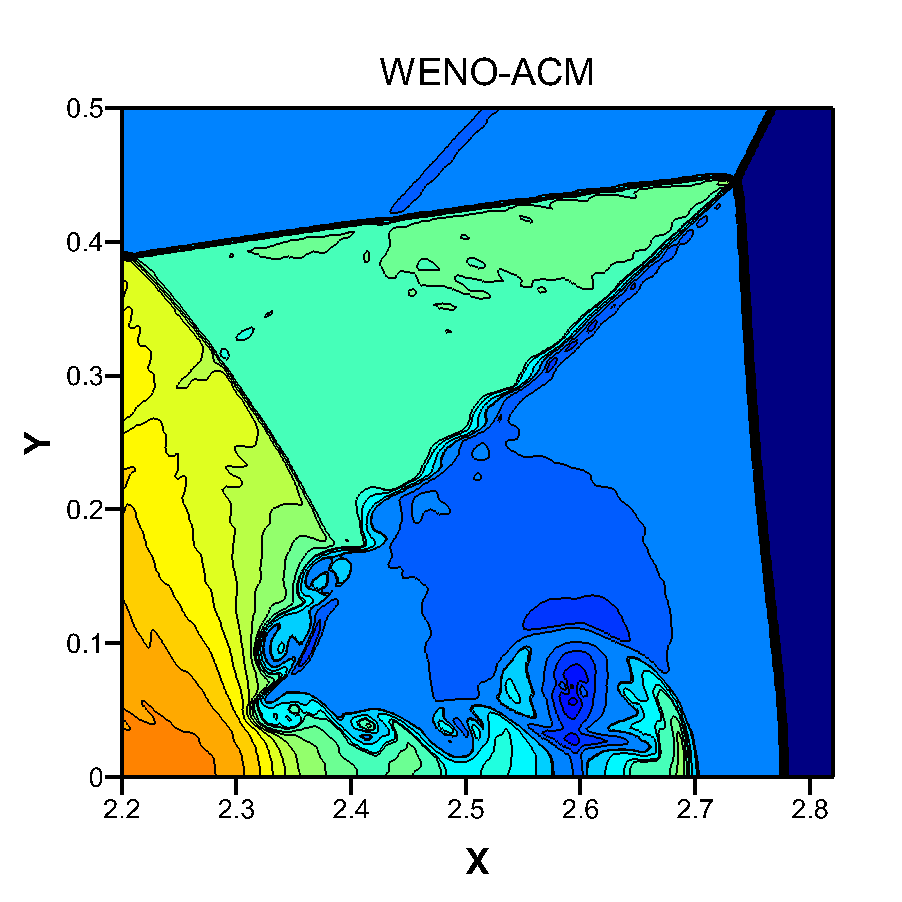}
\caption{The zoomed-in density plots for the DMR problem, computed
using the WENO-JS, WENO-M, WENO-PM6 and WENO-ACM schemes at output 
time $t = 0.2$ with a uniform mesh size of $2000\times500$.}
\label{fig:ex:DMR:ZoomedIn}
\end{figure}

In Fig. \ref{fig:ex:DMR}, we have shown the numerical results of
density obtained using the WENO-JS, WENO-M, WENO-PM6 and WENO-ACM 
schemes. Further, in Fig. \ref{fig:ex:DMR:ZoomedIn}, we have 
displayed the close-up view of the region around the double Mach
stems to observe more clearly the numerical solutions of all 
considered WENO schemes. In general, the global structure of the 
solution is very similar for different schemes, and all schemes are 
able to capture the companion structures behind the lower half of 
the right-moving reflection shock. However, the dissipation of the 
various schemes can be distinguished by the number and size of the 
small vortices generated along the slip lines. We can clearly see 
that the WENO-ACM and WENO-PM6 schemes capture more in number and 
bigger in size of the small vortices along the slip lines than the 
WENO-M and WENO-JS schemes, and it indicates that the resolving 
ability of the WENO-ACM and WENO-PM6 schemes is better than the 
other ones. In the solution of the WENO-JS scheme, the vortex rolled 
up near the slip lines is significantly dissipated due to the
excessive numerical dissipation. \minewB{It should be noted that 
all considered schemes suffered from the post shock oscillations in 
this test. The post shock oscillation seems slight in the solution 
of the WENO-JS schemes and it becomes more serious in the solution 
of the WENO-M scheme. Furthermore, it gets much more serious in the 
solutions of both the WENO-PM6 and WENO-ACM schemes. As mentioned in 
\cite{WENO-ZS}, these oscillations do not affect the ``ENO 
(essentially non-oscillatory)'' property of WENO schemes, but they 
are indeed responsible for the numerical residue to hang at the 
truncation error level of the scheme. As it is not the key point we 
are concerned about here, we refer to \cite{WENO-ZS} for more 
details, in which Zhang et al. proposed a systematic analysis and 
detailed discussion about this issue.}

\begin{example}
\bf{(Forward facing step problem, FFS)} 
\rm{This model problem was first presented by Woodward and Colella 
\cite{interactingBlastWaves-Woodward-Colella}. Recently, important 
details such as the physical instability and roll-up of the vortex 
sheet that emanates from the Mach stem have been successfully 
captured by various high order schemes 
\cite{ForwardFacingStep-Cockburn_Shu,ADER-WENO-1,ADER-WENO-2, 
WENO-eta,1DEuler_exact}. Our purpose is to prove that the WENO-ACM 
scheme is also able to successfully capture the roll-up of the 
vortex sheet and perform robustly on this severely stringent problem.

The setup of the problem is as follows. There is a step with a
height of $0.2$ length units located $0.6$ length units from the 
left-hand end of a wind tunnel, which is $1$ length unit wide and $3$
length units long. The computational domain of this problem is 
$\Omega = [0,0.6]\times[0,1]\cup[0.6,3]\times[0.2,1]$ and the flow 
is initialized by}
\label{ex:FFS}
\end{example}
\begin{equation*}
\big( \rho, u, v, p \big)(x, y, 0) = (1.4, 3.0, 0.0, 1.0), \quad
(x,y) \in \Omega.
\label{eq:initial_Euer2D:FFS}
\end{equation*}
Reflective boundary conditions are used along the walls of the wind
tunnel and the step, and inflow and outflow conditions are used at
the entrance and the exit of the wind tunnel respectively. In order
to handle the singularity at the left top corner of the step, the 
same technique used in \cite{interactingBlastWaves-Woodward-Colella} 
is employed. The computational domain is discretized with uniform 
mesh sizes of $900\times300$ and $1200\times400$, and the output 
time is chosen to be $t=4$. 

\begin{figure}[ht]
\centering
  \includegraphics[height=0.23\textwidth]
  {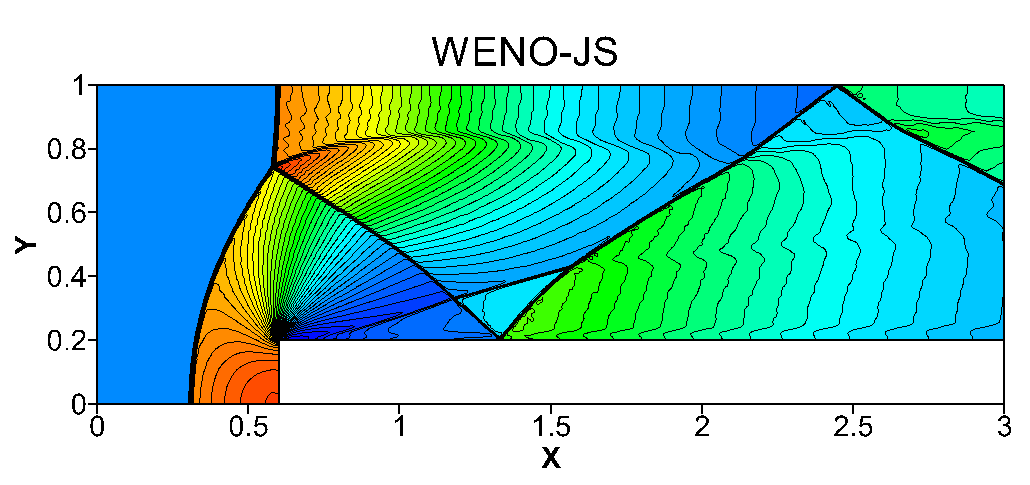}
  \includegraphics[height=0.23\textwidth]
  {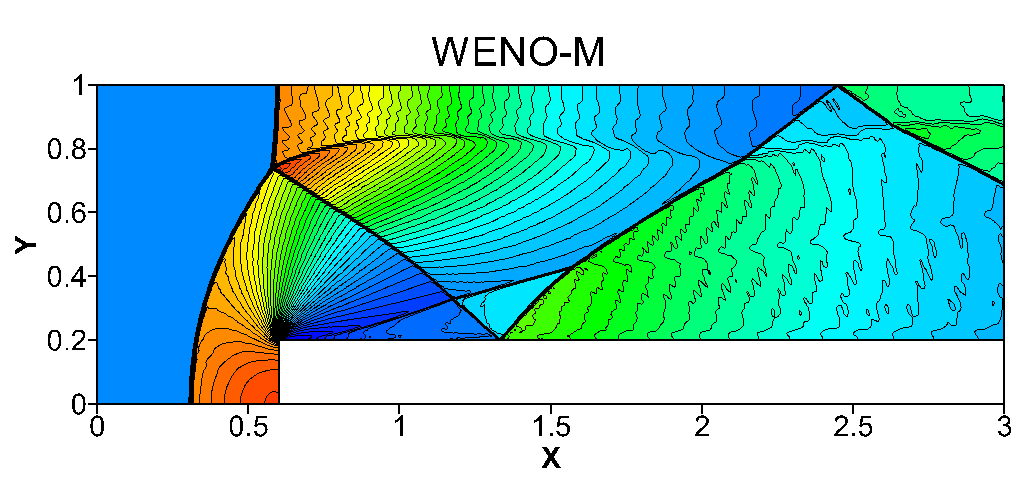}\\
  \includegraphics[height=0.23\textwidth]
  {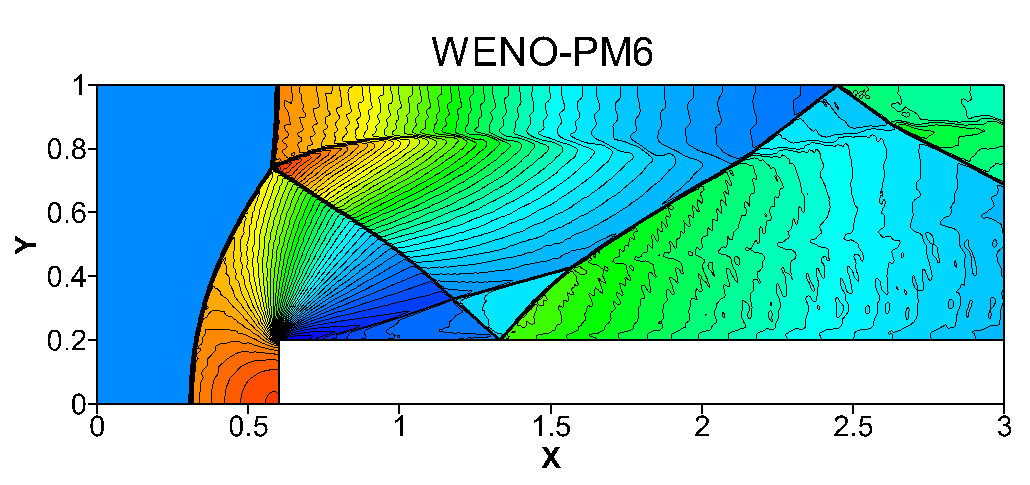}
  \includegraphics[height=0.23\textwidth]
  {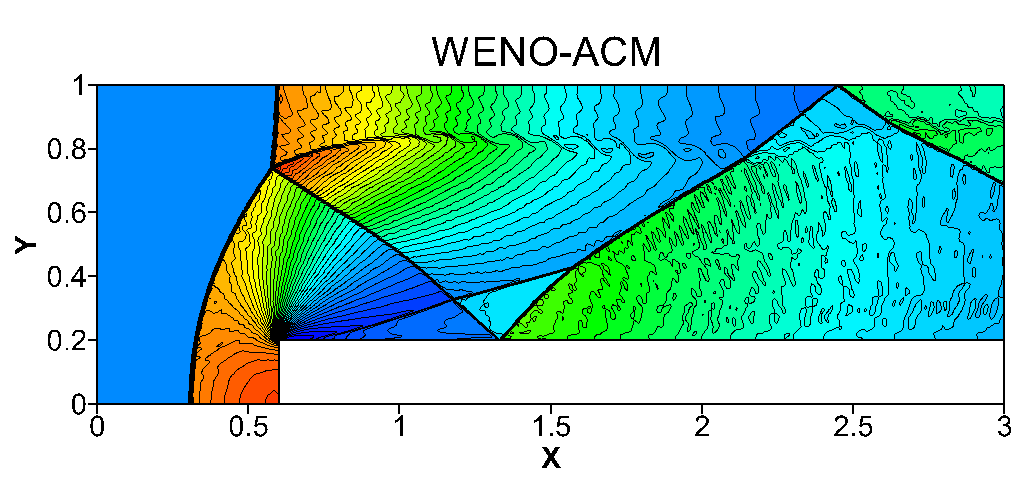}
\caption{Density plots for the forward facing step problem using 50
density contour lines with range from $0.32$ to $6.5$, computed
using the WENO-JS, WENO-M, WENO-PM6 and WENO-ACM schemes at output 
time $t = 4.0$ with a uniform mesh size of $900\times300$.}
\label{fig:ex:FFS:900x300}
\end{figure}

\begin{figure}[ht]
\centering
  \includegraphics[height=0.23\textwidth]
  {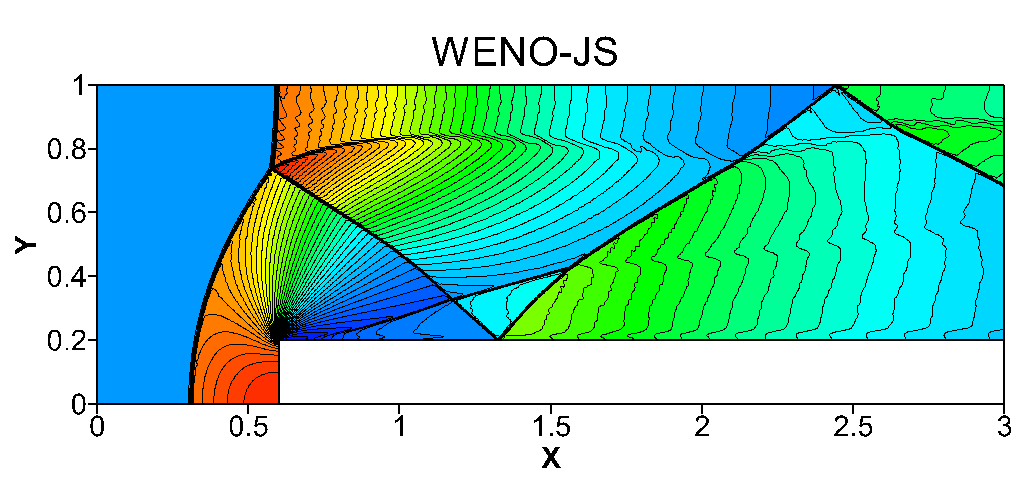}
  \includegraphics[height=0.23\textwidth]
  {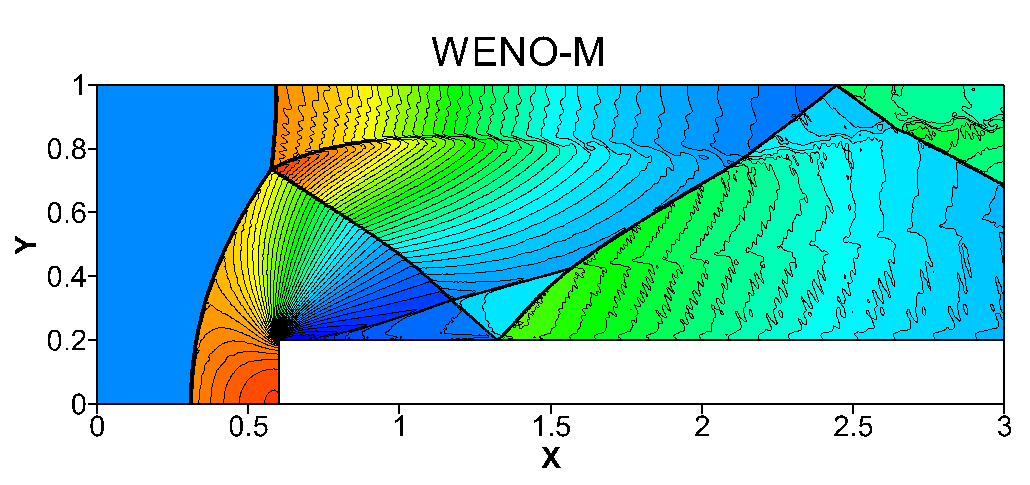}\\
  \includegraphics[height=0.23\textwidth]
  {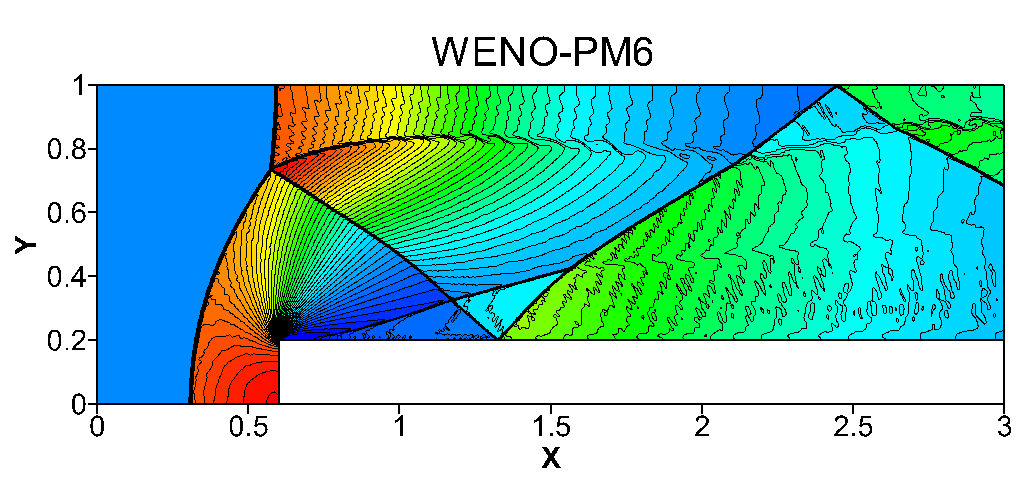}
  \includegraphics[height=0.23\textwidth]
  {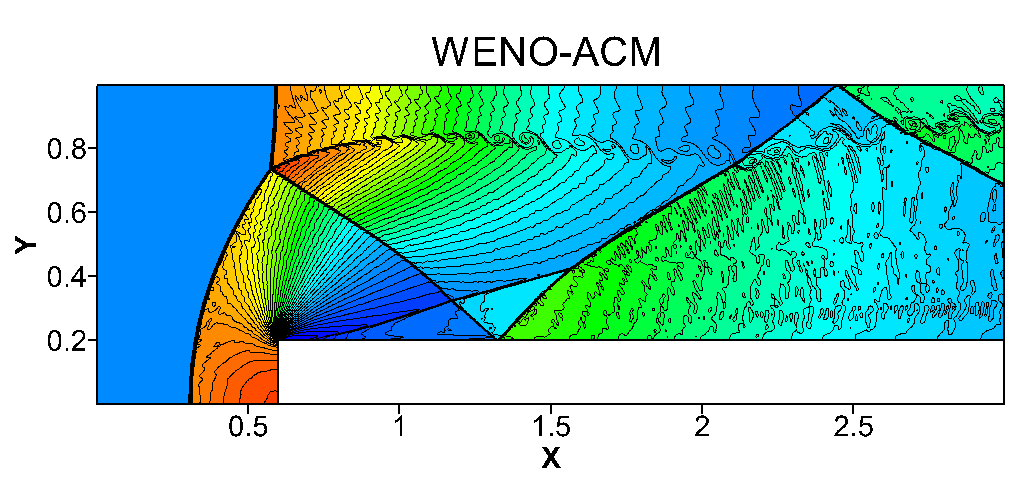}
\caption{Density plots for the forward facing step problem using 50
density contour lines with range from $0.32$ to $6.5$, computed
using the WENO-JS, WENO-M, WENO-PM6 and WENO-ACM schemes at output 
time $t = 4.0$ with a uniform mesh size of $1200\times400$.}
\label{fig:ex:FFS:1200x400}
\end{figure}

The density contours obtained by all considered schemes have been 
shown in Fig. \ref{fig:ex:FFS:900x300} and Fig. 
\ref{fig:ex:FFS:1200x400} for different computing mesh cells. We can 
observe that these considered schemes capture all the shocks 
properly with sharp profiles. From Fig. \ref{fig:ex:FFS:900x300}, we 
find that on the uniform mesh size of $900\times300$, the roll-up of 
the vortex sheet is clearly visible when the WENO-ACM scheme is 
used, while not observed in solutions of the other three considered 
WENO schemes. From Fig. \ref{fig:ex:FFS:1200x400}, we can see that 
with an increase of the mesh size from $900\times300$ to 
$1200\times400$, the roll-up of the vortex sheet becomes clearly 
visible when the WENO-PM6 scheme is used, while it is not 
particularly clear for the WENO-M scheme and still invisible for the 
WENO-JS scheme. Moreover, the roll-up of the vortex sheet is 
observed more evidently and maintains the most intensely when the 
WENO-ACM scheme is used. These demonstrate the advantage of the 
WENO-ACM scheme that has less dissipation and better resolution in 
capturing details of the complicated flow structures. \minewB{Again, 
we can see the post shock oscillations in solutions of all 
considered schemes as mentioned in the DMR example above. And very 
similar phenomenon is observed.}

\subsection{Computational cost comparison for 2D Euler problems}
\label{subsec:CostTest}
In this subsection, we compare the computational costs of the 
WENO-ACM scheme with the WENO-JS, WENO-M and WENO-PM6 schemes. The 
numerical experiments of the two-dimensional Euler system simulated 
earlier, that is, Example \ref{ex:shock-vortex} to Example 
\ref{ex:FFS} in subsection \ref{subsec:examples_2D_Euler}, are taken
as the test objects. These examples are tested with all conditions 
the same as in subsection \ref{subsec:examples_2D_Euler}. In 
addition, we also test each example with a different mesh size 
respectively. An in-house code developed in C++ is employed, running 
in serial mode under Windows system, and the CPU is Intel(R) 
Core(TM) i9-9880H. In order to rule out the effect of other 
operations such as the boundary treatment, initialization, etc., as 
applied in reference \cite{WENO-Z}, the CPU timing per Runge-Kutta 
step is considered here. And to mitigate the influence of 
randomness, each test is repeated three times under the same 
condition. 

Let $T(\mathrm{X})$ denote the \textit{cost} (the 
CPU time that a Runge-Kutta step consumes) of the WENO-X scheme and 
$P(\mathrm{X})$ denote the \textit{extra cost} brought by the 
mapping process of the WENO-X scheme compared to the WENO-JS scheme, 
that is, $P(\mathrm{X}) = \frac{T(\mathrm{X}) - T(\mathrm{JS})}{T(
\mathrm{JS})} \times 100\%$. In Table 
\ref{table:ComputationalCost:shock-vortex} to Table
\ref{table:ComputationalCost:FFS}, we have given $T(\mathrm{X})$ of 
all considered schemes and $P(\mathrm{X})$ (in brackets) of the 
WENO-M, WENO-PM6 and WENO-ACM schemes for Example
\ref{ex:shock-vortex} to Example \ref{ex:FFS}, respectively. In 
order to measure the \textit{reduced cost} of the mapping process of 
the WENO-ACM scheme compared to those of the WENO-M and WENO-PM6 
schemes, in the last two columns of each table, we have shown these 
\textit{reduced costs} in percetages computed by 
$\frac{P(\mathrm{M})- P(\mathrm{ACM})}{P(\mathrm{M})} = 
\frac{T(\mathrm{M})- T(\mathrm{ACM})}{T(\mathrm{M})-T(\mathrm{JS})}$ 
and $\frac{P(\mathrm{PM6})-P(\mathrm{ACM})}{P(\mathrm{PM6})} = 
\frac{T(\mathrm{PM6})-T(\mathrm{ACM})}{T(\mathrm{PM6})-
T(\mathrm{JS})}$ respectively.

From Table \ref{table:ComputationalCost:shock-vortex} to Table
\ref{table:ComputationalCost:FFS}, we can easily observe that: (1)
for each example with different mesh sizes, the 
\textit{cost} in the three tests has a certain degree of 
fluctuation; (2) on average, the 
\textit{extra costs} compared to the WENO-JS scheme 
are higher than $24\%$ for the WENO-M scheme, higher than $54\%$ for 
the WENO-PM6 scheme, while lower than $5\%$ for the WENO-ACM scheme; 
(3) the \textit{reduced costs} of the WENO-ACM 
scheme are more than $83\%$ and $93\%$ compared to the WENO-M scheme 
and the WENO-PM6 scheme, respectively.

\begin{table}[ht]
\begin{scriptsize}
\centering
\caption{CPU time (in seconds) and the extra computational cost compared to the
WENO-JS scheme (in percentage) per Runge-Kutta step of Example
\ref{ex:shock-vortex} as computed by considered WENO schemes.}
\label{table:ComputationalCost:shock-vortex}
\begin{tabular*}{\hsize}
{@{}@{\extracolsep{\fill}}llllllll@{}}
\hline
\space      & \space   & \space         & \space      & \space
& \space    & \multicolumn{2}{c}{Saved} \\
\cline{7-8}      
Grid size   & Test     & WENO-JS        & WENO-M      & WENO-PM6
& WENO-ACM  & Compared to WENO-M & Compared to WENO-PM6 \\
\hline
$200\times200$  &  Test 1 & 0.219  & 0.281(28.31\%) & 0.359(63.93\%) 
                          & 0.229(4.57\%) &  83.87\%  & 92.86\% \\
\space          &  Test 2 & 0.226  & 0.282(24.78\%) & 0.361(59.73\%) 
                          & 0.234(3.54\%) &  85.71\%  & 94.07\% \\
\space          &  Test 3 & 0.218  & 0.274(25.69\%) & 0.345(58.26\%) 
                          & 0.228(4.59\%) &  82.14\%  & 92.13\% \\
\space          & Average & 0.221  & 0.279(26.24\%) & 0.355(60.63\%) 
                          & 0.230(4.22\%) &  83.91\%  & 93.03\% \\
$400\times400$  &  Test 1 & 0.890  & 1.118(25.62\%) & 1.415(58.99\%) 
                          & 0.921(3.48\%) &  86.40\%  & 94.10\% \\
\space          &  Test 2 & 0.875  & 1.124(28.46\%) & 1.399(59.89\%) 
                          & 0.906(3.54\%) &  87.55\%  & 94.08\% \\
\space          &  Test 3 & 0.885  & 1.115(25.99\%) & 1.410(59.32\%) 
                          & 0.924(4.41\%) &  83.04\%  & 92.57\% \\
\space          & Average & 0.883  & 1.119(26.68\%) & 1.408(59.40\%) 
                          & 0.917(3.81\%) &  85.71\%  & 93.58\% \\
\hline
\end{tabular*}
\end{scriptsize}
\end{table}

\begin{table}[ht]
\begin{scriptsize}
\centering
\caption{CPU time (in seconds) and the extra computational cost compared to the
WENO-JS scheme (in percentage) per Runge-Kutta step of Example
\ref{ex:explosion} as computed by considered WENO schemes.}
\label{table:ComputationalCost:explosion}
\begin{tabular*}{\hsize}
{@{}@{\extracolsep{\fill}}llllllll@{}}
\hline
\space      & \space   & \space         & \space      & \space
& \space    & \multicolumn{2}{c}{Saved} \\
\cline{7-8}      
Grid size   & Test     & WENO-JS        & WENO-M      & WENO-PM6
& WENO-ACM  & Compared to WENO-M & Compared to WENO-PM6 \\
\hline
$200\times200$  &  Test 1 & 0.207  & 0.258(24.64\%) & 0.321(55.07\%) 
                          & 0.213(2.90\%) & 88.24\%   & 94.74\% \\
\space          &  Test 2 & 0.204  & 0.263(28.92\%) & 0.328(60.78\%) 
                          & 0.215(5.39\%) & 81.36\%   & 91.13\% \\
\space          &  Test 3 & 0.210  & 0.268(27.62\%) & 0.343(63.33\%) 
                          & 0.218(3.81\%) & 86.21\%   & 93.98\% \\
\space          & Average & 0.207  & 0.263(27.05\%) & 0.331(59.74\%) 
                          & 0.215(4.03\%) & 85.12\%   & 93.26\% \\
$400\times400$  &  Test 1 & 0.794  & 0.987(24.31\%) & 1.219(53.53\%) 
                          & 0.807(1.64\%) & 93.26\%   & 96.94\% \\
\space          &  Test 2 & 0.782  & 0.982(25.58\%) & 1.232(57.54\%) 
                          & 0.811(3.71\%) & 85.50\%   & 93.56\% \\
\space          &  Test 3 & 0.788  & 0.998(26.65\%) & 1.211(53.68\%) 
                          & 0.814(3.30\%) & 87.62\%   & 93.85\% \\
\space          & Average & 0.788  & 0.898(25.51\%) & 1.221(54.91\%) 
                          & 0.811(2.88\%) & 88.72\%   & 94.76\% \\
\hline
\end{tabular*}
\end{scriptsize}
\end{table}

\begin{table}[ht]
\begin{scriptsize}
\centering
\caption{CPU time (in seconds) and the extra computational cost compared to the
WENO-JS scheme (in percentage) per Runge-Kutta step of Example
\ref{ex:Riemann2D} as computed by considered WENO schemes.}
\label{table:ComputationalCost:Riemann2D}
\begin{tabular*}{\hsize}
{@{}@{\extracolsep{\fill}}llllllll@{}}
\hline
\space      & \space   & \space         & \space      & \space
& \space    & \multicolumn{2}{c}{Saved} \\
\cline{7-8}      
Grid size   & Test     & WENO-JS        & WENO-M      & WENO-PM6
& WENO-ACM  & Compared to WENO-M & Compared to WENO-PM6 \\
\hline
$600\times600$  &  Test 1 & 1.852  & 2.414(30.35\%) & 2.917(57.51\%) 
                          & 1.884(1.73\%) & 91.34\%   & 97.00\% \\
\space          &  Test 2 & 1.843  & 2.319(25.83\%) & 2.905(57.62\%) 
                          & 1.900(3.09\%) & 88.03\%   & 94.63\% \\
\space          &  Test 3 & 1.861  & 2.404(29.18\%) & 2.935(57.71\%) 
                          & 1.918(3.06\%) & 89.50\%   & 94.69\% \\
\space          & Average & 1.852  & 2.379(28.46\%) & 2.919(57.61\%) 
                          & 1.900(2.63\%) & 90.77\%   & 95.44\% \\
$1200\times1200$  
               &  Test 1 & 6.678  & 8.699(30.26\%) & 10.920(63.52\%) 
                         & 6.978(4.49\%) & 85.16\%   & 92.93\% \\
\space         &  Test 2 & 6.749  & 8.694(28.82\%) & 10.943(62.14\%) 
                         & 6.882(1.97\%) & 93.16\%   & 96.83\% \\
\space         &  Test 3 & 6.740  & 8.743(29.72\%) & 11.000(63.20\%) 
                         & 7.001(3.87\%) & 86.97\%   & 93.87\% \\
\space         & Average & 6.722  & 8.712(29.60\%) & 10.954(62.95\%)
                         & 6.954(3.44\%) & 88.37\%   & 94.53\% \\
\hline
\end{tabular*}
\end{scriptsize}
\end{table}

\begin{table}[ht]
\begin{scriptsize}
\centering
\caption{CPU time (in seconds) and the extra computational cost compared to the
WENO-JS scheme (in percentage) per Runge-Kutta step of Example
\ref{ex:DMR} as computed by considered WENO schemes.}
\label{table:ComputationalCost:DMR}
\begin{tabular*}{\hsize}
{@{}@{\extracolsep{\fill}}llllllll@{}}
\hline
\space      & \space   & \space         & \space      & \space
& \space    & \multicolumn{2}{c}{Saved} \\
\cline{7-8}      
Grid size   & Test     & WENO-JS        & WENO-M      & WENO-PM6
& WENO-ACM  & Compared to WENO-M & Compared to WENO-PM6 \\
\hline
$1000\times250$ &  Test 1 & 1.132  & 1.462(29.15\%) & 1.794(58.48\%) 
                          & 1.178(4.06\%) & 86.06\%   & 93.05\% \\
\space          &  Test 2 & 1.160  & 1.478(27.41\%) & 1.757(51.47\%) 
                          & 1.191(2.67\%) & 90.25\%   & 94.81\% \\
\space          &  Test 3 & 1.156  & 1.469(27.08\%) & 1.765(52.68\%) 
                          & 1.186(2.60\%) & 90.42\%   & 95.07\% \\
\space          & Average & 1.149  & 1.469(27.87\%) & 1.772(54.18\%) 
                          & 1.185(3.10\%) & 88.87\%   & 94.27\% \\
$2000\times500$ &  Test 1 & 4.299  & 5.610(30.50\%) & 6.749(56.99\%) 
                          & 4.462(3.79\%) & 87.57\%   & 93.35\% \\
\space          &  Test 2 & 4.278  & 5.618(31.32\%) & 6.739(57.53\%) 
                          & 4.435(3.67\%) & 88.28\%   & 93.62\% \\
\space          &  Test 3 & 4.307  & 5.648(31.14\%) & 6.766(57.09\%) 
                          & 4.502(4.53\%) & 85.46\%   & 92.07\% \\
\space          & Average & 4.295  & 5.625(30.98\%) & 6.751(57.20\%) 
                          & 4.466(4.00\%) & 87.10\%   & 93.01\% \\
\hline
\end{tabular*}
\end{scriptsize}
\end{table}

\begin{table}[ht]
\begin{scriptsize}
\centering
\caption{CPU time (in seconds) and the extra computational cost compared to the
WENO-JS scheme (in percentage) per Runge-Kutta step of Example
\ref{ex:FFS} as computed by considered WENO schemes.}
\label{table:ComputationalCost:FFS}
\begin{tabular*}{\hsize}
{@{}@{\extracolsep{\fill}}llllllll@{}}
\hline
\space      & \space   & \space         & \space      & \space
& \space    & \multicolumn{2}{c}{Saved} \\
\cline{7-8}      
Grid size   & Test     & WENO-JS        & WENO-M      & WENO-PM6
& WENO-ACM  & Compared to WENO-M & Compared to WENO-PM6 \\
\hline
$900\times300$  &  Test 1 & 1.873  & 2.497(33.32\%) & 3.142(67.75\%) 
                          & 1.945(3.84\%) & 88.46\%   & 94.33\% \\
\space          &  Test 2 & 1.880  & 2.502(33.09\%) & 3.124(66.17\%) 
                          & 1.926(2.45\%) & 92.60\%   & 96.30\% \\
\space          &  Test 3 & 1.868  & 2.465(31.96\%) & 3.147(68.47\%) 
                          & 1.928(3.21\%) & 89.95\%   & 95.31\% \\
\space          & Average & 1.874  & 2.488(32.79\%) & 3.138(67.46\%) 
                          & 1.933(3.17\%) & 90.34\%   & 95.31\% \\
$1200\times400$ &  Test 1 & 3.221  & 4.308(33.75\%) & 5.372(66.78\%) 
                          & 3.360(4.32\%) & 87.21\%   & 93.54\% \\
\space          &  Test 2 & 3.238  & 4.359(34.62\%) & 5.498(69.80\%) 
                          & 3.393(4.79\%) & 86.17\%   & 93.14\% \\
\space          &  Test 3 & 3.230  & 4.288(32.76\%) & 5.423(67.89\%) 
                          & 3.381(4.67\%) & 85.73\%   & 93.11\% \\
\space          & Average & 3.229  & 4.318(33.71\%) & 5.431(68.16\%) 
                          & 3.378(4.59\%) & 86.37\%   & 93.26\% \\                          
\hline
\end{tabular*}
\end{scriptsize}
\end{table}


\section{Conclusions}
\label{secConclusions} 
In order to reduce the computational cost introduced by mapping
processes of the mapped WENO-M \cite{WENO-M} and WENO-PM6
\cite{WENO-PM} schemes on the premise of retaining their advantages, 
we have devised a new mapped WENO scheme named WENO-ACM for 
hyperbolic conservation laws by introducing an approximate 
constant mapping function. It is theoretically and numerically 
demonstrated that the WENO-ACM scheme achieves the optimal 
convergence orders at critical points as the original mapped WENO-M 
scheme does. The new approximate constant mapping function satisfies 
the two additional properties first proposed in the WENO-PM6 scheme, 
and these properties ensure that the WENO-ACM scheme is able to 
generate comparable or better numerical solutions compared with the 
WENO-PM6 scheme, which has lower dissipation and provides higher 
resolution results than the classic WENO-JS and WENO-M schemes. 
Extensive numerical tests with two dimensional Euler equations show 
that the WENO-ACM scheme can reduce the cost of the mapping process 
by more than $83\%$ compared to the WENO-M scheme and by more than
$93\%$ compared to the WENO-PM6 scheme, making the extra 
computational cost reduced from more than $24\%$ for the WENO-M 
scheme and more than $54\%$ for the WENO-PM6 scheme to a more 
acceptable value of no more than $5\%$.



\bibliographystyle{model1b-shortjournal-num-names}
\bibliography{../refs_clear}

\end{document}